\documentclass[11pt, a4paper]{article}
\usepackage{times}
\usepackage{a4wide}
\usepackage[british]{babel}
\usepackage{enumerate, longtable}
\usepackage{amsmath, amscd, amsfonts, amsthm, amssymb, latexsym, comment, stmaryrd, graphicx}
\usepackage[all]{xy}
\usepackage[T1]{fontenc}
\usepackage[latin1]{inputenc}
\usepackage{booktabs}
\usepackage{cases}
\usepackage{csquotes}
\usepackage{dsfont}
\usepackage{enumitem}
\usepackage{float}
\usepackage[bottom]{footmisc}
\usepackage{mathabx}
\usepackage{mathrsfs}
\usepackage{mathtools}
\usepackage{pdfpages}
\usepackage{tikz-cd}
\usepackage[normalem]{ulem}
\usepackage{url}
\usepackage{xcolor}

\newtheorem{thm}{Theorem}[section]
\newtheorem{lem}[thm]{Lemma}
\newtheorem{defi}[thm]{Definition}
\newtheorem{rem}[thm]{Remark}

\newtheorem{prop}[thm]{Proposition}

\newtheorem*{thm*}{Theorem}
\newtheorem{definition}[thm]{Definition} 
\newtheorem{lemma}[thm]{Lemma}
\newtheorem{proposition}[thm]{Proposition}
\newtheorem{theorem}[thm]{Theorem}

\newtheorem{remark}[thm]{Remark}
\newtheorem*{theorem*}{Theorem}

\newcommand{\End}{{\rm End}}
\newcommand{\Hom}{{\rm Hom}}

\newcommand{\GL}{\mathrm{GL}}

\newcommand{\Sp}{\mathrm{Sp}}
\newcommand{\SL}{\mathrm{SL}}
\newcommand{\SO}{\mathrm{SO}}
\newcommand{\SU}{\mathrm{SU}}

\DeclareMathOperator{\Gal}{Gal}

\DeclareMathOperator{\Ker}{ker}

\newcommand{\Ind}{{\rm Ind}}

\newcommand{\PSL}{\mathrm{PSL}}

\newcommand{\N}{\mathbb{N}}
\newcommand{\Z}{\mathbb{Z}}
\newcommand{\Q}{\mathbb{Q}}

\newcommand{\Mod}[1]{\ (\mathrm{mod}\ #1)}

\makeatletter
\newtheorem*{rep@theorem}{\rep@title}
\newcommand{\newreptheorem}[2]{%
\newenvironment{rep#1}[1]{%
 \def\rep@title{#2 \ref{##1}}%
 \begin{rep@theorem}}%
 {\end{rep@theorem}}}
\makeatother
\newreptheorem{theorem}{Theorem}

\begin{document}

\selectlanguage{british}

\title{Automorphy of $\mathrm{GL}_2\otimes \mathrm{GL}_n$ in the self-dual case}
\author{
Sara Arias-de-Reyna\footnote{Departamento de \'Algebra, Facultad de Matem\'aticas, Universidad de Sevilla,
Avda.~Reina Mercedes s/n.
Apdo.~1160.~41080.~Sevilla,
Spain, sara$\_$arias@us.es Partially supported by project MTM2016-75027-P, funded by
the Ministerio de Econom\'ia y Competitividad of Spain and projects P 20$\_$01056 and US-1262169,
funded by the Fondo Europeo de Desarrollo Regional (FEDER) and the Consejer\'ia de Econom\'ia,
Conocimiento, Empresas y Universidad de la Junta de Andaluc\'ia.}, 
Luis Dieulefait\footnote{Departament d'\`Algebra i Geometria,
Facultat de Matem\`atiques,
Universitat de Barcelona,
Gran Via de les Corts Catalanes, 585,
08007 Barcelona, Spain, ldieulefait@ub.edu}, 
Josu P\'{e}rez\footnote{Departament d'\`Algebra i Geometria,
Facultat de Matem\`atiques,
Universitat de Barcelona,
Gran Via de les Corts Catalanes, 585,
08007 Barcelona, Spain, josu.perez.zarra@ub.edu. Amb el suport del programa d'ajuts predoctorals FI-SDUR del Departament de Recerca i Universitats de la Generalitat de Catalunya i el cofinan\c{c}ament pel Fons Social Europeus Plus. Expedient 2023FISDU00582.
   }}

\begin{figure}
	\centering
	\begin{minipage}{.5\textwidth}
  		\centering
  		\includegraphics[width=.5\linewidth]{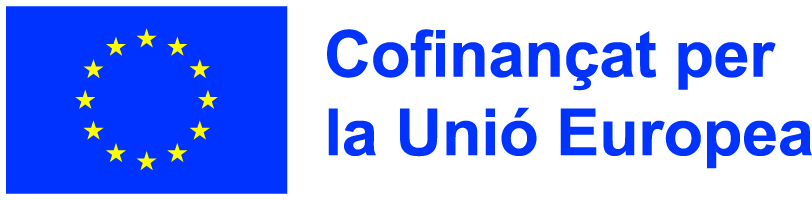}
	\end{minipage}%
	\begin{minipage}{.5\textwidth}
 		 \centering
		 \includegraphics[width=.5\linewidth]{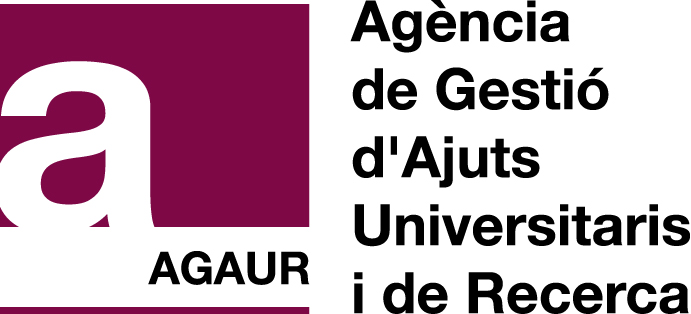}
	\end{minipage}
\end{figure}

\maketitle

\begin{quote}
\small \center \it A Pilar Bayer, en agradecimiento a m\'as de 3 d\'ecadas a cargo del STNB.
\end{quote}

\begin{abstract}
 In this paper we establish a new case of Langlands functoriality. 
 More precisely, we prove that the tensor product of the compatible system 
 of Galois representations attached to a level-1 classical modular form and 
 the compatible system attached to an $n$-dimensional RACP automorphic representation 
 of $\GL_n$ of the adeles of $\mathbb{Q}$ is automorphic, for any positive integer $n$, 
 under some natural hypotheses (namely regularity and irreducibility),   and a mild restriction on the level of the $n$-dimensional representation.

\end{abstract}

\section{Introduction}

The main goal in this paper is to prove a new case of Langlands functoriality, namely we will establish,
under some conditions, the automorphy of the tensor product of a classical modular form and an $n$-dimensional
automorphic representation, for every positive integer $n$. Previous results in this direction include
the case of $n=2$, solved by Ramakrishnan \cite{Ram}, the case of $n=3$ by Kim and Shahidi \cite{KiSh2002}, and
the automorphy of the $m$-fold tensor product of classical modular forms by Dieulefait \cite{Di-mfold}.

The case of functoriality we address stands as one of the most fundamental conjectures comprising the
so-called Langlands functoriality, together with other questions such as automorphy of symmetric powers or base change. A drawback of our result is that we require the classical modular form to be of level one. We are currently
 working on a generalisation for any classical modular form of odd level.

The main reference that we will follow through this paper is \cite{BLGGT}. In particular, we will work with
automorphic representations over totally real or CM fields that are RACP (regular algebraic cuspidal polarised) in
the sense of Section 2.1 of \cite{BLGGT} and the strongly compatible systems attached to them (see thm. 2.1.1 of \cite{BLGGT}).
Such Galois representations are RATOP (regular algebraic totally odd polarised) in the sense of Section 2.1 of loc. cit.
We say, as in Section 5 of loc. cit., that a compatible system of Galois representations is irreducible if there exists
a density one set $L$ of rational primes such that, for every prime in the field of coefficients of the system dividing a prime in $L$,
the corresponding Galois representation is irreducible.
We refer the reader to \cite{BLGGT} for all the relevant definitions. 

In order to define the tensor product of two compatible systems, we once and for all fix for each rational prime $\ell$ a prime $\lambda \subset \overline{\mathbb{Z}}$ above it (equivalently, we fix isomorphisms $\iota_\ell\colon\overline{\mathbb{Q}}_\ell\simeq\mathbb{C}$ for all primes $\ell$). This determines in each number field $F/\mathbb{Q}$ a choice of a prime $\lambda_F$ over $\ell$, in such a way that if $F\subset E$ then $\lambda_E\cap F=\lambda_F$.

Now, given two compatible systems of Galois representations of the absolute Galois group of a number field $F$, say $\{R_\lambda\}$ and $\{S_\lambda\}$ we denote by $\{(R\otimes S)_\lambda\}$ the compatible system indexed by the primes $\lambda$ of the compositum of the fields of definition of $R$ and $S$, obtained by $(R\otimes S)_{\lambda}=R_{\lambda}\otimes S_{\lambda}$. That this is well defined follows from our previous choice of primes.

The main tools that we will take from \cite{BLGGT} are the following:
two automorphy lifting theorems (in what follows, abbreviated to ALT), one for the ordinary case and the other for
the potentially diagonalisable case (Section 2 and 4 of loc. cit. respectively) and the result on changing the level and weight of an
automorphic representation (Section 4 of loc. cit.). This last result will be applied, in particular, to different kind of level raisings,
which in turn play a key role in order to ensure largeness of the residual images.

With that being said, the main result of this paper is:

\begin{thm}\label{thm:main} Let $f\in S_k(1)$ be a cuspidal eigenform and let $\{\rho_p(f)\}$ be the attached compatible system of Galois representations. Let $\pi$ be a RACP automorphic representation of $\mathrm{GL}_n(\mathbb{A}_{\mathbb{Q}})$  of level coprime to 3, where $n$
is a positive integer, and let $\{\rho_p(\pi)\}$ be the attached compatible system of Galois representations. We assume that
$\{\rho_p(f)\otimes r_p(\pi)\}$ is a regular and irreducible compatible system.
Then the system $\{\rho_p(f)\otimes r_p(\pi)\}$ is automorphic, i.e. there exists a RACP automorphic representation
$f \otimes \pi$ of $\mathrm{GL}_{2n}(\mathbb{A}_{\mathbb{Q}})$ such that the compatible system attached to $f \otimes \pi$ is isomorphic to  $\{\rho_p(f)\otimes r_p(\pi)\}$. 
\end{thm}

The notation $f \otimes \pi$ is justified, because due to known compatibilities of the compatible systems considered 
with the Local Langlands correspondence, the $L$-function of this automorphic form agrees with the product of 
the $L$-functions corresponding to $f$ and $\pi$.

\begin{rem}\label{rem:Irreducibilidad}
The irreducibility assumption for the tensor product is expected to hold in many cases. For example it is easy to see, using Ribet's large image result for classical modular forms, that if $n >3$ and the residual images of the system  $\{\rho_p(\pi)\}$ contain $\mathrm{Sp}_n(\mathbb{F}_{p})$
or $\Omega^{\pm}_n(\mathbb{F}_{p})$ for a density one set of primes then the compatible system $\{(\rho_p(f)\otimes r_p(\pi))\}$ is irreducible. This happens in particular in the case $n=4$, provided that $\pi$ is a lift of a genus $2$ Siegel modular form as in \cite{DZ}, thanks to Theorem 3.1 in \cite{DZ}. Moreover, it is expected that a similar result should hold for higher values of $n$, provided that the automorphic form is ``genuine", meaning that it does not come via lifting from a smaller reductive group, other than the symplectic or the orthogonal. Thus, as long as we restrict to ``genuine" self-dual automorphic forms, the irreducibility condition should always hold. \\
Just to give an example, if $\Delta$ denotes the cuspform of level $1$ and weight $12$ and $g$ denotes the genus $2$ Siegel cuspform of level $1$ and weight $20$, it is easy now to check all the required conditions to conclude that the theorem applies and the automorphic form $\Delta \otimes g$ is well-defined.
\end{rem}

\begin{rem}\label{rem:Prime3}
 In Theorem \ref{thm:main} we assumed that the level of the $n$-dimensional automorphic representation $\pi$ is coprime to 3. The reason for introducing this hypothesis is that, at a certain point in the proof, we will consider a modular form $g\in S_{2}(27)$ with complex multiplication by the quadratic field $\mathbb{Q}(\sqrt{-3})$, and we need to ensure that the image of the mod 13 Galois representation attached to $g$ is irreducible, after restriction to a certain field extension depending on $r_{\lambda}(\pi)$, $\lambda\vert 13$. 

We prove the main theorem under this assumption, but in principle the same techniques could be used to prove a similar result if the level of $\pi$ was not coprime to $3$. The strategy consists in considering a suitable infinite family of modular forms with complex multiplication, and replacing the modular form $g$ by another modular form in the family, chosen in such a way that the conductor of its field of complex multiplication is coprime to the level of $\pi$.
 
\end{rem}

\begin{rem}\label{rem:regularity} For some applications, such as the program proposed in \cite{ClozelThorneI}, it is enough to prove residual automorphy of a tensor product such as the one considered in the theorem, modulo a suitable prime $p$. It is worth remarking that in this situation, the regularity assumption can be removed: in fact, if it is not satisfied by the given pair $(f, \pi)$ one can replace $f$ by another cuspform $f'$ with the same mod $p$ Galois representation and sufficiently large weight (using, for example, the method of changing the weight in Section 3) so that the tensor product for the pair $(f', \pi)$ is regular. Then, if the rest of the conditions are satisfied, the theorem can be applied to establish automorphy of the second pair, thus residual automorphy (modulo $p$) of the given one.
\end{rem}

Together with the results from \cite{BLGGT} that we already mentioned, let us briefly list some tools
that we will use in our proof. We will make an essential use of the so-called Harris trick (see \cite{PRIMS} and \cite{BLGGT})
both to reduce to the case of even dimension $n$ and to deduce automorphy of the tensor product in the case of a CM classical
modular form. We will also use the safe chain of congruences linking a classical level one modular forms to a CM one, obtained in \cite{Di-Symm5} 
(including the Corrigenda to this paper, where some Steps of the chain are changed) and \cite{Di-mfold},
and also some variants of the ALT in \cite{Di-mfold}, which are specifically designed for tensor products. While the proof of the ALT in \cite{Di-mfold} are only sketched, in this article we expand on the details (see Theorems \ref{thm:mixto} and \ref{thm:pd}). These theorems rely on a generalization of a theorem of Thorne (Theorem 7.1 in \cite{Th2012}) to the semistable case, generalization of which we also present a proof in Section 3.
In the process of
changing the weight of the automorphic representation, a result on existence of weight zero lifts (which is the main result of
\cite{GHLS}) will be essential. Furthermore, we need several results on adequacy to be applied to the residual images, including the
case $p\vert n$, which are proved in \cite{GHT2016}. These results enable a generalisation of some (not all) ALT to the case $p\vert n$.
Last but not least, we will need two results that allow us to guarantee large
images of the residual representations in dimension $n$ under suitable ramification conditions in the compatible system, one of them
is proved in \cite{KLS08} and the other one is a result of Guralnick, whose proof we include due to the lack of a suitable reference.

The general idea of the proof of Theorem \ref{thm:main} is to build a chain of congruences connecting the given tensor product of
compatible systems to another tensor product such that the $2$-dimensional component is CM. For this last case, automorphy of the tensor product follows
(as in Harris trick) by applications of base change. We will call this tensor product the seed, because automorphy of any of the tensor products
 considered is eventually reduced to the automorphy of this one. In each of these congruences, either one or both of the two components in the
tensor product is changed in order to manipulate its level, weight, or both. The main technical difficulty lies on the fact that at
each of these congruences we want to apply a suitable ALT ``from right to left", in other words, we want to propagate automorphy
from the last element in the chain to the first one to obtain our goal. This requires, in particular, keeping track of local conditions at both sides
of the congruence, making sure that all tensor products considered are regular, and ensuring that the residual images are sufficiently
large.
We follow the $2$-dimensional safe chain referred above
whenever possible, but at some steps we will need to intersperse moves in the $n$-dimensional component in order that some ALT
can be applied.

We now proceed to give a brief description of the organization of the paper. In Section 2 we collect the results from group
theory that we will need. Most of what is contained in the section, including Subsection 2.2, has been proved by Guralnick. In Section 3 we collect the ALT that
we will use later on, and some results on changing the weights and level. Some of these ALT have been studied in other articles, while others are new theorems of this paper (the so called, ``Mixed'', ``Harris'', and ``semistable'' cases, of which a version already appeared in \cite{Di-mfold}). In Section 4 we recall step by step the safe chain from  \cite{Di-Symm5} and its Corrigenda and \cite{Di-mfold}, pointing out the relevant information for our construction. In Section 5 we begin the proof of Theorem \ref{thm:main} with three very technical preliminary steps, in which the
$n$-dimensional component is manipulated. Observe that from the beginning we are required to base change to a suitable solvable CM field $F$
and the automorphic $n$-dimensional representations that will appear in the safe chain from this point on are only defined over $F$.
This is harmless for our purposes because of solvable base change (which is used in both directions). Finally, Section 6 explains how
the series of congruences that safely link to the seed is constructed. \\

Acknowledgements: We want to thank Robert Guralnick and Toby Gee for contributing to this project in essential ways. Among other things, R.G. proved most of what is included in section 2, and also some of his recent papers on adequacy written with several collaborators arose in part from our conversations on what was needed for this paper. Similarly, the paper by T.G. and collaborators on the existence of weight $0$ lifts arose from his interest in this project. We also thank Florian Herzig for a careful reading and several comments on some parts of the paper. We deeply thank the anonymous referee for many suggestions on where to add explanations and how to improve the presentation of the results in this work.

\section{Some Group Theory}

In this section we collect some tools from group theory that will be applied to ensure that the residual image of certain representations
is large, allowing us to apply ALT and weight and level change. In the first subsection we include a theorem, communicated to us by Guralnick,
and we recall the notion of $(n, p)$-group of Khare, Larsen and Savin and the Main Theorem of section 2 of their paper \cite{KLS08}. Next, in Subsection \ref{subsec:AdequatePlus} we introduce
the notion of adequate$+$ subgroups, and include several properties. The contents of this subsection are due to Guralnick.

\subsection{Subgroups containing an almost simple group}\label{sec:subgroups}

The following theorem and its proof have been communicated to us by R. Guralnick.

\begin{thm}[Guralnick]\label{thm:Guralnick} Let $k$ be an algebraically closed field of characteristic $p$, and let $G$ be a finite (or Zariski closed) subgroup
of $\mathrm{GL}_n(k)$. Assume that $p>n+2$ and assume that $G$ contains no non-trivial normal subgroup consisting of unipotent elements.
If $G$ contains unipotent elements of every possible Jordan form, then $G$ contains $\mathrm{SL}_n(\mathbb{F}_p)$ or $\mathrm{SU}_n(\mathbb{F}_p)$.
\end{thm}

\begin{proof} We give the proof for the case that $G$ is finite (if $G$ is positive dimensional, the same proof applies). Call $V$ the space on which $\mathrm{GL}_n(k)$ acts, and
let $H\subset G$ be the subgroup of $G$ generated by its unipotent elements. We assume without loss of generality that $G=H$. The assumptions imply that $G$ acts completely reducibly on $V$, by \cite[Theorem A]{Gu}. The fact that $G$ contains a unipotent element with a
single Jordan block of size $n$ implies that $G$ is irreducible.

It also follows from \cite[Theorem B]{Gu} that either $G$ is a central product of quasi-simple finite groups of Lie type or $G$ is the sporadic simple group $J_1$ and $p=11$.
Since the Sylow $11$-subgroup
of $J_1$ is cyclic, only one possible Jordan form is represented in $G$. So we may assume that $G$ is a central product of quasi-simple finite groups of Lie type. Let $X$
be the corresponding product of simple algebraic groups; thus $X$ is connected and semisimple.

By \cite[Theorem B]{SS} (see also the proof of \cite[Theorem C]{SS}), it follows that $X$ contains $\mathrm{SL}_n(k)$. In particular, this implies that $G$ is simple
and must be a form of $\mathrm{SL}_n$ over a finite field, whence
$G=\mathrm{SL}_n(\mathbb{F}_q)$ or $\mathrm{SU}_n(\mathbb{F}_q)$, where $q$ is some power of $p$. The result follows.
\end{proof}

\begin{rem}
If $n>10$, one could use \cite{GS} to remove the condition $p>n+2$.
\end{rem}

Now we recall a group theoretical result from \cite{KLS08}. Before stating the theorem, we recall the following definitions:
\begin{itemize}

\item Let $\Gamma$ be a group, and $d$ a positive integer. We define $\Gamma^d$ as the intersection of all normal subgroups of $\Gamma$ of index at most $d$.

\item Let $n\geq 2$ be an integer and $p\equiv 1 \pmod{n}$ be a prime. A group of type $(n, p)$ is defined as a non-abelian homomorphic image of an extension
of $\mathbb{Z}/(n)$ by $\mathbb{Z}/(2p)$ such that $\mathbb{Z}/(n)$ acts faithfully on $\mathbb{Z}/(p)$.

\end{itemize}

\begin{thm}[Khare, Larsen, Savin]\label{thm:KLS}
Let $n\geq 2$ be an integer. There exist constants $d(n)$ and $p(n)$ depending only on $n$ such that,
if $d>d(n)$ is an integer, $p>p(n)$ and $\ell$ are distinct primes,
and $\Gamma\subset \mathrm{GL}_n(\mathbb{F}_{\ell^w})$ is a finite group such that $\Gamma^d$
contains a group of type $(n, p)$, then there exist $\alpha\geq w$, $g\in\mathrm{GL}_n(\mathbb{F}_{\ell^\alpha})$
and $1\leq k\leq \alpha$ such that $g^{-1}\Gamma g$ is  of the following form: a group containing
$\mathrm{SL}_n(\mathbb{F}_{\ell^k})$, $\mathrm{SU}_n(\mathbb{F}_{\ell^k})$,  $\mathrm{Sp}_n(\mathbb{F}_{\ell^k})$
or $\Omega^{\pm}_n(\mathbb{F}_{\ell^k})$ and contained in its normaliser in $\mathrm{GL}_n(\mathbb{F}_{\ell^\alpha})$.

Moreover, in each of the four cases, the quasi-simple group contained in $g^{-1}\Gamma g$ contains an element of order $p$.
\end{thm}

\begin{rem}
See Theorem 2.2 of \cite{KLS08}. The last statement is not included in the theorem, but it follows from its proof.
\end{rem}

\subsection{Adequate$+$ subgroups}\label{subsec:AdequatePlus}

J. Thorne introduced in \cite{Th2012} the notion of adequacy for a subgroup of $\mathrm{GL}_n(k)$, where $k$ is a finite field of odd
characteristic $p$. His definition implied that $p\nmid n$. Later this notion was extended in \cite{Thorne2016} to cover the case $p\mid n$.
In this paper, whenever we consider an adequate group, we will be referring to this extended notion. Note that if $p$ does not divide $n$,
then the two notions coincide.

We now introduce a stronger notion, which was proposed to us by Guralnick. The reason for this new notion is that it behaves well with respect to tensor products.

\begin{defi}\label{def:weakadequacy}  Let $k\subset \overline{\mathbb{F}}_p$ be a finite field, $V$ a finite dimensional $k$-vector space. We say that
$G\subset \GL(V)$ is \emph{adequate$+$} if it satisfies the following
conditions:
\begin{enumerate}
 \item $\mathrm{Ext}^1_G(V,V)=0$,
 \item $H^1(G,k)=0$,
 \item $H^2(G,k)=0$,
 \item (Weakly adequate) $\mathrm{End}(V)$ is spanned by the semisimple elements $g\in G$.
\end{enumerate}
\end{defi}

\begin{lem}\label{lem:condH2}
Let $G\subset \GL_n(k)$. Then the following are equivalent:

\begin{enumerate}
\item $G$ is adequate and $H^2(G,k)=0$;

\item $G$ is adequate$+$.
\end{enumerate}
\end{lem}

\begin{proof}
Both implications are proved in \cite{GHT2016} in Section 1. For example, condition (ii) in the definition of adequate follows from (1) and (3) above.
\end{proof}

For future use, we record the following result:

\begin{lem}\label{lem:SL2}
 Let $\ell\geq 7$ be a prime number. Then $\mathrm{SL}_2(\mathbb{F}_{\ell})$ in its natural $2$-dimensional representation is adequate$+$. The same holds for
$\mathrm{SL}_2(\mathbb{F}_{3^r})$ for any $r>2$.
\end{lem}

\begin{proof} Adequacy follows from Corollary 9.4 of \cite{GHT2016}. Moreover the $H^2$ condition is satisfied (see e.g.~\cite[Theorem 4.1]{BNP}).
\end{proof}

\begin{thm}\label{thm:adequate}
Let $G\subset \mathrm{GL}_n(\overline{\mathbb{F}}_{\ell})$ be one of the four groups $\mathrm{SL}_n(\mathbb{F}_{\ell^k})$,
$\mathrm{SU}_n(\mathbb{F}_{\ell^k})$,  $\mathrm{Sp}_n(\mathbb{F}_{\ell^k})$
or $\Omega^{\pm}_n(\mathbb{F}_{\ell^k})$. If  $k$ is sufficiently large, then $G$ is adequate$+$.
\end{thm}

\begin{proof}
The fact that the group $G$ is adequate is a particular case of \cite[Theorem 11.5]{GHT2016}.
By Lemma \ref{lem:condH2}, it suffices to prove that $H^2(G, \mathbb{F}_{\ell^k})=0$ with $k$
sufficiently large, which is well-known (see e.g.~\cite[Theorem 4]{Hiller}).
\end{proof}

\begin{prop}\label{prop:DirectProduct} Let  $k\subset \overline{\mathbb{F}}_p$ be a finite field,
$V_i$  finite dimensional $k$-vector spaces, and for each $i=1, \dots, r$, let $G_i\subset \GL(V_i)$
be adequate$+$ subgroups. Then the image of the direct product $G_1\times \cdots \times G_r$ in $\GL(V_1\otimes \cdots \otimes V_r)$ is adequate$+$.
\end{prop}

\begin{proof} First of all, note that  the four conditions in Definition \ref{def:weakadequacy} for a group $G$ to be adequate$+$ are preserved
when we replace $G$ by a quotient $G/N$, whenever $N$ is a finite, normal subgroup of $G$ of order coprime to $p$.
In our situation, the natural map $\GL(V_1)\times \cdots \times \GL(V_r)\rightarrow \GL(V_1\otimes \cdots \otimes V_r)$
has a finite kernel of order coprime to the characteristic of $k$, thus it suffices to prove that the 
direct product $G_1\times \cdots \times G_r$ acting on $V_1\otimes \cdots \otimes V_r$ 
satisfies the four conditions in Definition \ref{def:weakadequacy}.

Condition (4) in Definition \ref{def:weakadequacy} is preserved by direct product, since for any semisimple elements $g, h$, the tensor product $g\otimes h$ is semisimple.
Set $G=G_1\times \cdots \times G_r$ and $V=V_1\otimes \cdots \otimes V_r$. By K\"unneth formula, 
$\mathrm{Ext}^1_G(V, V)=H^1(G, V^{*}\otimes V)=
\bigoplus_i H^1(G_i, V_i^{*}\otimes V_i)\otimes 
H^0(\prod_{j\not=i} G_j,  \bigotimes_{j\not=i}V_j^{*}\otimes V_j)=0$.
Also by K\"unneth formula, $H^j(G, k)=0$ for  $j=1, 2$ follows from $H^j(G_i, k)=0$ for  $j=1, 2$.
\end{proof}

In our situation, we will have two (or more) representations $\rho_i:G_{F}\rightarrow \GL_n(\overline{\mathbb{F}}_p)$, which are not linearly disjoint, but may share some small quotient.
To deal with this situation, we will utilise the following elementary lemma.

\begin{lem}\label{lem:PrimeToP} Let $H\lhd G\subset \GL(V)$ be subgroups. Assume that $H\subset \GL(V)$ is adequate (resp. adequate$+$)
and $[G:H]$ has order prime to $p$. Then $G\subset \GL(V)$ is adequate (resp. adequate$+$).
\end{lem}

\begin{proof}
Since $H\subset G$, the fact that Condition (4) in the definition of adequate$+$
is satisfied for $H$ implies that it is also satisfied for $G$.
The rest follows from the fact that the restriction maps in cohomology are injective from the coprimeness assumption.
\end{proof}

A similar reasoning as in the proof of Lemma \ref{lem:PrimeToP} yields the following result.

\begin{lem}\label{lem:ForDihedral} Let  $k\subset \overline{\mathbb{F}}_p$ be a finite field, $V$ a finite dimensional $k$-vector space,
and $H\subset \mathrm{GL}(V)$ be a group of order prime to $p$ such that $H\subset \mathrm{GL}(V)$ is absolutely irreducible. 
Then $H\subset \mathrm{GL}(V)$ is adequate$+$.
\end{lem}

In order to apply Theorem \ref{thm:KLS}, we will need the following result:

\begin{prop}\label{prop:adequate} $n\geq 2$, $w\geq k\geq 1$ integers, $\ell$ a prime number (which we assume greater than
$3$ when $n=8$), and
let $G\subset \GL_n(\mathbb{F}_{\ell^w})$ be a group containing $\mathrm{SL}_n(\mathbb{F}_{\ell^k})$, $\mathrm{SU}_n(\mathbb{F}_{\ell^k})$,  $\mathrm{Sp}_n(\mathbb{F}_{\ell^k})$
or $\Omega^{\pm}_n(\mathbb{F}_{\ell^k})$ and contained in its normaliser in $\mathrm{GL}_n(\mathbb{F}_{\ell^w})$.
Then the index of the corresponding quasi-simple classical group in $G$ is coprime to $\ell$. Moreover, if  $k$ is sufficiently large, then $G$ is adequate$+$.
\end{prop}

\begin{proof}
We begin by proving the first claim. The case when $G$ contains $\mathrm{Sp}_n(\mathbb{F}_{\ell^k})$ is proved in \cite[Corollary 2.6]{KLS08}, and
we will give a proof that covers all the cases following a similar reasoning. Since the outer automorphism groups of classical groups are
well known (cf.~\cite{ATLAS}), it suffices to prove that field automorphisms do not belong to the normaliser of the corresponding classical group (note
that we are excluding the only case where graph automorphisms can occur, namely the case of characteristic $\ell=3$ for type $D_4$).
Since $G$ acts absolutely irreducibly, it is enough to show that the representation is not invariant under any field automorphism.
This follows from the fact that in all four cases, the classical group contains an element whose trace does not live in any proper subfield of $\mathbb{F}_{\ell^k}$,
which can be seen by restricting to a small subgroup such as $\SL_2$ for the first three cases and $\SO_3$ for the orthogonal case.

The second claim follows automatically from Theorem \ref{thm:adequate} and Lemma \ref{lem:PrimeToP}.
\end{proof}

\begin{rem}\label{rem:TensorProduct} Let $F$ be a field, and
assume we have two representations $\rho_1, \rho_2$ of the absolute Galois group $G_F$ into finite-dimensional $k$-vector spaces
$V_1$ and $V_2$ respectively, where $k\subset \overline{\mathbb{F}}_p$ is a finite field. Let $K_i$ denote the field cut out by 
$\rho_i$, for $i=1, 2$. Let further $K\subset K_1\cap K_2$ be a finite Galois extension of $F$ with
degree coprime to $p$, such that the extension $K_1/K$ and $K_2/K$ are linearly disjoint, and such that the image of $\rho_i(G_{K})$
is adequate+. Then combining 
Proposition \ref{prop:DirectProduct} and Lemma \ref{lem:PrimeToP}, we can conclude that the image of the tensor product of
 $\rho_1$ and $\rho_2$ is also adequate$+$. This fact will be used repeatedly in the paper.

We expand on a concrete situation that will arise in Section \ref{sec:Diagrams}, Remark \ref{rem:AdequateImage}.
We have a tensor product of a $2$-dimensional representation $\overline{\rho}:G_F\rightarrow \GL_2(k)$ 
and a $2n$-dimensional representation 
$\overline{r}:G_F\rightarrow \GL_{2n}(k)$, for some $n>1$ different from $4$ and some finite field $k\subset \overline{\mathbb{F}}_{p}$, in such a way that
the image of the restriction of $\overline{\rho}$ and $\overline{r}$ to $G_{F(\zeta_{p})}$ is adequate$+$. More precisely, the image
of $\overline{\rho}\vert_{G_{F(\zeta_p)}}$ contains a subgroup of the form $\SL_2(\mathbb{F}_{p^r})$, and the image of $\overline{r}\vert_{G_{F(\zeta_p)}}$ contains a classical group
of the form $\SL_{2n}(\mathbb{F}_{p^s})$, $\SU_{2n}(\mathbb{F}_{p^s})$, $\Sp_{2n}(\mathbb{F}_{p^s})$ or $\Omega^{\pm}_{2n}(\mathbb{F}_{p^s})$, 
and is contained in its normaliser in $\GL_{2n}(\mathbb{F}_{p^w})$ for some (big enough) exponents $r, s, w$.
Now, if we restrict to the extension $M$ of $F(\zeta_p)$ cut out by the determinant of $\overline{\rho}$, we can assume that the image of 
$\overline{\rho}$ coincides with the almost simple group $\SL_2(\mathbb{F}_{p^r})$, which must be linearly disjoint with the extension of 
$M$
cut out by $\overline{r}$. Both $\overline{\rho}(G_M)$ and $\overline{r}(G_M)$ are adequate$+$ (cf.~Lemma \ref{lem:SL2} and 
Proposition \ref{prop:adequate}), thus 
we are able to apply 
Proposition \ref{prop:DirectProduct} and Lemma \ref{lem:PrimeToP} to conclude that $(\overline{\rho}\otimes \overline{r})(G_{F(\zeta_p)})$ is adequate$+$.

\end{rem}


\section{Automorphy Lifting Theorems and change of level and weight}\label{sec:ALT}

As explained in the introduction, the strategy to the proof of our main theorem is to construct chains of congruences between tensor products of modular forms and automorphic forms. In order to propagate automorphy from one end of the chain to the other, it is crucial that each of the congruences making up the chain preserves automorphy. To ensure that this is the case, we will make use of several automorphy lifting theorems (ALT for short).

In the first part of this section we gather the statements and proofs of the ALT that we will use in future sections. These are mostly taken from \cite{BLGGT}, or inspired from it. The rest of the section is devoted to the proof of a minimal automorphy lifting theorem, an extension of a theorem of Thorne \cite{Th2012} which plays an important role in the proofs of Theorems \ref{thm:mixto} and \ref{thm:pd}.

\subsection{Statements of main ALT}
\label{sec:statementsALT}

We will apply several automorphy lifiting theorems in the paper. To facilitate
the reading of the paper, we will label them  as we state them.

 The first couple of results are borrowed from \cite{BLGGT}. We recall their statements here for the convenience of the reader.

\begin{thm*}[ALT-ORD]
 Suppose that $F$ is a CM (or totally real) field, that $\ell$
is an odd prime and that $n \in \mathbb{Z}_{\geq 1}$. Let $(r, \mu)$ be an $n$-dimensional, algebraic,
polarized $\ell$-adic representation of $G_F$ satisfying the following properties:
\begin{enumerate}
\item the reduction $\overline{r}$ is irreducible and $\overline{r}(G_{F(\zeta_l)}) \subset GL_n(\overline{\mathbb{F}}_{\ell})$ is adequate;
\item $\zeta_{\ell}\not\in F$; 
\item $r$ is ordinary at all primes above $\ell$;
\item $(\overline{r}, \overline{\mu})$ is ordinarily automorphic.
\end{enumerate}
Then $(r, \mu)$ is ordinarily automorphic. If $r$ is also crystalline (resp.~potentially
crystalline), then $(r, \mu)$ is ordinarily automorphic of level prime to $\ell$ (resp.~potentially level prime to $\ell$).\end{thm*}

This result is precisely Theorem 2.4.1 of \cite{BLGGT}, for an odd prime $\ell$, in the case when $\ell$ does not divide the dimension $n$. Corollary 7.3 of \cite{Thorne2016}
ensures that this result remains valid in the case when $\ell$ divides
the dimension $n$, with the more general definition of adequacy in \cite{Thorne2016}.

The next result is Theorem 4.2.1 of \cite{BLGGT} for the potentially diagonalisable situation.

\begin{thm*}[ALT-PD] Let $F$ be an imaginary CM field with maximal totally real subfield $F^{+}$ and let $c$ denote the nontrivial element of $\mathrm{Gal}(\overline{F}/F^{+})$. Suppose that that $\ell$ is an odd prime, and that $(r, \mu)$ is a regular algebraic, irreducible, $n$-dimensional, polarized representation of $G_F$. Let $\overline{r}$ denote the semi-simplification of the reduction of $r$ and let $d$ denote the maximal dimension of an irreducible subrepresentation of the restriction of $\overline{r}$ to the closed subgroup of $G_F$ generated by all Sylow pro-$\ell$-subgroups. Suppose that $(r, \mu)$ enjoys the following properties:
\begin{enumerate} 
 \item $r\vert_{G_{F_v}}$ is potentially diagonalizable (and so in particular potentially crystalline) for all $v\vert \ell$.
 \item The restriction $\overline{r}\vert_{G_{F(\zeta_{\ell})}}$ is irreducible, $\ell \geq 2(d + 1)$, and $\zeta_{\ell}\not\in F$.
 \item  $(\overline{r}, \overline{\mu})$ is potentially diagonalizably automorphic.\end{enumerate}
Then $(r, \mu)$ is potentially diagonalizably automorphic (of level potentially prime
to $\ell$).\end{thm*}

\begin{rem} Whenever we apply (ALT-PD), both representations will be potentially diagonalisable at all primes $v\vert\ell$ of the base field $F$.
Note that the assumption $\ell\geq 2(d+1)$ can be relaxed as in \cite{AppendixDiGe}, but not in the situation when $\ell\vert n$. In any case, we will apply
this theorem only for $\ell\geq 2(d+1)$.
\end{rem}

Next, we want to introduce some automorphy lifting theorems, based on results of \cite{Di-mfold}, which are specifically designed to prove automorphy of tensor products. But, before stating and proving these theorems, we discuss one of their main hypothesis: regularity of tensor products.

First, some notation on Hodge-Tate weights. If $(r,\mu)$ is a regular algebraic polarized $\ell$-adic $n$-dimensional representation of $G_F$, then it is de-Rham at all places above $\ell$ (see Section 2.1 in \cite{BLGGT}). In particular, for each embedding $\tau\colon F \hookrightarrow \overline{\Q}$ we can attach to $r$ a multiset $\mathrm{HT}_\tau(r)$ of $n$ integers, which are called the Hodge-Tate weights of $r$ with respect to $\tau$. One may write this multiset explicitly as $\{a_{\tau,1},\ldots,a_{\tau,n}\}$, although we will omit $\tau$ from the notation when the Hodge-Tate weights of $r$ are independent of the choice of the embedding $\tau$ (this will often be the case later in Sections 5 and 6).

Let $(r,\mu)$ and $(s,\eta)$ be two regular algebraic polarized $\ell$-adic Galois representations of $G_F$, of respective dimensions $n$ and $m$. Let $\tau\colon F \hookrightarrow \overline{\Q}$ be an embedding, and let $\{a_{\tau,1},\ldots,a_{\tau,n}\}$ and $\{b_{\tau,1},\ldots,b_{\tau,m}\}$ be, respectively, the HT weights of $r$ and of $s$ with respect to $\tau$. Then, the HT weights of $r\otimes s$ with respect to $\tau$ are 
\begin{equation*}
\{a_{\tau,i}+b_{\tau,j}\colon i=1,\ldots,n; j=1,\ldots,m\}.
\end{equation*}
If we want $r\otimes s$ to be regular, we need all these numbers to be different for each $\tau$. This motivates the following definitions, which are similar to the notions defined in \cite{Di-mfold}.

\begin{definition}
\label{def:max}
Let $n\geq1$. Let $(r,\mu)$ be a regular, algebraic, polarized $\ell$-adic $n$-dimensional representation of $G_F$. Let $\tau\colon F\hookrightarrow\overline{\Q}_\ell$, and suppose that the Hodge-Tate numbers of $r$ are
\begin{equation*}
\textup{HT}_{\tau}(r) = \{a_{\tau,1},\ldots,a_{\tau,n}\}.
\end{equation*}
Define the minimal difference in $r$, denoted ${\min}(r)$, to be $\min_{\tau,i,j}\{|a_{\tau,i}-a_{\tau,j}|\}$. Similarly define the maximal difference in $r$, denoted ${\max}(r)$, to be $\max_{\tau,i,j}\{|a_{\tau,i}-a_{\tau,j}|\}$.
\end{definition}

\begin{lemma}
\label{lemma:regularity}
Let $(r,\mu)$ and $(s,\eta)$ be as in the previous definition, and assume that ${\max}(s)<{\min}(r)$. Then $r\otimes s$ is regular. Moreover, if $r$ and $s$ are ordinary at a place $\nu|\ell$ of $F$, then so is $r\otimes s$.
\end{lemma}
\begin{proof}
Algebraicity of $r\otimes s$ follows from Fontaine's formalism. To see it is also regular, choose an embedding $\tau\colon F\hookrightarrow\overline{\Q}_\ell$, and write
\begin{align*}
	\textup{HT}_{\tau}(r) &= \{a_{\tau,1},\ldots,a_{\tau,n}\}, \\
	\textup{HT}_{\tau}(s) &= \{b_{\tau,1},\ldots,b_{\tau,m}\},
\end{align*}
for the Hodge-Tate numbers of $r$ and $s$ with respect to $\tau$. Then, the Hodge-Tate numbers of $r\otimes s$ with respect to $\tau$ are
\begin{equation*}
	\textup{HT}_{\tau}(r\otimes s) = \{a_{\tau,i}+b_{\tau,j} \mid 1\leq i \leq n, 1\leq j \leq m\}.
\end{equation*}
Regularity means that these numbers are distinct, or equivalently, that $a_{\tau,i}-a_{\tau,i'} \neq b_{\tau,j'} - b_{\tau,j}$ for any choice of indeces. This follows from the assumption ${\max}(s)<{\min}(r)$.

Assume now that both $r$ and $s$ are ordinary at an $\ell$-adic place $\nu$ of $F$ (see page 528 of \cite{BLGGT} for the definition of ordinary representations). In particular, this means that for each $i=1,\ldots,n$, the $i$-th graded piece of the filtration attached to $r|_{G_{F_\nu}}$ is one dimensional, and acted upon by a character $\chi_i$. Similarly, for each $j=1,\ldots,m$, the $j$-th graded piece of the filtration attached to $s|_{G_{F_\nu}}$ is one dimensional, and acted upon by a character $\phi_j$. Moreover, there exists an open subgroup $U\subset F_\nu^\times$ such that if $\alpha\in U$, then
\begin{align*}
(\chi_i\circ \textrm{Art}_{F_\nu})|_U (\alpha) = \prod_{\tau\colon F_\nu\hookrightarrow \overline{\mathbb{Q}}_\ell}(\tau\alpha)^{-a_{\tau,i}}, \\
(\phi_j\circ \textrm{Art}_{F_\nu})|_U (\alpha) = \prod_{\tau\colon F_\nu\hookrightarrow \overline{\mathbb{Q}}_\ell}(\tau\alpha)^{-b_{\tau,j}},
\end{align*}
where $\textrm{Art}_F$ denotes the local Artin reciprocity map. In addition, for each embedding $\tau\colon F_\nu\hookrightarrow\overline{\Q}_\ell$,
\begin{align*}
a_{\tau,1} > \ldots > a_{\tau,n}, \\
b_{\tau,1} > \ldots > b_{\tau,m}.
\end{align*}

To prove that $r\otimes s$ is also ordinary at $\nu$, we need to show it satisfies the same properties. Using Fontaine's formalism and the assumption that ${\max}(s)<{\min}(r)$, a simple calculation shows that the graded pieces of $r\otimes s|_{G_{F_\nu}}$ are

\begin{align*}
\textrm{Gr}^{1}(r\otimes s|_{G_{F_\nu}}) &= \textrm{Gr}^{1}(r|_{G_{F_\nu}}) \otimes \textrm{Gr}^{1}(s|_{G_{F_\nu}}), \\
\textrm{Gr}^{2}(r\otimes s|_{G_{F_\nu}}) &= \textrm{Gr}^{1}(r|_{G_{F_\nu}}) \otimes \textrm{Gr}^{2}(s|_{G_{F_\nu}}), \\
&\ldots \\
\textrm{Gr}^{m}(r\otimes s|_{G_{F_\nu}}) &= \textrm{Gr}^{1}(r|_{G_{F_\nu}}) \otimes \textrm{Gr}^{m}(s|_{G_{F_\nu}}), \\
\textrm{Gr}^{m+1}(r\otimes s|_{G_{F_\nu}}) &= \textrm{Gr}^{2}(r|_{G_{F_\nu}}) \otimes \textrm{Gr}^{1}(s|_{G_{F_\nu}}), \\
&\ldots \\
\textrm{Gr}^{2m}(r\otimes s|_{G_{F_\nu}}) &= \textrm{Gr}^{2}(r|_{G_{F_\nu}}) \otimes \textrm{Gr}^{m}(s|_{G_{F_\nu}}), \\
&\ldots \\
\textrm{Gr}^{nm}(r\otimes s|_{G_{F_\nu}}) &= \textrm{Gr}^{n}(r|_{G_{F_\nu}}) \otimes \textrm{Gr}^{m}(s|_{G_{F_\nu}}).
\end{align*}
In particular, these graded pieces are also one dimensional, and acted upon by the products $\chi_{i}\phi_{j}$. Furthermore,
\begin{equation*}
(\chi_i\phi_j\circ \textrm{Art}_{F_\nu})|_U (\alpha) = \prod_{\tau\colon F_\nu\hookrightarrow \overline{\mathbb{Q}}_\ell}(\tau\alpha)^{-(a_{\tau,i}+b_{\tau,j})}.
\end{equation*}
Lastly, we need to check that the exponents appearing in these formulas are decreasing as we move along the filtration. This means that if we fix an embedding $\tau\colon F_\nu\hookrightarrow\overline{\Q}_\ell$, then we want the exponent of $(\tau\alpha)$ under the action of $\chi_1\phi_1$ to be smaller than that of $\chi_1\phi_2$, which in turn should be smaller than that of $\chi_1\phi_3$, etc. In our case this follows from the assumption ${\max}(s)<{\min}(r)$, which implies that 
\begin{equation*}
a_{\tau,1}+b_{\tau,1} > \ldots > a_{\tau,1}+b_{\tau,m} > a_{\tau,2}+b_{\tau,1} > \ldots > a_{\tau,2}+b_{\tau,m}> \ldots > a_{\tau,n}+b_{\tau,m},
\end{equation*}
as desired.
\end{proof}

A particular situation that will arise often in practice is the case where $r$ is the restriction of a $2$-dimensional representation arising from a modular form, and $s$ is an $n$-dimensional representation whose Hodge-Tate weights with respect to an embedding $\tau$ are of the form
\begin{equation*}
\{0,a,2a,\ldots,(n-1)a\},
\end{equation*}
for some integer $a>0$ which is independent of $\tau$. In this case, the HT weights of $r$ with respect to $\tau$ are $\{0,k-1\}$, where $k$ is the weight of the modular form, and so, in order to check that $r\otimes s$ is regular, it is enough, by the previous lemma, to check that $a>k-1$. Let us give a name to this particular scenario.

\begin{defi}\label{def:C-very-spread}
 Let $\{0, k_1, \dots, k_{m-1}\}$ be a strictly increasing $m$-tuple of natural numbers
 in arithmetic progression, and $C\geq 2$ some natural number. We will say that
 $\{0, k_1, \dots, k_{m-1}\}$ is \emph{$C$-very spread} if $k_1>2C$.
\end{defi}

\begin{rem}
\label{rem:Regularity_tensor}
Note that, by the previous discussion, whenever $r$ and $s$ are $\ell$-adic, regular, algebraic representations of $G_F$, of respective dimensions $2$ and $n$, and there exists positive integers $k<C$ such that:
\begin{itemize}
\item $\mathrm{HT}_\tau(r)=\{0, k\}$ for any embedding $\tau\colon F\hookrightarrow \overline{\Q}_\ell$,
\item $\mathrm{HT}_\tau(s)$ is $C$-very spread for any embedding $\tau\colon F\hookrightarrow \overline{\Q}_\ell$,
\end{itemize}
then $r\otimes s$ is regular algebraic. The reason why we use $2C$ instead of $C$ in the definition will be explained in Remark \ref{rem:Constant2}.

\end{rem}

Next we will state and prove the following variations of Theorems 2.1 and 2.2 of \cite{Di-mfold} (in particular, these variants will apply to the case when $\ell$ divides the dimension of the representation). The key ingredient that makes these proofs possible is the following theorem, a minimal automorphy lifting theorem which extends Theorem 7.1 of \cite{Th2012} (together with its correction and extension to the $\ell| n$ case given in Corollary 7.3 of \cite{Thorne2016}) to semistable representations. The proof of this theorem is relegated to Sections \ref{section_1}, \ref{section_2} and \ref{section_3}, but we state it here for convenience of the reader.

\begin{theorem}[ALT-SEMISTABLE]
\label{alt_semistable}
Let $n\in\Z_{\geq1}$, and fix an odd prime $\ell$. Let $F$ be an imaginary CM field with totally real subfield $F^+$, and assume that $\zeta_\ell\notin F$. Suppose that $(r,\mu)$ is an $n$-dimensional, algebraic, polarized $\ell$-adic representation of $G_F$ satisfying the following properties:
\begin{enumerate}
\item $\overline{r}$ is irreducible and maps $G_{F(\zeta_\ell)}$ to an adequate subgroup (adequate as in Definition 2.20 of \cite{Thorne2016}, so in particular it could be that $\ell\mid n$).
\item $(\overline{r},\overline{\mu})$ is automorphic, arising from a regular algebraic, cuspidal, polarized automorphic representation $(\pi,\chi)$ such that at all places $\nu$ of $F$
\begin{equation*}
r|_{G_{F_\nu}} \sim r_\ell(\pi)|_{G_{F_\nu}}
\end{equation*}
(that is, the representations are connected, as defined in Definition \ref{def:connected}). In particular, locally at $\nu\mid\ell$ both representations are potentially semistable with equal Hodge-Tate weights.
\end{enumerate}
Then $(r,\mu)$ is automorphic.
\end{theorem}

\begin{remark}
To the authors knowledge, the first instance of an automorphy lifting theorem in the semistable case is Theorem 7.1 of \cite{Calegari}, who proves it in the case where $\pi$ has slightly regular weight. Here we remove this assumption, altough we claim no great originality, since the proof is essentially that of Theorem 7.1 of \cite{Th2012} (with the correction and extension to the $\ell\mid n$ case as in Corollary 7.3 of \cite{Thorne2016}).

As Calegari and other experts were aware at the time, the strong form of local-global compatibility, together with a finiteness result on inertia types, would allow to extend the results of Thorne to the potentially semistable case (without the slightly regular weight assumption). Since the publishing of \cite{Calegari}, the strong form of local-global compatibility was proven by Caraiani in \cite{Caraiani2012} and \cite{Caraiani}, making this extension finally possible.
\end{remark}

\begin{thm}[ALT-MIXED]\label{thm:mixto} Let $F$ be a CM field, let $\ell$ be an odd prime number such that $\zeta_{\ell}\not\in F$,
and let  $(R, \mu_0)$, $(R', \mu'_0)$ (resp.~$(S, \mu_1)$, $(S', \mu'_1)$) be
$n_0$-dimensional (resp.~$n_1$-dimensional) regular, algebraic, polarized $\ell$-adic representations
of $G_F$ which are irreducible.
Assume that the tensor product $R'\otimes S'$ is automorphic and assume that $(R, \mu_0)$ and $(R', \mu'_0)$ (resp.~$(S, \mu_1)$ and $(S', \mu'_1)$) are congruent mod $\ell$.

Assume furthermore the following technical conditions:

\begin{enumerate}
 \item $R\otimes S$ and $R'\otimes S'$ are regular and irreducible,
 \item For every  $v\vert \ell$ place of $F$, $R\vert_{G_{F_v}}$, $R'\vert_{G_{F_v}}$ are ordinary and $S\vert_{G_{F_v}}$, $S'\vert_{G_{F_v}}$
 are potentially diagonalisable.
 \item The image of $\overline{(R\otimes S)}\vert_{G_{F(\zeta_{\ell})}}$ is adequate$+$.
\end{enumerate}
Then $R\otimes S$ is automorphic.

\end{thm}

\begin{proof} 
We expand on the details of the proof outlined in Theorem 2.1 in \cite{Di-mfold}, which combines ideas of the proofs of Proposition 4.1.1 and Theorem 4.2.1 of \cite{BLGGT}. As in loc. cit., the theorem will follow from Theorem 2.4.1 of loc. cit. after replacing $S$ and $S'$ by appropiate ordinary representations $\xi$ and $\xi'$.

Let us start by making some preliminary reductions. Since soluble base change preserves automorphy (see Lemma 2.2.2 of loc. cit.), it is enough to prove the theorem after replacing $F$ by a soluble CM extension. So using Lemma A.2.1 of loc. cit. we replace $F$ by a soluble extension which is linearly disjoint from $\overline{F}^{\Ker(\overline{R\otimes S})}(\zeta_\ell)$ and such that:
\begin{enumerate}[itemsep=0.5em]
\item $F/F^+$ is unramified at all finite primes,
\item all primes dividing $\ell$, and primes at which $S$ or $S'$ ramify, are split over $F^+$,
\item if $\nu$ is a place in $F$ above $\ell$, then $F_\nu$ contains a primitive $\ell$-th root of unity, and the restrictions $\overline{S}|_{G_{F_\nu}}$ and $\overline{S'}|_{G_{F_\nu}}$ are both trivial.
\end{enumerate}

With these reductions, we are ready to construct the ordinary representation $\xi$ that will replace $S$. Choose an integer $t>\max(R\otimes S)$, $\max(R)$ and $\max(S)$ (see Definition \ref{def:max}). The reason to take $t$ this large is to guarantee regularity and ordinarity of some tensor products. Given $\tau\colon F\hookrightarrow \overline{\Q}_\ell$, set
\begin{equation*}
\textup{H}_{\tau} = \{0,t,\ldots,(n_1-1)t\}.
\end{equation*}
Notice that for primes above $\ell$ the reduction $\overline{S}$ has an ordinary and crystalline lift $1\oplus\epsilon^{-t}\oplus\ldots\oplus\epsilon^{(1-n_1)t}$, so Proposition 3.2.1 of loc. cit, gives us a continuous homomorphism 
\begin{equation*}
\xi \colon G_{F} \rightarrow \mathrm{GL}_{n_1}(\overline{\Q}_\ell)
\end{equation*}
lifting $\overline{S}$ and such that:
\begin{enumerate}[itemsep=0.5em]
\item if $\nu|\ell$ then $\xi$ is ordinary and crystalline at $\nu$, with Hodge-Tate numbers $\textup{H}_{\tau}$,
\item if $\nu\nmid \ell$ and $S$ ramifies at $\nu$, then $S$ and $\xi$ are connected at $\nu$,
\item if $\nu\nmid \ell$ and $S$ is unramified at $\nu$, then $\xi$ is also unramified at $\nu$ (in particular, this implies that $S$ and $\xi$ are connected at $\nu$, see p. 524 of \cite{BLGGT}).
\end{enumerate}

At this point we would like to conclude that $R\otimes S$ is automorphic iff so is $R\otimes \xi$. However, since $S$ and $\xi$ are not connected at primes dividing $\ell$, we can not apply a minimal automorphy lifting theorem right away. In order to do so, we first need to apply Harris tensor trick (tensoring by induced characters to deduce automorphy).

Note that $S$ is potentially diagonalizable at places $\nu|\ell$. Also, $\xi$ is potentially diagonalizable at places $\nu|\ell$ by Lemma 1.4.3 of loc. cit. Thus, repeating the arguments in the proof of Proposition 4.1.1 of loc. cit., we may find soluble extensions $M$ and $F_1$ of $F$, and continuous characters $\theta,\theta'\colon G_M\rightarrow \overline{\Q}_\ell^\times$ such that:
\begin{enumerate}[itemsep=0.5em]
\item There is a congruence $\textup{Ind}_{G_M}^{G_F}\theta|_{G_{F_1}}\equiv \textup{Ind}_{G_M}^{G_F}\theta'|_{G_{F_1}} \Mod \ell$.
\item The Hodge-Tate weights of $\textup{Ind}_{G_M}^{G_F}\theta|_{G_{F_1}}$ match those of $\xi$, and the Hodge-Tate weights of $\textup{Ind}_{G_M}^{G_F}\theta'|_{G_{F_1}}$ match those of $S$.
\item The representations $(S\otimes \textup{Ind}_{G_M}^{G_F}\theta)|_{G_{F_1}}$ and $(\xi\otimes \textup{Ind}_{G_M}^{G_F}\theta')|_{G_{F_1}}$ are connected at all finite places (note that these representations are regular by our choice of $t$ and Lemma \ref{lemma:regularity}).
\item The representation $\overline{(R\otimes S\otimes \textup{Ind}_{G_M}^{G_F}\theta)}|_{G_{F_1}} \Mod \ell$ is irreducible and maps $F_1(\zeta_\ell)$ to an adequate subgroup.
\end{enumerate}
Consider the tensor products
\begin{equation*}
(R\otimes S\otimes \textup{Ind}_{G_M}^{G_F}\theta)|_{G_{F_1}}, \quad  (R\otimes \xi\otimes \textup{Ind}_{G_M}^{G_F}\theta')|_{G_{F_1}}.
\end{equation*}
By construction, these two representations are congruent modulo $\ell$, and moreover, they are also connected. Indeed, since we are tensoring two connected representations by the same representation, the resulting representations are still connected. Furthermore, by our choice of $t$, both tensor products are regular by Lemma \ref{lemma:regularity}. Thus, we can apply Theorem \ref{alt_semistable} (minimal automorphy lifting theorem in the semistable case) to conclude that one of these representations is automorphic if and only if so is the other. Lemmas 2.2.4, 2.2.1 and 2.2.2 of \cite{BLGGT} (respectively, these say that induction preserves automorphy, that tensoring by CM characters preserves automorphy, and, that soluble base change preserve automorphy) now imply that $R\otimes S$ is automorphic if and only if so is $R\otimes \xi$.

To conclude the proof, we repeat this process with $R'$ and $S'$, so there exists an ordinary and crystalline representation $\xi'$ such that $R'\otimes \xi'$ is automorphic if and only if so is $R'\otimes S'$. Hence, by assumption, $R'\otimes \xi'$ is automorphic. Then, the representations $R\otimes \xi$ and $R'\otimes \xi'$ are both ordinary (by Lemma \ref{lemma:regularity} and the choice of $t>\max(R)$ for $R\otimes\xi$, similarly for $R'\otimes\xi'$) and congruent modulo $\ell$, so by Theorem 2.4.1 of loc. cit. it follows that $R\otimes \xi$ is automorphic, which implies that $R\otimes S$ is automorphic, as desired.

\end{proof}

\begin{rem} We will apply this theorem in the case $\ell\in \mathcal{L}$
(see Section \ref{sec:PreliminarySteps} for the definition of $\mathcal{L}$); in particular we will always have $\ell\geq 7$.
\end{rem}

\begin{thm}[ALT-HARRIS]\label{thm:pd} Let $F$ be a CM field, let $\ell$ be an odd prime number such that $\zeta_{\ell}\not\in F$,
and let  $(R, \mu_0)$, $(R', \mu'_0)$ (resp. $(S, \mu_1)$, $(S', \mu'_1)$) be
$n_0$-dimensional (resp. $n_1$-dimensional) regular, algebraic,  polarized $\ell$-adic representations
of $G_F$ which are irreducible. Assume that the tensor product $R'\otimes S'$ is automorphic and assume that $(R, \mu_0)$ and $(R', \mu'_0)$ (resp. $(S, \mu_1)$ and $(S', \mu'_1)$)
are congruent mod $\ell$. Assume furthermore the following technical conditions:

\begin{enumerate}
\item $\ell\nmid n_0$.
 \item $R\otimes S$ and $R'\otimes S'$ are regular and irreducible.
 \item For every  $v\vert \ell$ place of $F$, $R\vert_{G_{F_v}}$, $R'\vert_{G_{F_v}}$
 are potentially diagonalizable and for every $v$ not dividing
 $\ell$, $R\vert_{G_{F_v}}$ and $R'\vert_{G_{F_v}}$ are connected.
 \item  $S\vert_{G_{F_v}}$ and $S'\vert_{G_{F_v}}$ are connected, for all finite place $v$ of $F$.
 \item The image of $\overline{(R\otimes S)}\vert_{G_{F(\zeta_{\ell})}}$ is adequate$+$.
 \item $R\otimes R'\otimes S$ is regular.
\end{enumerate}
Then $R\otimes S$ is automorphic.

\end{thm}

\begin{proof}
We expand on the details of the proof outlined in Theorem 2.2 of \cite{Di-mfold}, which follows closely Proposition 4.1.1 of \cite{BLGGT}.

Let us start by making some preliminary reductions. Since soluble base change preserves automorphy (Lemma 2.2.2 of loc. cit.), we can replace $F$ by a soluble CM extension. Using Lemma A.2.1 of loc. cit., we replace $F$ by a soluble extension, linearly disjoint from $\overline{F}^{\Ker(\overline{R\otimes S})}(\zeta_\ell)$, and such that:
\begin{enumerate}[itemsep=0.5em]
\item $F/F^+$ is unramified at all finite primes,
\item all finite places which divide $\ell$ or at which $R$ or $R'$ ramify are split over $F^+$,
\item if $\nu$ is a place of $F$ dividing $\ell$ or at which $R$ or $R'$ ramify, then $\overline{R}|_{G_{F_\nu}}$ and $\overline{R}|_{G_{F_\nu}}$ are trivial.
\item if $\nu$ is a place of $F$ dividing $\ell$, then both $R|_{G_{F_\nu}}$ and $R'|_{G_{F_\nu}}$ are diagonalizable.
\end{enumerate}
Arguing as in the proof of Proposition 4.1.1 of loc. cit., we now take soluble extensions $M$ and $F_1$ of $F$, and continuous characters $\theta,\theta'\colon G_M\rightarrow \overline{\Q}_\ell^\times$ such that: 

\begin{enumerate}
\item There is a congruence $\textup{Ind}_{G_M}^{G_F}\theta|_{G_{F_1}}\equiv \textup{Ind}_{G_M}^{G_F}\theta'|_{G_{F_1}} \Mod \ell$.
\item The Hodge-Tate weights of $\textup{Ind}_{G_M}^{G_F}\theta|_{G_{F_1}}$ match those of $R'$, and the Hodge-Tate weights of $\textup{Ind}_{G_M}^{G_F}\theta'|_{G_{F_1}}$ match those of $R$.
\item The representations $(\textup{Ind}_{G_M}^{G_F}\theta \otimes R)|_{G_{F_1}}$ and $(\textup{Ind}_{G_M}^{G_F}\theta' \otimes R')|_{G_{F_1}}$ are connected at all finite places.
\item The representation $(\textup{Ind}_{G_M}^{G_F}\theta \otimes R)|_{G_{F_1}}$ is irreducible and maps $F_1(\zeta_\ell)$ to an adequate subgroup.
\end{enumerate}
 
Consider now the tensor products
\begin{equation*}
 (\textup{Ind}_{G_M}^{G_F}\theta \otimes R\otimes S)|_{G_{F_1}}, \quad  (\textup{Ind}_{G_M}^{G_F}\theta' \otimes R'\otimes S')|_{G_{F_1}}.
\end{equation*}
By construction, these representations are congruent, and moreover, they are connected at all finite places (so in particular potentially semistable). Indeed, the latter follows because we are tensoring connected representations by connected representions. Also, they are regular (this follows from point 2 above and hypothesis 6). Thus, we meet the conditions to apply Theorem \ref{alt_semistable} (automorphy lifting for connected potentially semistable representations), and conclude that one of them is automorphic if and only if the other one is.

Lemmas 2.2.4, 2.2.1 and 2.2.2 of \cite{BLGGT} (respectively, these say that induction preserves automorphy, that tensoring by CM characters preserves automorphy, and, that soluble base change preserve automorphy) let us conclude that $R\otimes S$ is automorphic if and only if so is $R'\otimes S'$, concluding the proof.

\end{proof}

\begin{rem} Note that in the ALT-HARRIS it is allowed that the dimension of the representation $S$ is divisible by $\ell$.
\end{rem}

\begin{rem}\label{rem:Constant2} In the applications in the following sections, $R$ and $R'$ will be $2$-dimensional representations attached to
modular forms belonging to the chain $\mathcal{C}$ constructed in \cite{Di-Symm5}, hence their weights $\{0, k\}$ and
$\{0, k'\}$ will satisfy that $k$ and $k'$ will be
smaller than a certain constant $C_0$, and furthermore $k\not=k'$. We will work with representations $S$ with $C_0$-very
spread weights. It follows that condition 6 in the theorem will be satisfied (this is the reason for the $2$ multiplying
the constant in Definition \ref{def:C-very-spread}). \end{rem}

At certain points in the paper, we will need to perform change of level and weight of automorphic
representations. The following result will be applied to modular forms.

\begin{thm}[Big weight]\label{thm:BigWeight}
Let $f$ be a newform of weight $k$ and let $p>5$ be a prime number. Assume that the  image of
$\bar{\rho}_p(f)\vert_{\mathbb{Q}(\zeta_p)}$ is irreducible. Assume that $\rho_p(f)$ locally at
$p$ is potentially diagonalisable. We twist
$\overline{\rho}_{p}(f)$ in such a way that the Serre weight $k$ satisfies $k\leq p+1$. Assume $k\not= p$,
and fix a positive integer $D$. Then, up to twisting by a finite order character ramified only at $p$,
$\rho_{p}(f)$ is congruent to another modular Galois representation $\rho_{p}(f')$ for a newform $f'$
of weight $k'\in [D, p^2-1+D]$. $f'$ has level prime to $p$ and $\rho_{p}(f')$ is (crystalline)
and potentially diagonalisable at $p$.

Moreover, we can choose $f'$ such that the $p$-adic Galois representation $\rho_{p}(f)$ and $\rho_{p}(f')$
connect locally at each $\ell\not=p$.
\end{thm}

\begin{rem}\label{rem:k_distinto_de_p}
In all applications of this theorem, we are going to have a newform $f$ of even weight such that $p$
is unramified in the nebentypus of $f$. From this it is clear that the ``minimal'' residual Serre
weight will be even, so in particular the condition $k\not=p$ is satisfied.
\end{rem}

\begin{rem}\label{rem:3.11}
In the Fontaine Laffaille situation, the potentially Barsotti-Tate situation or the ordinary crystalline case,
the representation is potentially diagonalisable (see \cite{BLGGT} and \cite[Lemma 4.4.1]{GeeKisin}
for the potentially Barsotti-Tate case).
\end{rem}

\begin{proof}[Proof of Theorem \ref{thm:BigWeight}]
The proof is similar to the one of Lemma 4.1 in \cite{Di-mfold} (compare also with Theorem 0.1 in \cite{Blanco2021}). We just indicate the differences with that proof.
Call $\overline{\rho}$ the twist of $\bar{\rho}_{p}(f)$ of minimal Serre weight. We need to find a modular lift of
$\bar{\rho}$ attached to a newform $f'$ with the required conditions. We consider two cases.

\begin{itemize}

\item Case 1: Assume $k(\bar{\rho})=p+1$. In this case, we know that there is a modular lift of weight $p+1$,
ordinary at $p$, and by Hida theory, we can also find a crystalline ordinary lift (thus potentially diagonalisable)
corresponding to a modular form $f'$ of weight $k'\in [D, p^2-1+D]$.

\item Case 2: Assume $k(\bar{\rho})\leq p-1$. In this case we know that there is a modular crystalline lift of weight $k\leq p-1$,
thus it falls in the Fontaine Laffaille case. If $\overline{\rho}$ is reducible,  this lift is ordinary (because we are in the
Fontaine-Laffaille range), hence we can apply Hida Theory
like in Case 1 to find a crystalline ordinary lift corresponding to
a modular form $f'$ whose weight is comprised in the interval $[D, p^2-1+D]$.

If $\overline{\rho}$ is irreducible, we proceed as in \cite{Di-mfold}. Namely, using Lemma 4.19 of \cite{BLGG}, we can find a
potentially diagonalisable and crystalline modular lift of arbitrarily large weight $k'$. In fact, since the only restriction
on this weight is given by the fact that the residual Serre weight is fixed, the arguments in Lemma 4.19 of loc.~cit. (using the fact that a fundamental character of level two has order $p^2-1$) allows us to pick
$k'\in [D, (p^2-1)+ D]$.

\end{itemize}
\end{proof}

We will apply Theorem 4.4.1 of \cite{BLGGT} to produce automorphic lifts of residually automorphic
representations with prescribed local types several times. We recall the setup. Let $F$ be a CM field with maximal totally real subfield $F^{+}$,
 $n$ a positive integer and $\ell$ an odd rational prime, satisfying that $\ell>2(n+1)$, $\zeta_{\ell}\not\in F$ and, for
all primes $v\vert \ell$ of $F^{+}$, $v$ is completely split in $F/F^{+}$. Fix an embedding $\iota:\overline{\mathbb{Q}}_{\ell}\simeq \mathbb{C}$.
Let $S$ be a finite set of primes of $F^{+}$, which are split in $F/F^{+}$, and containing all those places of $F^+$ above $\ell$.
Given a continuous representation $\overline{r}:G_F\rightarrow \GL_n(\overline{\mathbb{F}}_{\ell})$,
with $\overline{r}\vert_{G_{F(\zeta_{\ell})}}$ irreducible, and such that for some algebraic character $\mu$ of $G_{F^+}$ the representation $(\overline{r}, \overline{\mu})$ is a polarised
mod $\ell$ representation unramified outside $S$, which is either ordinarily automorphic or potentially diagonalisably automorphic, then
Theorem 4.4.1 of \cite{BLGGT} ensures the existence of a regular, algebraic, cuspidal, polarised automorphic representation $(\pi, \chi)$ of
$\GL_n(\mathbb{A}_F)$, satisfying the following conditions:

\begin{enumerate}
\item $\overline{r}_{\ell, \iota}(\pi)\simeq \overline{r}$;
\item $r_{\ell, \iota}(\chi)\varepsilon^{1-n}_{\ell}=\mu$;
\item The level of $\pi$ is potentially prime to $\ell$;
\item $\pi$ is unramified outside $S$;
\item For any prefixed choice of a prime $\widetilde{v}$ of $F$ above each of the primes $v\in S$ above $\ell$, and for each prefixed lift
$\rho_v:G_{F_{\widetilde{v}}}\rightarrow \GL_n(\mathcal{O}_{\overline{\mathbb{Q}}_{\ell}})$ of $\overline{r}_{\ell, \iota}(\pi)\vert_{G_{F_{\widetilde{v}}}}$ which is
potentially diagonalisable and such that, for each embedding $\tau:F_{\widetilde{v}}\rightarrow \overline{\mathbb{Q}}_{\ell}$, the Hodge Tate weights $\mathrm{HT}_{\tau}(\rho_v)$
are all different, it holds
$$r_{\ell, \iota}(\pi)\vert_{G_{F_{\widetilde{v}}}}\sim \rho_v$$
\item  For any prefixed choice of a prime $\widetilde{v}$ of $F$ above each of the primes $v\in S$ which are not above $\ell$, and for each prefixed lift
$\rho_v:G_{F_{\widetilde{v}}}\rightarrow \GL_n(\mathcal{O}_{\overline{\mathbb{Q}}_{\ell}})$ of $\overline{r}_{\ell, \iota}(\pi)\vert_{G_{F_{\widetilde{v}}}}$, it holds
$$r_{\ell, \iota}(\pi)\vert_{G_{F_{\widetilde{v}}}}\sim \rho_v$$
\end{enumerate}

\begin{rem}
On the one hand, assume that $v\vert \ell$ is a place of $F^{+}$. Since the relationship ``connects to'' is an equivalence relationship
(cf.~\cite[$(2)$ in page 530]{BLGGT}), and we chose a lift $\rho_v$ which is potentially diagonalisable, we have that
$r_{\ell, \iota}(\pi)\vert_{G_{F_{\widetilde{v}}}}$ is also potentially diagonalisable.
Moreover, this relationship preserves the restriction of the Weil-Deligne representation to the inertia group at $\widetilde{v}$
(cf.~\cite[$(3)$ in page 530]{BLGGT}). Hence if we choose $\rho_v$ to be crystalline, then $r_{\ell, \iota}(\pi)\vert_{G_{F_{\widetilde{v}}}}$
is crystalline too, thus $\ell$ does not appear in the level of $\pi$.

On the other hand, assume that $v\nmid \ell$ is a place of $F^{+}$ belonging to $S$. Theorem 4.4.1 of \cite{BLGGT} ensures the existence of
$\pi$ such that $r_{\ell, \iota}(\pi)\vert_{G_{F_{\widetilde{v}}}}$ connects to a prefixed lift $\rho_v$ of
$\overline{r}\vert_{G_{F_{\widetilde{v}}}}$. By Lemma 1.3.2 of \cite{BLGGT},  $r_{\ell, \iota}(\pi)\vert_{G_{F_{\widetilde{v}}}}$ will belong
to a unique irreducible component of the corresponding deformation ring. Moreover the same reasoning
shows that we can always take the lift $\rho_v$ belonging to a unique irreducible component of the
corresponding deformation ring (which we will always do, without mentioning it explicitly). Thus, with
this choice $\rho_1=r_{\ell, \iota}(\pi)\vert_{G_{F_{\widetilde{v}}}}$ and  $\rho_2=\rho_v$ are actually
strongly connected (i.e.~$\rho_1\rightsquigarrow \rho_2$ and $\rho_2\rightsquigarrow \rho_1$, cf.~\cite[pag 524]{BLGGT}).

In particular, \cite[page 524, $(6)$]{BLGGT} shows that the local types $(r\vert_{I_v}, N_1)$ of $\rho_1$ and $(r\vert_{I_v}, N)$ of $\rho_2$ at $v$ agree.

\end{rem}

Next, we want to record a special situation in which,
applying Theorem A in \cite{GHLS}, one can produce a weight $0$ lift (this means, under our conventions on HT weights, that the HT weights of the lift with respect to any embedding $\tau$ are $\{0,1,\ldots,n-1\}$).

\begin{thm}\label{thm:Weight0}[Weight $0$ Thm] Let $F$ be a CM field, $F^{+}$ its maximal totally real subfield and $\ell$ a
rational prime, $n$ a positive integer such that the following conditions are satisfied: $\ell>2(n+1)$,
$\ell$ is unramified in $F$ and all primes  of $F^+$ above $\ell$ split in $F/F^{+}$.

Let $(\pi, \chi)$ be a regular algebraic polarised cuspidal automorphic representation of $\mathrm{GL}_n(\mathbb{A}_F)$,
and consider the set $S$ of primes of $F^{+}$ above $\ell$ and the primes dividing the level of $\pi$. We assume that
all $v\in S$ are split in $F/F^+$. For each $v\in S$, fix a prime $\tilde{v}$ of $F$ above $v$.
Let $r=r_{\ell, \iota}(\pi)$ and $\mu=r_{\ell,\iota}(\chi)
$. Assume that $r$ is Fontaine Laffaille locally at each $v\vert \ell$ (in particular, $\ell$ is unramified for $\pi$).
Assume that $\bar{r}\vert_{G_{F(\zeta_{\ell})}}$ is absolutely irreducible. For each $v\in S$, we will fix a lift $\rho_v$ of $\overline{r}\vert_{G_{F_{\tilde{v}}}}$ as follows:
\begin{itemize}
\item If $v\nmid \ell$, $\rho_v=r\vert_{G_{F_{\tilde{v}}}}$.

\item If $v\vert \ell$, we take a lift $\rho_v$ which is potentially diagonalisable and is of Hodge type $0$ (whose existence follows from Theorem A of \cite{GHLS}).

\end{itemize}

Then there is an automorphic RACP lift $(\pi', \chi')$ such that:

\begin{enumerate}

\item $\overline{r}_{\ell, \iota}(\pi')\simeq \overline{r}$.

\item The level of $\pi'$ is potentially prime to $\ell$, and furthermore $\pi'$ is unramified outside $S$.

\item For each $v\in S$, $r_{\ell, \iota}(\pi')\vert_{G_{F_{\widetilde{v}}}}\sim \rho_v$. In particular, $r_{\ell, \iota}(\pi')$ has weight $0$.

\end{enumerate}

\end{thm}

\subsection{Proof of ALT Semistable: Galois Deformation Theory}
\label{section_1}

%
%
In this section we recall the deformation theory of Galois representations that will be used when patching. The conventions that we follow are those of  \cite{BLGGT}. In particular $F$ denotes a CM field, or a totally real field, and $F^+$ denotes its maximal real subfield. Let $S$ be a finite set of 
primes in $F^+$ which are split in $F$. For each prime in $S$ pick a prime in $F$ lying above it, and denote by $\tilde{S}$ the set of these places. Write $F_S$ for the maximal extension of $F$ unramified outside $S$. Denote by $G_{F^+,S}$ the Galois group $\Gal(F_S/F^+)$, and similarly write $G_{F,S}$ for $\Gal(F_S/F)$.

Let $\ell$ be a rational prime. The $\ell$-adic representations of $G_{F,S}$ that we are concerned with are polarized, in the sense of Section 2.1 of \cite{BLGGT}. More precisely, let $\mathcal{G}$ be the semidirect product of $\GL_n\times\GL_1$ with $\{1,\jmath\}$ given by
\begin{equation*}
\jmath(g,a)^\jmath\jmath^{-1} = (a^tg^{-1},a),
\end{equation*}
and let $\nu\colon\mathcal{G}\rightarrow\GL_1$ be the map sending $(g,a)$ to $a$ and $\jmath$ to $-1$. A pair $(r,\mu)$ of continuous homomorphisms
\begin{equation*}
r \colon G_{F,S} \rightarrow \GL_n(\overline{\Q}_\ell), \quad \mu \colon G_{F^+,S} \rightarrow \overline{\Q}_\ell^\times
\end{equation*}
is said to be \emph{polarized} if there exists a continuous homomorphism $\rho\colon G_{F^+,S}\rightarrow \mathcal{G}_n$ such that $\nu\circ\rho=\chi$, and such that composition with projection onto $\GL_n$, then restriction to $G_{F,S}$ gives $r$. Polarized representations with values in finite fields are defined similarly.

Let $L$ be a finite extension of ${\Q}_\ell$ with ring of integers $\mathcal{O}$, maximal ideal $\lambda$, and residue field $\mathds{F}$. Consider a representation $\overline{r}\colon G_{F^+,S}\rightarrow\mathcal{G}_n(\mathds{F})$ such that $\overline{r}^{-1}(\GL_n\times\GL_1)=G_{F,S}$, and choose a character $\chi\colon G_{F^+,S}\rightarrow \mathcal{O}^\times$ such that $\nu\circ\overline{r}=\overline{\chi}$. The lifts of $\overline{r}$ that we will consider will be valued in the category of complete Noetherian local $\mathcal{O}$-algebras with residue field $\mathds{F}$, which will be denoted $\mathcal{C}_\mathcal{O}$.

%
%
\label{local deformation rings}
Let us start by recalling the local deformation theory for $\overline{r}$. Let $v\in S$, and denote by $\overline{r}|_{G_{F_{\tilde{v}}}}$ the composite 
\begin{equation*}
G_{F_{\tilde{v}}} \rightarrow G_{F,S} \rightarrow \mathcal{G}_n(\mathds{F}) \rightarrow \GL_n(\mathds{F}).
\end{equation*}

\begin{definition}
A lifting of $\overline{r}|_{G_{F_{\tilde{v}}}}$ to an object $R$ in $\mathcal{C}_\mathcal{O}$ is a continuous homomorphism $r\colon G_{F_{\tilde{v}}} \rightarrow \GL_n(R)$ with $r \equiv \overline{r}|_{G_{F_{\tilde{v}}}} \Mod{\mathfrak{m}_R}$. Two liftings of $\overline{r}|_{G_{F_{\tilde{v}}}}$ to $R$ are said to be equivalent if they are conjugate by an element of $1+\textup{M}_n(\mathfrak{m}_R)$.
\end{definition}

\begin{proposition}
The functor $\mathcal{D}_\nu^\Box$ from $\mathcal{C}_\mathcal{O}$ to \textup{SETS} mapping a ring $R$ to the set of liftings of $\overline{r}|_{G_{F_{\tilde{v}}}}$ to $R$ is represented by an object $R_\nu^\Box$ of $\mathcal{C}_\mathcal{O}$, which is called the local universal deformation ring at $\nu$.
\end{proposition}
\begin{proof}
This goes back to the work of Mazur. See Proposition 1.3.1 of \cite{Boeckle2013} for a proof.
\end{proof}

%
%
\begin{definition}
\label{local deformation problems}
A local deformation problem is a subfunctor $\mathcal{D}_\nu$ of $\mathcal{D}_\nu^\Box$ satisfying the following properties:
\begin{enumerate}[itemsep=0.5em]
\item $\overline{r}|_{G_{F_{\tilde{v}}}}$ lies in $\mathcal{D}_\nu(\mathds{F})$.
\item $\mathcal{D}_\nu$ is closed under equivalence of liftings.
\item if $r\in D_\nu(R)$ and $f\colon R\rightarrow S$ is an morphism in $\mathcal{C}_\mathcal{O}$, then $f\circ r\in D_\nu(S)$.
\item if $f\colon R\hookrightarrow S$ is an injective morphism in $\mathcal{C}_\mathcal{O}$ and $r\colon G_{F_{\tilde{v}}} \rightarrow \GL_n(R)$ is a lifting such that $r\in D_\nu(S)$, then $r\in D_\nu(R)$.
\item suppose given $r_1\in \mathcal{D}_\nu(R_1)$ and $r_2\in \mathcal{D}_\nu(R_2)$, and closed ideals $I_1\subset R_1$ and $I_2\subset R_2$ such that there is an isomorphism $R_1/I_1\simeq R_2/I_2$ under which $r_1 \Mod{I_1} \equiv r_2 \Mod{I_2}$. Let $R_3$ denote the subring of $R_1\times R_2$ consisting of pairs with same image in $R_1/I_1\simeq R_2/I_2$. Then $r_1\times r_2\in D_\nu(R_3)$.
\item if $r_j\in D_\nu(R_j)$ is an inverse system, then $\varprojlim r_j \in D_\nu(\varprojlim R_j)$.
\end{enumerate}
\end{definition}

\begin{lemma}
Let $I$ be an $1+\textup{M}_n(\mathfrak{m}_{R_\nu^\Box})$-invariant ideal of $R_\nu^\Box$ such that $I\subsetneqq\mathfrak{m}_{R_\nu^\Box}$, and the quotient $R_\nu^\Box/I$ is reduced. Then the subfunctor $\mathcal{D}_\nu$ of $\mathcal{D}_\nu^\Box$ defined by liftings of $\overline{r}|_{G_{F_{\tilde{v}}}}$ to $R$ whose corresponding map $R_\nu^\Box\rightarrow R$ has kernel containing $I$ is a local deformation problem. Moreover, every local deformation problem arises in this way.
\end{lemma}
\begin{proof}
This is Lemma 2.2.3 of \cite{CHT2008} with the corrigenda in Lemma 3.2 of \cite{PRIMS}.
\end{proof}

In what remains of this section we will recall from \cite{Th2012} the local deformation problems that will appear during the patching argument. These are divided in two classes, depending on whether $\nu|\ell$ or not. The main difference of this paper with loc. cit. is that the local deformation rings we consider at $\nu|\ell$ parametrize potentially semistable representations instead of crystalline representations.

Let us start by considering the local deformation problems at places $\nu|\ell$. The potentially semistable deformation rings are described by specifying a Hodge-Tate type and the extension which kills a specific inertia type. To be more precise, for each $\tau\in\Hom_{\Q_\ell}(F_{\tilde{\nu}},\overline{\Q}_\ell)$ let $\lambda_\tau=\{\lambda_{\tau,1},\ldots,\lambda_{\tau,n}\}$ denote a multiset of $n$ integers, which we assume to satisfy
\begin{equation*}
\lambda_{\tau,1} < \lambda_{\tau,2} < \ldots < \lambda_{\tau,n}.
\end{equation*}
Then we will refer to the collection $\lambda=\{\lambda_\tau\}_\tau$ as an $\ell$-adic Hodge-Tate type. Let $E/F_\nu$ be a finite extension.

%
%
\begin{theorem}
\label{semistable_rings}
There exists a quotient of $R_\nu^\Box$, denoted $R_\nu^{E-\textup{st}}$, uniquely characterized by the property that a homomorphism $\zeta\colon R_\nu^\Box\rightarrow \overline{\Q}_\ell\,$ factors through $R_\nu^{E-\textup{st}}$ if and only if the corresponding representation is semistable after base change to $E$ and has $\ell$-adic Hodge-Tate type $\lambda$.
\end{theorem}
\begin{proof}
This follows from Corollary 2.6.2 of \cite{Kisin2008}.
\end{proof}

 An inertia type is a representation $\theta\colon I_{F_{\tilde{\nu}}} \rightarrow \GL_n(\overline{\Q}_\ell)$ with open kernel. Let $\Theta=\{\theta_1,\ldots,\theta_k\}$ denote a finite set of inertia types.

%
%
\begin{theorem}
\label{semistable_rings}
There exists a quotient of $R_\nu^\Box$, denoted $R_\nu^{\Theta, \lambda}$, uniquely characterized by the property that a homomorphism $\zeta\colon R_\nu^\Box\rightarrow \overline{\Q}_\ell$ factors through $R_\nu^{\Theta,\lambda}$ if and only if the corresponding representation is potentially semistable of type inside $\Theta$, and of $\ell$-adic Hodge-Tate type $\lambda$. Furthermore, $R_\nu^{\Theta,\lambda}$ is reduced, $\ell$-torsion free, and equidimensional of dimension
\begin{equation*}
1+n^2+[F_\nu\colon \Q_\ell]n(n-1)/2.
\end{equation*}
\end{theorem}
\begin{proof}
When $k=1$, the existence and uniqueness is proven in Theorem 2.7.6 of \cite{Kisin2008}, and the dimension formula follows from Theorem 3.3.4 of loc. cit.

When $k>1$, let $E/F_{\nu}$ be a finite extension such that the inertia subgroup $I_E\subset I_{F_{\tilde{\nu}}}$ is contained in $\ker{\theta_i}$ for all $i$. Let  $R_\nu^{E,\lambda}$ be as in the previous theorem. The proof of Theorem 2.7.6 of Kisin constructs each $R_\nu^{\theta_i,\lambda}$ as the quotient of $R_\nu^{E,\lambda}$ corresponding to the union of some components of $\textup{Spec}(R_\nu^{E,\lambda})$. If we denote these unions by $\mathcal{C}^{i,\lambda}$, then it suffices to define $R_\nu^{\Theta, \lambda}$ to be the quotient of $R_\nu^{E,\lambda}$ corresponding to $\cup \mathcal{C}^{i,\lambda}.$
\end{proof}

Let us now consider two local deformation problems at places $\nu\nmid\ell$. The first deformation problem is given by $R_\nu^{\textup{fl}}$, the maximal reduced, $\ell$-torsion free quotient of $R_\nu^\Box$, and we denote it by $\mathcal{D}_\nu^{\textup{fl}}$ (see Section 2.3.4 of \cite{Thorne2016}).

The other deformation problem we will consider is the one corresponding to Taylor-Wiles primes. Let $\nu$ be a finite place of $F^+$, not in $S$ (so $\overline{r}|_{G_{F_{\tilde{v}}}}$ is unramified at $\nu$), which splits in $F$, and  satisfies that
\begin{equation*}
\N\nu \equiv 1 \Mod{\ell}.
\end{equation*}
Since  $\overline{r}|_{G_{F_{\tilde{v}}}}$ is unramified, it can be written as $\overline{\phi_\nu}\oplus\overline{\psi_\nu}$ where $\overline{\psi_\nu}$ is an eigenspace of Frobenius corresponding to an eigenvalue $\alpha_\nu$, on which Frobenius acts semisimply.

\begin{definition}
Denote by $\mathcal{D}_\nu^{\textup{TW}}$ the subfunctor of $\mathcal{D}_\nu^\Box$ which maps a ring $R$ to the set of liftings of $\overline{r}|_{G_{F_{\tilde{v}}}}$ which are $1+\textup{M}_n(\mathfrak{m}_R)$-conjugate to one of the form ${\phi_\nu}\oplus{\psi_\nu}$, where $\phi_\nu$ is unramified and $\psi_\nu$ maps inertia to a subset of the scalar matrices.
\end{definition}

That this functor defines a local deformation problem is the content of Lemma 4.2 of \cite{Th2012}. With this we finish our review of the local deformation theory of Galois representation, and we move onto the global one.

%
%
\label{section global deformation rings}
\begin{definition}
A lifting of $\overline{r}$ to an object $R$ of $\mathcal{C}_\mathcal{O}$ is a continuous homomorphism $r\colon G_{F^+,S}\rightarrow \mathcal{G}_n(R)$ with $r\equiv \overline{r} \Mod{\mathfrak{m}_R}$ and $\nu\circ r = \chi$. Two liftings are said to be equivalent if they are conjugate by an element of $1+\textup{M}_n(\mathfrak{m}_R)$.

Let $T\subseteq S$. By a $T$-framed lifting of $\overline{r}$ to $R$ we mean a tuple $(r;\alpha_\nu)_{\nu\in T}$ where $r$ is a lifting of $\overline{r}$ and $\alpha_\nu\in1+\textup{M}_n(\mathfrak{m}_R)$. Two $T$-framed liftings $(r;\alpha_\nu)_{\nu\in T}$ and $(r';\alpha'_\nu)_{\nu\in T}$ are equivalent if there is an element $\beta\in1+\textup{M}_n(\mathfrak{m}_R)$ such that $r'=\beta r \beta^{-1}$ and $\alpha'_\nu=\beta\alpha_\nu$.
\end{definition}

Given a collection of local deformation problems $\mathcal{D}_\nu$ for $\nu\in S$, we define a global deformation problem to be the tuple
\begin{equation*}
\mathcal{D} = (F/F^+,S,\tilde{S},\mathcal{O},\overline{r},\chi,\{\mathcal{D}_\nu\}_{\nu\in S}).
\end{equation*}
A $T$-framed lifting $(r;\alpha_\nu)_{\nu\in T}$ is of type $\mathcal{D}$ if $r|_{G_{F_{\nu}}}$ is in $\mathcal{D}_\nu$ for all $v\in S$. We denote by $\mathcal{D}^{\Box,T}$ the functor from $\mathcal{C}_\mathcal{O}$ to SETS which maps an $\mathcal{O}$-algebra $R$ to the set of equivalence classes of $T$-framed liftings of $\overline{r}$ of type $\mathcal{D}$. If $T=S$ we denote it by $\mathcal{D}^{\Box}$, and if $T=\varnothing$, by $\mathcal{D}$.

\begin{proposition}
Assume that $\overline{r}|_{G_{F,S}}$ is absolutely irreducible. Then for any $T\subseteq S$ the functor $\mathcal{D}^{\Box,T}$ is represented by an object $R_\mathcal{D}^{\Box,T}$ in $\mathcal{C}_\mathcal{O}$.
\end{proposition}
\begin{proof}
This is Proposition 2.2.9 in \cite{CHT2008}.
\end{proof}
The rings corresponding to the functors $\mathcal{D}^{\Box,T}$, $\mathcal{D}^{\Box}$ and $\mathcal{D}$ defined above are denoted by $R_\mathcal{D}^{\Box,T}$, $R_\mathcal{D}^{\Box}$ and $R_\mathcal{D}$ respectively.

Given a deformation problem $\mathcal{D}$ as above, and a finite set of Taylor-Wiles primes $Q$ disjoint from $S$, we extend $\mathcal{D}$ to a new deformation problem $\mathcal{D}_Q$ by considering Taylor-Wiles deformation at primes in $Q$. More precisely,
\begin{equation*}
\mathcal{D}_Q = (F/F^+,S\cup Q,\tilde{S}\cup\tilde{Q},\mathcal{O},\overline{r},\chi,\{\mathcal{D}_\nu\}_{\nu\in S}\cup\{\mathcal{D}^{TW}_\nu\}_{\nu\in Q}).
\end{equation*}
The previous proposition implies that if $T\subseteq S$, then the deformation problems $\mathcal{D}_Q^{\Box,T}$, $\mathcal{D}_Q^{\Box}$ and $\mathcal{D}_Q$ are representable, and we denote its representing rings by $R_Q^{\Box,T}$, $R_Q^{\Box}$ and $R_Q$ respectively. 

To end this section, let us recall that global deformation rings are quotients of power series rings over local deformation rings, and that this fact is at the heart of the patching argument. Since in our proof of the potentially semistable automorphy lifting theorem (see Theorem 3.3) we make no modifications to the proof of Theorem 6.8 in \cite{Th2012} regarding this aspect, we will make no further comment on this topic, and instead refer the interested reader to Section 4 of \cite{Th2012} (with its corrigenda in Section 7 of \cite{Thorne2016}).

%
%
\label{section_connected_deformations}

The following definition is an adaptation of Definition 3.6 in \cite{Thorne2016} to the potentially semistable case, and it formalizes the idea of when two local representations are very similar. Let $\ell$ and $p$ be prime numbers, and let $F_\nu/\Q_p$ be a finite extension. Let $\lambda$ be an $\ell$-adic Hodge-Tate type, as considered in Section \ref{local deformation rings}. Consider a continuous representation $\rho_1\colon G_{F_\nu}\rightarrow \GL_n(\overline{\Z}_\ell)$ with residual representation $\overline{\rho}_1\colon G_{F_\nu}\rightarrow \GL_n(\overline{\mathds{F}}_\ell)$.

\begin{definition}
\label{def:connected}
In the context above, let $\rho_2\colon G_{F_\nu}\rightarrow \GL_n(\overline{\Z}_\ell)$ be another continuous representations with residual representation $\overline{\rho}_2\colon G_{F_\nu}\rightarrow \GL_n(\overline{\mathds{F}}_\ell)$. Assume that $\overline{\rho}_1\simeq\overline{\rho}_2$. Then $\rho_1$ \emph{connects} with $\rho_2$, denoted by $\rho_1\sim\rho_2$, if
\begin{enumerate}
\item $\ell\neq p$, and $\rho_1$ and $\rho_2$ define points on a common irreducible component of $\textup{Spec}(R^{\textup{fl}}\otimes\overline{\Q}_\ell)$.
\item $\ell=p$, $\rho_1$ and $\rho_2$ are semistable over a finite extension $E/F_\nu$ of the same Hodge-Tate type $\lambda$, and, they define points on a common irreducible component of $\textup{Spec}(R^{E,\lambda}\otimes\overline{\Q}_\ell)$.
\end{enumerate}
\end{definition}

This definition makes sense by Lemma 1.2.2 in \cite{BLGGT}, and it is independent of the choice of isomorphism  $\overline{\rho}_1 \simeq \overline{\rho}_2$.

\subsection{Proof of ALT Semistable: Galois Representations and Hecke Algebras}
\label{section_2}

In this section we recall the main properties of Galois representations attached to maximal ideals of the Hecke algebra, emphasizing their potential semistability. Essentially, this is the content of Theorem 4.1 of \cite{Thorne2016}, which we repeat here for convenience of the reader.

The main theorem is the following, which attaches Galois representations to regular algebraic, cuspidal, polarized automorphic representations (see Section 2 of \cite{BLGGT} for the definitions). Note that these representations are potentially semistable at primes $\nu|\ell$.

\begin{theorem}
\label{automorphic_representations}
Let $(\pi,\chi)$ be a regular algebraic, cuspidal, polarized automorphic representation of $\GL_n(\mathds{A}_F)$. Then there is a continuous semisimple representation
	\begin{equation*}
	r_{\ell,\iota}(\pi)\colon G_F \rightarrow \GL_n(\overline{\Q}_\ell)
	\end{equation*}
and an integer $w$ with the following properties:
	\begin{enumerate}[itemsep=0.5em]
		\item $(r_{\ell,\iota}(\pi), \epsilon_\ell^{1-n}r_{\ell,\iota}(\chi))$ is a totally odd, polarized $\ell$-adic representation.
		
		\item $r_{\ell,\iota}(\pi)$ is de Rham, and if $\tau$ is an embedding $F\hookrightarrow \overline{\Q}_\ell$, then
			\begin{equation*}
			\textup{HT}_\tau(r_{\ell,\iota}(\pi)) = \{\lambda_{\iota\tau,1}+n-1,\lambda_{\iota\tau,2}+n-2,\ldots,\lambda_{\iota\tau,n}\}.
			\end{equation*}
			
		\item If $\nu$ is a finite place of $F$, then
			\begin{equation*}
				\iota \textup{WD}(r_{\ell,\iota}(\pi)|_{G_{F_\nu}})^{\textup{F-ss}} \simeq \textup{rec}(\pi_\nu \otimes | \textup{det} |_\nu^{(1-n)/2}).
			\end{equation*}
			If $\nu\nmid \ell$ these Weil-Deligne representations are pure of weight $\omega$. When $\nu\mid\ell$, then $r_{\ell,\iota}(\pi)$ is potentially semistable at $\nu$, and if $\pi_\nu$ is unramified, then it is crystalline.
	\end{enumerate}
\end{theorem}

(here $\epsilon_\ell$ denotes the $\ell$-adic cyclotomic character). The proof is due to many people. In particular, existence is proved in Theorem 3.2.3 of \cite{Chenevier2013}, and the strong form of local-global compatibility is due to \cite{Caraiani2012} and \cite{Caraiani} (see Theorem 2.1.1 of \cite{BLGGT} for a more exhaustive list of references).

The previous theorem is the main tool for attaching Galois representations to Hecke algebras. Before introducing these Hecke algebras and their representations, let us recall the notation that we will use to speak about automorphic forms on definite unitary groups. The following is a brief summary of the notation used in Section 4 of \cite{Thorne2016}, to which we refer the reader for more details.

Fix a coefficient field $L\subset\overline\Q_\ell$ with ring of integers $\mathcal{O}$. Let us assume that $F/F^+$ is everywhere unramified, that primes $\nu$ of $F^+$ above $\ell$ are split in $F$, and that $n\left[F^+\colon\Q\right]\equiv0\Mod4$. Then there exists a unitary group over $F^+$, denoted $G$, which is split by $F$. This group can be extended to an affine group scheme over $\mathcal{O}_{F^+}$, also denoted $G$, which satisfies the following property: for any place $\nu$ of $F^+$ split as $\nu=\omega\omega^c$ in $F$, there is an isomorphism
\begin{equation*}
\phi_\omega\colon G(\mathcal{O}_{F^+_\nu})\simeq\GL_n(\mathcal{O}_{F_\omega}).
\end{equation*}

Let us denote by $S_\lambda(U,A)$ the space of automorphic forms of weight $\lambda$ (an $\ell$-adic Hodge-Tate type which determines the highest weights of algebraic representations of $\GL_n$), level $U$ (an open compact subgroup $\prod_\nu U_\nu\subset G(\mathds{A}_{F^+}^\infty)$), and with coefficients in $A$ (an $\mathcal{O}$-algebra). Define
\begin{equation*}
\mathcal{A}_\lambda = \varinjlim_U S_\lambda(U,\overline\Q_\ell).
\end{equation*}
Then each irreducible submodule $\pi \subset \mathcal{A}_\lambda$ corresponds to the finite part of an automorphic representation of $G(\mathds{A}_{F^+})$, denoted by  $\otimes_\nu \pi_\nu$, and with an infinite part that can be described explicitly in terms of $\lambda$.

Denote by $\mathds{T}_\lambda^T(U,\mathcal{O})$ the Hecke algebra of weight $\lambda$, level $U$, with coefficients in $\mathcal{O}$, and outside of $T$ (a finite set of places of $F^+$). This is the $\mathcal{O}$-subalgebra of $\End_{\mathcal{O}}(S_\lambda(U,\mathcal{O}))$ generated by the Hecke operators $T_\omega^i$, where $i=1,\ldots,n$ and $\omega$ runs throught the finite places of $F$ which are split over $F$ and prime to $T$.

Lastly, denote by $\delta_{F/F^+}$ the nontrivial character of the quadratic extension $F/F^+$.

%
%
\begin{proposition}
\label{representations_submodules}
Let $U$ be an open compact subgroup of $G(\mathds{A}_{F^+}^\infty)$, and let $T$ be a finite set of finite places of $F^+$ such that $\mathds{T}^T_\lambda(U,\mathcal{O})$ is defined. If $\pi$ is an irreducible submodule of $\mathcal{A}_\lambda$ such that $\pi^U\neq 0$, then there exists a unique continuous semi-simple representation
\begin{equation*}
r_{\ell,\iota}(\pi)\colon G_F \rightarrow \GL_n(\overline{\Q}_\ell)
\end{equation*}
satisfying the following properties:
\begin{enumerate}[itemsep=0.5em]
\item  $r_{\ell,\iota}(\pi)^c \simeq r_{\ell,\iota}(\pi)^\vee \epsilon^{1-n}$,
\item If $r_{\ell,\iota}(\pi)$ is irreducible then $(r_\ell(\pi),\epsilon^{1-n}\delta_{F/F^+}^n)$ is polarized, and in this case, if $\nu$ is a finite place of $F^+$ which splits in $F$ as $\nu=\omega\omega^c$, then $\pi_\nu\circ \phi_\omega^{-1}$ is generic.
\item $r_{\ell,\iota}(\pi)$ is unramified above finite inert places $\nu$ of $F^+$ for which $U_\nu$ is a hyperspecial maximal compact subgroup.
\item $r_{\ell,\iota}(\pi)$ is de Rham, and if $\tau$ is an embedding $F\hookrightarrow \Q_\ell$, then
	\begin{equation*}
	\textup{HT}_\tau(r_{\ell,\iota}(\pi)) = \{\lambda_{\iota\tau,1}+n-1,\lambda_{\iota\tau,2}+n-2,\ldots,\lambda_{\iota\tau,n}\}.
	\end{equation*}
\item If $\nu$ is a finite place of $F^+$ split in $F$ as $\nu=\omega\omega^c$, then there is an isomorphism
	\begin{equation*}
	\iota\textup{WD}(r_{\ell,\iota}(\pi)|_{G_{F_\omega}})^{\textup{F-ss}} \simeq \textup{rec}_{F_\omega}(\pi_\nu \circ \phi_\omega^{-1} \otimes | \textup{det} |_\nu^{(1-n)/2}).
	\end{equation*}
	In particular, if $\nu|\ell$ then $r_{\ell,\iota}(\pi)$ is potentially semistable, and if $\pi_\nu$ is unramified then $r_{\ell,\iota}(\pi)$ is crystalline.
\end{enumerate}
\end{proposition}

\begin{proof}
The second assertion follows from Theorem 3.2 in \cite{Thorne2016}. The rest follows from Theorem 2.3 in \cite{Guerberoff2011} together with local-global compatibility in its strong form, as proven in \cite{Caraiani2012} and \cite{Caraiani}. 
\end{proof}

%
%
\begin{proposition}
\label{residual_hecke}
Let $\mathfrak{m}\subset\mathds{T}^T_\lambda(U,\mathcal{O})$ be a maximal ideal. Then, there exists a unique continuous semisimple representation
	\begin{equation*}
	\overline{r}_\mathfrak{m} \colon G_F \rightarrow \GL_n(\mathds{T}^T_\lambda(U,\mathcal{O})/\mathfrak{m}),
	\end{equation*}
satisfying the following properties:
\begin{enumerate}[itemsep=0.5em]
\item $\overline{r}_\mathfrak{m}^c \simeq \overline{r}_\mathfrak{m}^\vee \epsilon^{1-n}$.
\item $\overline{r}_\mathfrak{m}$ is unramified at all but finitely many places.
\item $\overline{r}_\mathfrak{m}$ is unramified above finite inert places $\nu$ of $F^+$ for which $U_\nu$ is a hyperspecial maximal compact subgroup.
\item If $\nu\not\in T$ is a finite place of $F^+$ then $\overline{r}_\mathfrak{m}$ is unramified at places above $\nu$. Furthermore, if $\nu$ splits as $\omega\omega^c$ in $F$, then $\overline{r}_\mathfrak{m}(\textup{Frob}_\omega)$ has characteristic polynomial
	\begin{equation*}
	X^n + \ldots + (-1)^j(\mathbb{N}\omega)^{j(j-1)/2}T_\omega^{j}X^{n-j} + \ldots + (-1)^n(\mathbb{N}\omega)^{n(n-1)/2}T_\omega^{n}.
	\end{equation*}
\item If $\overline{r}_\mathfrak{m}$ is irreducible then the pair $(\overline{r}_\mathfrak{m},\epsilon^{1-n}\delta_{F/F^+}^n)$ is polarized.
\end{enumerate}
\end{proposition}
\begin{proof}
The proof is the same as that of Proposition 3.1 in \cite{Guerberoff2011}. It essentially frollows from reducing modulo $\mathfrak{m}$ the representation constructed in Proposition \ref{representations_submodules} above.
\end{proof}

%
%
\begin{proposition}
\label{lifting_maximal_hecke}
Let $\mathds{T}^T_\lambda(U,\mathcal{O})$ and $\mathfrak{m}$ be as before, and suppose further that $\overline{r}_\mathfrak{m}$ is absolutely irreducible. Then $\overline{r}_\mathfrak{m}$ admits an extension to a continuous homomorphism $\overline{r}_\mathfrak{m} \colon G_{F^+} \rightarrow \mathcal{G}_n(\mathds{T}^T_\lambda(U,\mathcal{O})/\mathfrak{m})$.
Furthermore, fixed any such extension, there exists a unique continuous lifting
\begin{equation*}
r_\mathfrak{m} \colon G_{F^+} \rightarrow \mathcal{G}_n(\mathds{T}^T_\lambda(U,\mathcal{O})_\mathfrak{m})
\end{equation*}
satisfying the following properties:
\begin{enumerate}[itemsep=0.5em]
\item $\nu\circ r_\mathfrak{m}= \epsilon^{1-n}\delta^{n}_{F/F^+}$.
\item $r_\mathfrak{m}$ is unramified at almost all places.
\item $r_\mathfrak{m}$ is unramified above finite inert places $\nu$ of $F^+$ for which $U_\nu$ is a hyperspecial maximal compact subgroup.
\item If $\nu$ is a finite place of $F^+$ which is not in $T$, then ${r}_\mathfrak{m}$ is unramified at places above $\nu$. Furthermore, if $\nu$ splits as $\omega\omega^c$ in $F$, then ${r}_\mathfrak{m}(\textup{Frob}_\omega)$ has characteristic polynomial
	\begin{equation*}
	X^n + \ldots + (-1)^j(\mathbb{N}\omega)^{j(j-1)/2}T_\omega^{j}X^{n-j} + \ldots + (-1)^n(\mathbb{N}\omega)^{n(n-1)/2}T_\omega^{n}.
	\end{equation*}
\end{enumerate}
\end{proposition}
\begin{proof}
The proof is the same as that of Theorem 4.1 in \cite{Thorne2016}. Let us recall the construction for future reference. For each minimal prime $\mathfrak{p}\subset\mathds{T}^T_\lambda(U,\mathcal{O})_\mathfrak{m}$ we can find an irreducible subrepresentation 
$\pi\subset\mathcal{A}_\lambda$ such that $\pi^U\neq0$. The representation $r_{\ell,\iota}(\pi)$ given by Proposition \ref{representations_submodules} can be conjugated to have values in the ring of integers of some finite extension of the field of fractions of $\mathds{T}^T_\lambda(U,\mathcal{O})/\mathfrak{p}$. Moreover, this can be done in such a way that after reduction by the maximal ideal we recover $\overline{r}_\mathfrak{m}$. Let $r_\mathfrak{p}$ denote this conjugation of $r_{\ell,\iota}(\pi)$, and let $\mathcal{O}_\mathfrak{p}$ denote the ring of integers where it has values. Then, consider the product of these representations (see Theorem 4.1 of \cite{Thorne2016} for a precise definition)
\begin{equation*}
\times_\mathfrak{p}\; r_\mathfrak{p} \colon G_{F^+} \rightarrow \mathcal{G}_n(\times_\mathfrak{p}\; \mathcal{O}_\mathfrak{p}).
\end{equation*}
Since $\mathds{T}^T_\lambda(U,\mathcal{O})_\mathfrak{m}$ is reduced, it injects in the product of the rings $\mathcal{O}_\mathfrak{p}$, and after applying Lemma 2.4 of loc. cit., we can conjugate this representation to have values in $\mathds{T}^T_\lambda(U,\mathcal{O})_\mathfrak{m}$.
\end{proof}

\subsection{Proof of ALT Semistable: Patching}
\label{section_3}

In this section we extend the minimal automorphy lifting theorems of \cite{Th2012} to the semistable case. The proofs are the same as in \cite{Th2012} (taking into account the corrections and extensions to the case $\ell|n$ in \cite{Thorne2016}). There is just one small adjustment that needs to be done to Proposition 6.8 of \cite{Th2012} (or to Theorem 4.2 of \cite{Thorne2016}), in order to account for the fact that the semistable deformation rings need not be smooth. This is explained in the proof of Proposition \ref{patching} below.

Let us set up the context. Let $T$ be a finite set of places of $F^+$ which are split in $F$, and assume that $T$ contains the places above $\ell$. Let $U$ be an open compact subgroup of $G(\mathds{A}_{F^+}^\infty)$ such that
\begin{enumerate}
\item $U_\nu$ is arbitrary if $\nu\in T$ (in particular, if $\nu|\ell$),
\item $U_\nu = G(\mathcal{O}_{F^+_\nu})$ if $\nu\notin T$ is split in $F$,
\item$ U_\nu$ is a hyperspecial maximal compact of $G(F_\nu^+)$ if $\nu$ is inert in $F$.
\end{enumerate}
Let $L$ be a finite extension of ${\Q}_\ell$ with ring of integers $\mathcal{O}$, maximal ideal $\lambda$, and residue field $\mathds{F}$. Choose an $\ell$-adic Hodge-Tate type $\lambda$. In this setting $\mathds{T}^T_\lambda(U,\mathcal{O})$ is well defined, and we can consider a maximal ideal $\mathfrak{m}$ inside of it. 

In order to define our global deformation problem, we recall that if $\pi$ is an irreducible submodule of $\mathcal{A}_\lambda$ corresponding to some minimal prime in $\mathds{T}^T_\lambda(U,\mathcal{O})$, then $\pi^U\neq0$. In particular, $\pi_\nu^{U_\nu}\neq0$ at $\nu|\ell$, so the work of Bernstein-Zelevinsky \cite{Bernstein1977} and Zelevinsky \cite{Zelevinsky1980}, together with the strong form of local-global compatibility, proven by Caraiani in \cite{Caraiani}, imply that there are a finite number of possibilities for the inertia type of $r_{\ell,\iota}(\pi)$ at $\nu|\ell$. Let $\Theta$ denote the finite set containing these possible inertia types.

Define the following global deformation problem, whose local deformation problems $\mathcal{D}_\nu$ at $\nu\in T$ correspond to $R_\nu^\textup{fl}$ if $\nu\nmid \ell$, and to $R_\nu^{\Theta, \lambda}$ if $\nu|\ell$. That is,
\begin{equation*}
\mathcal{D} = (F/F^+,T,\tilde{T},\mathcal{O},\overline{r}_\mathfrak{m},\epsilon^{1-n}\delta^{n}_{F/F^+},\{\mathcal{D}_\nu\}_{\nu\in S}),
\end{equation*}
where $\epsilon$ denotes the $\ell$-adic cyclotomic character, and $\delta_{F/F^+}$ denotes the nontrivial character of the quadratic extension $F/F^+$.

\begin{lemma}
Assume that $\overline{r}_\mathfrak{m}$ is absolutely irreducible. Then, the lifting $r_\mathfrak{m}$ is of type $\mathcal{D}$, and so, it comes from a surjective homomorphism $R_\mathcal{D}\rightarrow \mathds{T}^T_\lambda(U,\mathcal{O})_\mathfrak{m}$.
\end{lemma}
\begin{proof}
This follows from the construction of $r_\mathfrak{m}$ in Proposition \ref{lifting_maximal_hecke} and the choice of $\Theta$. Indeed, $r_\mathfrak{m}$ consists of a product of representations given by Proposition \ref{representations_submodules}, which are clearly of type $\mathcal{D}$. The claim now follows from item 5 in the definition of local deformation problems (i.e. Definition \ref{local deformation problems}).
\end{proof}

The patching argument also requires to have maps $R_Q\rightarrow \mathds{T}^T_\lambda(U_1(Q),\mathcal{O})_\mathfrak{m}$ for sets of Taylor-Wiles primes $Q$ (see the proof of Theorem 6.8 of \cite{Th2012} for the definition of $U_1(Q)$ and the details). That these exists follows from the argument in loc. cit, noting that since $U_1(Q)_\nu=U_\nu$ at $\nu|\ell$ then we still have the same possibility for the inertia types at $\nu|\ell$.

The next proposition is an extension of Theorem 6.8 in \cite{Th2012} (or Proposition 7.2 in \cite{Thorne2016}) to the potentially semistable case. We remark that the notion of adequate subgroup that we use is that of Definition 2.20 of \cite{Thorne2016} (which includes the case $\ell|n)$.

%
%
\begin{proposition}
\label{patching}
Let $r\colon G_{F^+}\rightarrow\mathcal{G}_n(\mathcal{O})$ be another lifting of $\overline{r}_\mathfrak{m}$ of type $\mathcal{D}$. Assume $\zeta_\ell\not\in F$,  and,
\begin{enumerate}
\item $\overline{r}_\mathfrak{m}$ maps $G_{F(\zeta_\ell)}$ to an adequate subgroup.
\item there exists $f\colon \mathds{T}^T_\lambda(U,\mathcal{O}) \rightarrow \mathcal{O}$ such that $f\circ r_\mathfrak{m} \sim r$ at $\nu\in T$ (i.e. they are connected).
\end{enumerate}
Then, there exists another homomorphism $g\colon \mathds{T}^T_\lambda(U,\mathcal{O}) \rightarrow \mathcal{O}$ such that $g\circ r_\mathfrak{m}$ and $r$ are conjugate to each other in $\GL_n(\mathcal{O})$.
 \end{proposition}
\begin{proof}
The proof is the same as in Theorem 6.8 in \cite{Th2012}, with a small adjustment arising from the fact that semistable deformation rings need not be smooth. Indeed, at the end of the argument, we need to consider the irreducible components of the rings $R^{\textup{fl}}_\nu$ and $R_\nu^{\Theta, \lambda}$ that contain the restriction of $f\circ r_\mathfrak{m}$ to ${G_{F_\omega}}$ (where $\omega|\nu$ is a place of $F$). For the rest of the argument to work, we need to show that these components are unique, and this follows from the strong form of local-global compatibility (see \cite{Caraiani}) and Lemma 1.3.2 in \cite{BLGGT} (cf. with the proof in \cite{Th2012}, where uniqueness at primes $\nu|\ell$ followed from smoothness of the crystalline deformation rings).
\end{proof}

The semistable automorphy lifting theorem follows  from this proposition. Let us restate it as in Section \ref{sec:statementsALT} with the numbering it appeared there.

%
%
\begin{reptheorem}{alt_semistable}
Let $n\in\Z_{\geq1}$, and fix an odd prime $\ell$. Let $F$ be an imaginary CM field with totally real subfield $F^+$, and assume that $\zeta_\ell\notin F$. Suppose that $(r,\mu)$ is an $n$-dimensional, algebraic, polarized $\ell$-adic representation of $G_F$ satisfying the following properties:
\begin{enumerate}
\item $\overline{r}$ is irreducible and maps $G_{F(\zeta_\ell)}$ to an adequate subgroup (adequate as in Definition 2.20 of \cite{Thorne2016}, so in particular it could be that $\ell\mid n$).
\item $(\overline{r},\overline{\mu})$ is automorphic, arising from a regular algebraic, cuspidal, polarized automorphic representation $(\pi,\chi)$ such that at all places $\nu$ of $F$
\begin{equation*}
r|_{G_{F_\nu}} \sim r_\ell(\pi)|_{G_{F_\nu}}
\end{equation*}
(that is, the representations are connected, as defined in Definition \ref{def:connected}). In particular, locally at $\nu\mid\ell$ both representations are potentially semistable with equal Hodge-Tate weights.
\end{enumerate}
Then $(r,\mu)$ is automorphic.
\end{reptheorem}
\begin{proof}
This is deduced from the previous proposition in the same way Theorem 7.1 is deduced from Theorem 6.8 in \cite{Thorne2016} (cf. with Corollary 7.3 of \cite{Thorne2016}).
\end{proof}


\section{2-dimensional Chain}\label{sec:2-dim Chain}

Let $f\in S_k(N)$ be a modular form without complex multiplication, which is an eigenform for all Hecke operators $T_{q}$ with $q\nmid N$.
Recall that for each prime $p$, we have fixed a prime $\mathcal{P}\subset \overline{\mathbb{Z}}$ (equivalently, an isomorphism $\iota_p\colon\overline{\mathbb{Q}}_p\simeq\mathbb{C}$), hence a prime $\mathfrak{p}\vert p$ of the coefficient field $\mathbb{Q}_f$ of $f$. Let
$\rho_{\mathfrak{p}}(f):G_{\mathbb{Q}}\rightarrow \GL_2(\overline{\mathbb{Q}}_p)$
be the irreducible Galois representation of the absolute Galois group $G_{\mathbb{Q}}$ associated to $f$ and $\mathfrak{p}$ by classical work of
Deligne. Since our choice of prime $\mathcal{P}$ is fixed, we will denote this representation by $\rho_{p}(f)$.
By Ribet's theorem, the image of the residual representation
$\overline{\rho}_{p}(f)$ is large (i.e.~contains $\SL_2(\mathbb{F}_p)$) for $p$ sufficiently large.
We say the prime $p$ is \emph{exceptional} for $f$ if the image of $\overline{\rho}_{p}(f)$ is not large.
If $f$, $f'$ are two modular eigenforms, we will say that they are \emph{congruent modulo $p$} if the residual representations
$\overline{\rho}_p(f)$ and $\overline{\rho}_p(f')$ coincide.

Fix a modular form $f\in S_k(1)$. The goal of this section is to describe a chain connecting $f$ with a CM modular form $g$. By a chain, we mean 
a list of modular forms 
\begin{equation*}\mathcal{C}:=\{f=f_1, f_2 \dots, f_h=g\}\end{equation*}
 such that, if $f_i$ and
$f_{i+1}$ are two consecutive terms, then there exists a prime $p_i$ such that $f_i$ is congruent to either $f_{i+1}$,
a twist of $f_{i+1}$, or a Galois conjugate of $f_{i+1}$ modulo $p_i$. In addition, we want the chain to be "safe", meaning that at each congruence in the chain, a suitable ALT can be applied.
More precisely, we want that  the
following conditions hold at all congruences of the chain:

\begin{itemize}

 \item Local conditions at $p_i$: $\rho_{p_i}(f_i)$ and $\rho_{p_i}(f_{i+1})$ are either both potentially
 diagonalisable at  $p_i$, or both ordinary at  $p_i$.

 \item Residual image: the restriction of $\overline{\rho}_{p_i}(f_i)$ to $G_{ \mathbb{Q}(\zeta_{p_i}) }$ has adequate$+$ image. 	
\end{itemize}

In what follows, we present the construction of a particular "safe" chain, linking a given modular form $f\in S_k(1)$ with the CM form $g\in S_2(27)$. For the rest of the paper, this particular chain will be called the  \emph{skeleton chain}. The construction of the skeleton chain consists of 15 steps.
The first eleven steps are taken from \cite{Di-Symm5}, whereas the last four steps have been collected from  \cite{Di-mfold}. Observe that there is a Corrigenda to \cite{Di-Symm5}, where Step 7 of the chain is corrected and Step 3 and Step 4 are replaced by two steps called Step 3NEW and Step 4NEW. Following a suggestion of the referee, for a greater clarity and in order to keep this paper self-contained, we will present step by step a description of the skeleton chain, including in particular full details for the three steps addressed in the Corrigenda. 

It is important to remark, for future reference in this paper, that the congruences in this chain can be classified into three types, according to the local data of the two modular forms involved.

\begin{itemize}
\item \emph{Congruence of type A}: Both  $f_i$ and $f_{i+1}$ are potentially
diagonalisable at $p_i$, of different weights,
and connected at all $p\not=p_i$.
 \item \emph{Congruence of type B}: Both $f_i$, $f_{i+1}$ are ordinary at $p_i$, of different weights, and are connected at all $p\not=p_i$.
 \item \emph{Congruence of type C}: Both $f_i$, $f_{i+1}$ are potentially diagonalisable at $p_i$, but they are of different type at some $p\not=p_i$ (in particular, cases of level raising or lowering at $p$).
\end{itemize}

\begin{rem}\label{rem:CongruenceTypes} 
We will see that,  except at Step 12, whenever there is a congruence of type B in the skeleton chain, say between $f_{i}$ and $f_{i+1}$, we have that $f_{i}$ is a weight $2$ modular form, which
 is Steinberg at the prime $p_{i}$. Thus the representation $\rho_{p_i}(f_{i})$ is ordinary, but not potentially crystalline locally at $p_{i}$.
 In Section 6.3, we
 will need to build on the previous congruence between $f_{i-1}$ and $f_{i}$ modulo $p_{i-1}$. One can check that, with the above mentioned exception,
 whenever the congruence between
 $f_{i}$ and $f_{i+1}$ is of type B, the preceding congruence between $f_{i-1}$ and $f_i$ in the chain is a congruence of type C, 
 where the type is changed at $p_i$. 
\end{rem}

Almost all congruences in the skeleton chain are of type A. Whenever a type B congruence appears, we will draw a diagram to illustrate Remark \ref{rem:CongruenceTypes}. In Step 12 we also draw a diagram, although in this case the congruence of type B is not preceeded by a type C congruence (we will need a separate section, namely Section 6.4, to deal with this case). For completeness, we will also include a diagram in cases where an isolated type C congruence appears.

\begin{rem}\label{rem:image-chain} We will see in the description that, at all congruences, except for Step 15, we have that 
$\overline{\rho}_{p_i}(f_i)(G_{ \mathbb{Q} })$
contains $\SL_2(\mathbb{F}_{{p_i}^r})$, where either $p_i\geq 7$ or $p_i=3$,  with $r> 2$ in the latter case. By elementary group theory, the same holds for the restriction of $\overline{\rho}_{p_i}(f_i)$ to $\mathbb{Q}(\zeta_{p_i})$.
Thus adequacy$+$ of
$\overline{\rho}_{p_i}(f_i)(G_{ \mathbb{Q}(\zeta_{p_i}) })$ holds by Lemma \ref{lem:SL2}.\end{rem}

We proceed now to the detailed description of the skeleton chain.

\subsection{Step 1 (cf.~\cite[Section 3.1]{Di-Symm5})} The purpose of this step is to introduce a good dihedral prime in the level of the given modular form, and it will consist of three congruences. 

We start with $f:=f_1$. First, we pick a prime $r>k$ such that the image of $\overline{\rho}_r(f_1)$ is large. 
We have a mod $r$ congruence with a modular form $f_2$ of weight $2$ with level $r$ and some  nebentypus. 
Observe that this congruence is of type A because one of the representations is Fontaine-Laffaille and the 
other one is potentially Barsotti-Tate, thus they are both potentially diagonalisable (cf.~Remark \ref{rem:3.11}). 

Next, a suitable bound $B_0$ is fixed as in \cite[Section 3.1]{Di-Symm5} (denoted $B$ in loc.~cit.). In particular, $B_0$ is greater than $k$, $2r$ and $53$.  Once $B_0$ is fixed, we pick primes $t_0$ and $q_0$ (these are denoted $t$ and $q$ in loc.~cit.) satisfying the conditions specified in Section 3.1 of loc.~cit. (in particular, the residual image modulo $t_0$ is large). Modulo $t_0$ we produce a congruence with a modular form $f_3$ of weight $2$ and level $rq_0^2$, having $q_0$ as a good dihedral prime with respect to the bound $B_0$. This congruence is of type C, since both representations are Barsotti-Tate and thus potentially diagonalisable (cf. Remark \ref{rem:3.11}), and the type has changed at the prime $q_0$. 

From now on, up to Step 4NEW (included), we are going to work in characteristics $7\leq p<B_0$, and this implies that the residual images will be large because of the good dihedral prime $q_0$ in the level.

To conclude this step, we move back to characteristic $r$ and obtain a congruence, up to twist, with a form $f_4$ of weight $k'$ greater than $2$ and smaller than $r$, and level $q_0^2$. This is again a type A congruence.

 \small

  \begin{equation*}
  \hskip-2cm\xymatrix@C-1pc{f_{1}\in S_k(1) \ar@/_2pc/ @{-}[r]_{\text{mod } r\atop \raisebox{-3ex}{\text{\small TYPE A}}} & f_{2}\in S_2(\Gamma_1(r))
  \ar@/_2pc/ @{-}[r]_{\text{mod } t_0\atop \raisebox{-3ex}{\text{\small TYPE C}}} &
  f_{3}\in S_{2}(\Gamma_1(q_0^2r))\ar@/_2pc/ @{-}[r]_{\text{mod } r\atop \raisebox{-3ex}{\text{\small TYPE A}}} & f_{4}\in S_{k'}(q_0^2) \\}
 \end{equation*}

\normalsize

\subsection{Step 2 (cf.~\cite[Section 3.2]{Di-Symm5})}
In this step we apply the method of Weight Reduction via Galois Conjugation (WRGC) described in \cite{Dieulefait2009}. It consists of several iterations of a combination of congruences and Galois conjugation, specially devised to reduce the weight. 

We describe the generic step that will be iterated. We start with a form of weight $2<k<r$ and level $q_0^2$. If $k\leq 14$, we are done and proceed to Step 3NEW below. Otherwise, pick $p$ to be the smallest prime larger than $k$ (except when $k=32$, where we take $p=43$). We begin by taking a mod $p$ congruence with a weight $2$ modular form, with level $p q_0^2$ and nebentypus at $p$. The size of $B_0$ with respect to $r$ ensures that $p<B_0$, thus the residual image is large. This congruence is of type A. 

Here comes the key step of the WRGC method: it is shown in \cite[Section 3.2]{Di-Symm5} that, for a suitable Galois conjugation, one can pass to a conjugated form, also of weight two and level $pq_0^2$, and conjugated nebentypus, such that when reducing modulo $p$, the corresponding Serre weight (up to twist) is smaller than $k$ and greater than $2$. Thus we have a congruence, up to twist, with a form of level $q_0^2$ and weight $k'$ satisfying $2<k'<k$. The congruence is of type A, and again $p<B_0$ guarantees that the residual image is large.

After this iterative process, we end up with a modular form of weight $2<k\leq 14$, level $q_0^2$ and trivial nebentypus.

\subsection{Step 3NEW (cf.~Corrigenda to \cite{Di-Symm5} in Section 4 (b) of \cite{BD2023})}

This step is actually a preliminary step, which is required before 
the introduction of the Micro-Good Dihedral prime that will occur in Step 4NEW below. The goal is to adjust the weight to 16, 
and the level to $q_0^2$. Throughout all this step, residual images in all congruences are going to be large because of the good dihedral prime $q_0$. 

We start with a modular form $f_i$ of weight $2<k\leq 14$, level $q_0^2$ and trivial nebentypus. We perform a congruence mod 47, with a modular form $f_{i+1}$ of weight $2$,  level $47q_0^2$  and nebentypus $\mu=\omega^{k-2}$, where $\omega$ is the Teichm\"uller lift of the cyclotomic character (cf. \cite{Khare-Wintenberger-I} for the existence and modularity of this lift). The nebentypus $\mu$ has conductor $47$ and order 23. Note that this congruence is of type A. 

If we reduce (modulo $23$) the local $23$-adic Galois representation attached to $f_{i+1}$, the residual representation will be either unramified at $47$ or will have unipotent ramification at $47$. Accordingly, we divide into two cases:

\begin{itemize}

\item Case 1: There is unipotent ramification at $47$. In this case, we take a minimal modular lift, corresponding to a weight 2 modular form $f_{i+2}$ of level $47q_0^2$ that is Steinberg locally at 47.  We have a type C congruence, since we have changed the local type at $47$. 

Next, we move to characteristic 47, and consider $\overline{\rho}_{47}(f_{i+2})$. By considering the Hida family deforming this residual representation, we have a congruence with a modular form $f_{i+3}$ of weight $48$ and level $q_0^2$. Both forms in this congruence are ordinary at $47$, thus this is a type B congruence. 

 We draw a diagram of the congruences performed in Step 3NEW so far for this case. 

 \small

  \begin{equation*}
  \hskip-2cm\xymatrix@C-1pc{f_{i}\in S_k(q_0^2) \ar@/_2pc/ @{-}[r]_{\text{mod } 47\atop \raisebox{-3ex}{\text{\small TYPE A}}} & 
  f_{i+1}\in S_2(\Gamma_1(47\cdot q_0^2))
  \ar@/_2pc/ @{-}[r]_{\text{mod } 23\atop \raisebox{-3ex}{\text{\small TYPE C}}} &
  f_{i+2}\in S_{2}(47\cdot q_0^2)\ar@/_2pc/ @{-}[r]_{\text{mod } 47\atop \raisebox{-3ex}{\text{\small TYPE B}}} & f_{i+3}\in S_{48}(q_0^2) \\}
 \end{equation*}

\normalsize

\item Case 2: There is no ramification at $47$. In this case, if we take a minimal lift, it corresponds to a modular form $f_{i+2}$ of weight $2$ and  level $q_0^2$. This is a type C congruence (we are doing level lowering at $47$).

It is well known (and follows from multiplying by the Hasse invariant $E_{47-1}$) that there is a mod $47$ congruence between $f_{i+2}$  and a newform $f_{i+3}$ of level $q_0^2$ and weight $48$. Let us show that this is a type A congruence. The input form is Barsotti-Tate, hence potentially diagonalisable locally at $47$. For the other form, we consider two cases: if the residual representation locally at $47$ is described by powers of the cyclotomic character, then we are in the ordinary case (residually ordinary implies ordinary in this situation), and we are again moving in a Hida family, but in this case both representations are crystalline and ordinary, thus potentially diagonalisable. In the complementary case, the supersingular case, it is also true that the output representation is potentially diagonalisable. To show this, we rely on Proposition 3.13 in \cite{Kisin_Durham}, where it is proved that, in this case, the universal ring of local deformations which are crystalline of weight $p+1=48$ is a domain, thus it is irreducible, plus the fact that can easily be checked (and is a particular case of Lemma 4.1.19 of \cite{BLGG}) that in the supersingular weight $2$ case we can construct a local lift of weight $48$ which is induced from a crystalline character of the unramified quadratic extension of $\mathbb{Q}_{47}$. Thus, any crystalline deformation of this weight has to be connected to this particular one, therefore it is potentially diagonalisable.

The diagram corresponding to the congruences performed in Step 3NEW in Case 2 is the following:

 \small

  \begin{equation*}
  \hskip-2cm\xymatrix@C-1pc{f_{i}\in S_k(q_0^2) \ar@/_2pc/ @{-}[r]_{\text{mod } 47\atop \raisebox{-3ex}{\text{\small TYPE A}}} & f_{i+1}\in S_2(\Gamma_1(47q_0^2))
  \ar@/_2pc/ @{-}[r]_{\text{mod } 23\atop \raisebox{-3ex}{\text{\small TYPE C}}} &
  f_{i+2}\in S_{2}(q_0^2)\ar@/_2pc/ @{-}[r]_{\text{mod } 47\atop \raisebox{-3ex}{\text{\small TYPE A}}} & f_{i+3}\in S_{48}(q_0^2) \\}
 \end{equation*}

\normalsize

\end{itemize}

In both cases we end up with a modular form $f_{i+3}$ of weight $48$ and level $q_0^2$, so we can treat them in a joint way. Our aim is to reduce the weight from $48$ to $16$, and this will be performed in two steps.

First, we move to $p=53$, and modulo $p$ we have a congruence with a modular form $f_{i+4}$ of weight $2$  and level $53q_0^2$ with nebentypus $\omega^{46}$ of conductor $53$ and order $26$. This is a type A congruence, and the residual image is large because of the good dihedral prime $q_0$ and the fact that $53<B_0$. 

Following a method of Khare, we are going to work modulo the auxiliary prime $13$ (which divides the order of $\omega^{46}$) to modify the nebentypus. We reduce modulo $13$ and we take a minimal lift, which corresponds to a modular form $f_{i+5}$ of weight $2$ and level $53q_0^2$ (for existence and modularity of this lift, cf.~\cite{Khare-Wintenberger-I}). The character giving ramification at $53$ in this minimal lift, in other words, the nebentypus of the corresponding modular form $f_{i+5}$,  is $\omega^{26}$, a quadratic character. This is a type C congruence, since the local type at $53$ was modified, and again the residual image is large.

Next we reduce again modulo $53$. Here, the possible values for the residual Serre's weight (up to twist) are (according to \cite{Savitt_Corrigenda})  $26+2=28$ or $53+3-28=28$. By taking a modular minimal lift, we obtain a congruence, up to twist,  with a modular form $f_{i+6}$ of weight $28$ and level $q_0^2$. This is a type A congruence, and the residual image is large.

 \small

  \begin{equation*}
  \hskip-2cm\xymatrix@C-1pc{f_{i+3}\in S_{48}(q_0^2) \ar@/_2pc/ @{-}[r]_{\text{mod } 53\atop \raisebox{-3ex}{\text{\small TYPE A}}} & f_{i+4}\in S_2(\Gamma_1(53\cdot q_0^2))
  \ar@/_2pc/ @{-}[r]_{\text{mod } 13\atop \raisebox{-3ex}{\text{\small TYPE C}}} &
  f_{i+5}\in S_{2}(\Gamma_1(53\cdot q_0^2))\ar@/_2pc/ @{-}[r]_{\text{mod } 53\atop \raisebox{-3ex}{\text{\small TYPE A}}} & f_{i+6}\in S_{28}( q_0^2)\\}
 \end{equation*}

\normalsize

To reduce the weight further, we repeat the above steps, using different auxiliary primes. First, take $p=29$, and consider a modular weight two lift of level $29q_0^2$, and nebentypus $\omega^{26}$, which is a character of conductor $29$ and order $14$. This is a type A congruence, and the residual image is large.

Next, we reduce modulo $7$ and take a minimal modular lift corresponding to a modular form $f_{i+8}$ of weight $2$ and level $29 q_0^2$, with nebentypus $\omega^{14}$, a quadratic character. This is  a type C congruence, since the local type at $29$ was modified, and again the residual image is large.

Finally, we reduce modulo $29$, and applying \cite{Savitt_Corrigenda} we see that the residual Serre weight, up to twist, is $16$. Thus we get a congruence with a modular form $f_{i+9}$ of weight $16$ and level $q_0^2$. This is a type A congruence, and the residual image is large. 

 \small

  \begin{equation*}
  \hskip-2cm\xymatrix@C-1pc{f_{i+6}\in S_{28}(q_0^2) \ar@/_2pc/ @{-}[r]_{\text{mod } 29\atop \raisebox{-3ex}{\text{\small TYPE A}}} & f_{i+7}\in S_2(\Gamma_1(29\cdot q_0^2))
  \ar@/_2pc/ @{-}[r]_{\text{mod } 7\atop \raisebox{-3ex}{\text{\small TYPE C}}} &
  f_{i+8}\in S_{2}(\Gamma_1(29\cdot q_0^2))\ar@/_2pc/ @{-}[r]_{\text{mod } 29\atop \raisebox{-3ex}{\text{\small TYPE A}}} & f_{i+9}\in S_{16}( q_0^2)\\}
 \end{equation*}

\normalsize

\subsection{Step 4NEW (cf.~Corrigenda to \cite{Di-Symm5} in Section 4 (c) of \cite{BD2023})} 

The aim of this step is to introduce the Micro-Good Dihedral prime $43$ into the level of the modular form produced in Step 3NEW (see page 631 of \cite{Di-Symm5} for the definition of Micro-Good Dihedral, which we abbreviate by MGD in what follows). This step will be performed by a single congruence modulo $43$, doing supercuspidal level raising.

To performe this supercuspidal level raising, we want to apply Theorem 1.2 in \cite{BD2023}. Let us check that the conditions of this theorem are satisfied. Since $p=43$, $k=16$ and $N=q_0^2$, the condition $p>\max\{k+1,6\}$ and $p\nmid N$ are satisfied. As explained in  \cite{Di-Symm5}, the good-dihedral prime $q_0$ in the level implies that the residual image satisfies the hypothesis of Theorem 1.2 in \cite{BD2023}. Thus, we met all the conditions to apply this theorem, obtaining a modular weight-$2$ lift of this residual representation corresponding to a newform of level $43^2q_0^2$ which is supercuspidal at $43$. More specifically, it follows from Proposition 3.10 in loc. cit. that the $43$-adic Galois representation just constructed is of type $(2,\tau)$ with $\tau\simeq\omega_2^m\oplus\omega_2^{43\cdot m}$ and $m=k+(p+1)(p-2)=1820$. Note that the character $\omega_2^m$ has order $66$. Twisting by $\omega_2^a$ with $a\equiv 35 \Mod{42}$ we obtain another weight-$2$ newform which is supercuspidal at $43$ and whose local type $\tau\simeq\omega_2^{m+44\cdot a}\oplus\omega_2^{43\cdot (m+44\cdot a)}$ corresponds to a character of order $11$.

We conclude that, up to twist, we have a congruence modulo $43$ with a form $f''$ of weight $2$ and level $43^2q_0^2$, 
Good-Dihedral at $q_0$ with supercuspidal local parameter at $43$ given by a character of order $11$.  
Thus $43$ is a MGD prime for this form. Both forms involved in this congruence up to twist are potentially diagonalisable, 
because one of them is Fontaine-Laffaille and the other one is potentially Barsotti-Tate. This congruence is thus of type A.

\subsection{Step 5 (cf.~\cite[Section 3.5]{Di-Symm5})}

In this step we remove the good-dihedral prime $q_0$. In any congruence that follows, the largness of the residual image will no longer be ensured by the presence of the good-dihedral prime $q_0$, and has to be checked by other means.

We start with a modular form $f_i\in S_2(43^2q_0^2)$ and perform a congruence modulo $t_0$. The residual mod $t_0$ representation will either be unramified or has unipotent ramification at $t_0$ (cf.~\cite{Khare-Wintenberger-I}, where the same situation occurs). In addition, the residual image is large (this is proven in \cite[Section 3.5]{Di-Symm5}. The proof uses the presence of the MGD prime $43$). In the unramified case, the necessary and sufficient conditon for Steinberg level raising at $q_0$ is satisfied (cf. Section 3.5 of loc.~cit.). Hence, doing level raising if necessary, in both cases we have a congruence modulo $t_0$ with a modular form $f_{i+1}$ of weight $2$ and level $43^2q_0$, Steinberg locally at $q_0$. In both cases, the type at $q_0$ was changed, thus this is a type C congruence.

Next we consider the prime $q_0$, and the ordinary Hida family containing the $q_0$-adic representation attached to $f_{i+1}$. We specialize to weight $q_0+1$, thus obtaining a congruence with a modular form $f_{i+2}$ of weight $q_0+1$, ordinary at $q_0$, with level $43^2$. The residual image mod $q_0$ is large, as shown in \cite[Section 3.5]{Di-Symm5}. Since both forms are ordinary, this is a type B congruence. 

We draw a diagram of the congruences performed in this step.

 \small

  \begin{equation*}
  \xymatrix@C-1pc{f_{i}\in S_2(43^2q_0^2) \ar@/_2pc/ @{-}[r]_{\text{mod } t_0\atop \raisebox{-3ex}{\text{\small TYPE C}}} & f_{i+1}\in S_2(43^2q_0)
  \ar@/_2pc/ @{-}[r]_{\text{mod } q_0\atop \raisebox{-3ex}{\text{\small TYPE B}}} & f_{i+2}\in S_{q_0+1}(43^2)\\}
 \end{equation*}

  \normalsize

\subsection{Step 6 (cf.~\cite[Section 3.6]{Di-Symm5})}

This step is similar in nature to Step 2. The goal is to reduce the weight from $q_0+1$ to weight $2<k\leq 14$, and we will use again the WRGC method. The only significant difference is that the Good-Dihedral prime $q_0$ is no longer in the level, therefore a more sofisticated variant of 
the  WRGC method is applied, that takes care of ensuring that residual images are large all through the process. This variant uses the Micro-Good Dihedral prime $43$, and therefore throughout this step we avoid to work in characteristics $p=11$ and $p=43$. The reader can find all details in \cite[Section 3.6]{Di-Symm5}.

The generic step is exactly the same as in Step 2, thus all congruences are of type A. This step finishes with a modular form of level $43^2$ and weight $2<k\leq 14$.

\subsection{Step 7 (cf.~Corrigenda to \cite{Di-Symm5} in Section 4 (a) of \cite{BD2023})}

In Step 7 of the skeleton chain, as described in \cite[Section 3.7]{Di-Symm5}, the mistake that was corrected in the Corrigenda appears only once, and can easily be corrected applying ideas that we have already described in Step 3NEW above. Thus, we will summarize Step 7 following Section 3.7 of loc.~cit., except for the erroneous congruence, which we will replace by a different one (as was done in the Corrigenda), explaining it in full detail. 

The aim of  Step 7 is to accommodate the weight to $16$, like in Step 3NEW above. We will apply the same sequence of congruences as in Step 3NEW, the only difference is that, since the Good-Dihedral prime $q_0$ is no longer in the level, we need to check that the residual images are large with a different technique, involving the MGD prime $43$. In order to preserve this MGD prime in the level, we will avoid working in characteristics $p=11$ and $p=43$. 

We start with a modular form $f_i$ of weight $2<k\leq 14$, level $43^2$ and trivial nebentypus. We perform a congruence mod 47, with a modular form $f_{i+1}$ of weight $2$,  level $43^247$  and nebentypus $\mu=\omega^{k-2}$, which  has conductor $47$ and order 23. Note that this congruence is of type A. The residual image is large, as proven in \cite[Section 3.7]{Di-Symm5}.

 If we reduce (modulo $23$) the local $23$-adic Galois representation attached to $f_{i+1}$, we obtain a large residual image, as shown in \cite[Section 3.7]{Di-Symm5}.
The residual representation $\overline{\rho}_{23}(f_{i+1})$ is either unramified at $47$ or has unipotent ramification at $47$. Accordingly, we divide into two cases (exactly as in Step 3NEW above):

\begin{itemize}

\item Case 1: There is unipotent ramification at $47$. In this case, we take a minimal modular lift, corresponding to a weight 2 modular form $f_{i+2}$ of level $43^247$ that is Steinberg locally at 47.  We have a type C congruence, since we have changed the local type at $47$. 

Next, we move to characteristic 47 and consider $\overline{\rho}_{47}(f_{i+2})$. By considering the Hida family deforming the residual representation attached to the $47$-adic Galois representation corresponding to $f_{i+2}$, we have a congruence with a modular form $f_{i+3}$ of weight $48$ and level $43^2$. Both forms in this congruence are ordinary at $47$, thus this is a type B congruence. 

We need to check that the residual representation $\overline{\rho}_{47}(f_{i+2})$ has a large image. To do so, we will consider all possibilites 
for the residual image, following Dickson's classification of maximal subgroups of $\mathrm{PSL}(2, \mathbb{F})$ for a finite field $\mathbb{F}$ 
of positive characteristic $p$ (with $p=47$ in this case), namely reducible, dihedral, exceptional or large. First of all, since we are working 
in characteristic $47\not=11, 43$, the presence of the MGD prime ensures that the residual image is absolutely irreducible, ruling out the 
reducible case. Moreover, the local type at $43$ ensures that the projective residual image of inertia at $43$ has order $11$, and in the 
exceptional cases in Dickson's classification no element of order $11$ can occur. Therefore the only possibilities for the image of 
$\overline{\rho}_{47}(f_{i+2})$
are dihedral or large. Assume that the image is dihedral, and let $K/\mathbb{Q}$ be the quadratic field such that the restriction of 
the 
image of $\overline{\rho}_{47}(f_{i+2})$ to $G_K$ becomes reducible. The residual representation $\overline{\rho}_{47}(f_{i+2})$ can only be ramified 
at $43$ and $47$, and the same applies to $K/\mathbb{Q}$. However, the order of the residual image of the inertia group at $43$ is $11$, which is odd, 
thus the quadratic extension $K/\mathbb{Q}$ cannot ramify at $43$. Thus $K/\mathbb{Q}$ ramifies at $p=47$, and we can apply Lemma 3.1 of 
\cite{Di-Symm5} to deduce that one of the equalities $47=2k-1$ or $47=2k-3$ holds, where $k$ is the Serre weight of the residual representation 
$\overline{\rho}_{47}(f_{i+2})$. 
But in our case the Serre weight is $k=2$ or $k=48$ (this is so because the modular form $f_{i+2}$ has weight  two and is Steinberg at $p=47$), which do not satisfy the equalities above. Hence the residual image is not dihedral, and we conclude that it must be large.

 We draw a diagram of the congruences performed in Step 7 so far for this case. 

 \small

  \begin{equation*}
  \hskip-2cm\xymatrix@C-1pc{f_{i}\in S_k(43^2) \ar@/_2pc/ @{-}[r]_{\text{mod } 47\atop \raisebox{-3ex}{\text{\small TYPE A}}} & f_{i+1}\in S_2(\Gamma_1(43^2 47))
  \ar@/_2pc/ @{-}[r]_{\text{mod } 23\atop \raisebox{-3ex}{\text{\small TYPE C}}} &
  f_{i+2}\in S_{2}(43^247)\ar@/_2pc/ @{-}[r]_{\text{mod } 47\atop \raisebox{-3ex}{\text{\small TYPE B}}} & f_{i+3}\in S_{48}(43^2) \\}
 \end{equation*}

\normalsize

\item Case 2: There is no ramification at $47$.
Note that this is the point where the mistake occurred in \cite[Section 3.7]{Di-Symm5}: it was erroneously claimed that in Case 2 
a non-minimal lift with Steinberg ramification at $47$ always exists. Since this need not be the case, what we do is to 
take a minimal lift, corresponding to a modular form $f_{i+2}$ of weight $2$ and  level $43^2$. This is a type C congruence 
(we are doing level lowering at $47$), and we have already pointed out, before splitting the proof into two cases, why the 
residual image is large.

It is well known (and follows from multiplying by the Hasse invariant $E_{47-1}$) that there is a mod $47$ congruence between 
$f_{i+2}$ and a newform $f_{i+3}$ of level $43^2$ and weight $48$. As in Step 3NEW above, it can be shown that the $47$-adic 
Galois representation attached to $f_{i+3}$ is potentially diagonalisable locally at $47$, thus this is a type A congruence. 
It remains to check that the image of $\overline{\rho}_{47}(f_{i+2})$ in this congruence is large. The reasoning is similar 
to the one for Case 1. Namely, we need to check that the residual image cannot be reducible, dihedral or exceptional. The same 
argument as above, involving the presence of the MGD prime $43$, together with the fact that the ramification at $43$ has order 
$11$, rules out the reducible and exceptional cases. Since the order of ramification at $43$ is odd, the quadratic field 
extension $K/\mathbb{Q}$ involved in the dihedral case must ramify at $47$ only. An application of Lemma 3.1 of \cite{Di-Symm5}, 
together with the fact that the Serre weight of the residual representation $\overline{\rho}_{47}(f_{i+2})$ is $k=2$ (the input form is a form of weight $2$ and level $43^2$), ensures that the image cannot be dihedral. Thus we have large image, as was to be proven.

The diagram corresponding to the congruences performed so far in Step 7 in Case 2 is the following:

 \small

  \begin{equation*}
  \hskip-2cm\xymatrix@C-1pc{f_{i}\in S_k(43^2) \ar@/_2pc/ @{-}[r]_{\text{mod } 47\atop \raisebox{-3ex}{\text{\small TYPE A}}} & f_{i+1}\in S_2(\Gamma_1(43^2 47))
  \ar@/_2pc/ @{-}[r]_{\text{mod } 23\atop \raisebox{-3ex}{\text{\small TYPE C}}} &
  f_{i+2}\in S_{2}(43^2)\ar@/_2pc/ @{-}[r]_{\text{mod } 47\atop \raisebox{-3ex}{\text{\small TYPE A}}} & f_{i+3}\in S_{48}(43^2) \\}
 \end{equation*}

\normalsize

\end{itemize}

In both cases we end up with a modular form $f_{i+3}$ of weight $48$ and level $43^2$, so we can treat them jointly from this point on. We still need to reduce the weight from $48$ to $16$, exactly like in Step 3NEW.
The next congruences are exactly the same as in Step 3NEW (hence we skip the explanations concerning the existence of the modular forms involved in the congruences), with the only difference that instead of the Good-Dihedral prime $q_0$, the MGD prime $43$ must suffice to ensure large image. We will check at each congruence that this is the case.

As in Step 3NEW, we move to $p=53$, and modulo $p$ we have a congruence with a modular form $f_{i+4}$ of weight $2$  and level $43^253$ with nebentypus $\omega^{46}$ of conductor $53$ and order $26$. This is a type A congruence. Using the presence of the MGD prime $43$ and Lemma 3.1 of \cite{Di-Symm5}, it can be proven that the residual image is large, exactly as in the congruences above.

Next, following a method of Khare to reduce the weight, we are going to work modulo the auxiliary prime $13$ (which divides the order of $\omega^{46}$) to modify the nebentypus. We reduce modulo $13$ and we take a minimal lift, which corresponds to a modular form $f_{i+5}$ of weight $2$ and level $43^2 53$. The character giving ramification at $53$ in this minimal lift, in other words, the nebentypus of $f_{i+5}$,  is $\omega^{26}$, a quadratic character. 
This is a type C congruence, since the local type at $53$ was modified. Concerning the image of the residual representation 
$\overline{\rho}_{13}(f_{i+4})$, note that by the presence of the MGD prime $43$ (and the fact that the ramification at $43$ has order $11$), the 
only possibilities are dihedral or large. Assume that the residual image is dihedral, and call $K/\mathbb{Q}$ the quadratic field such that the 
residual representation becomes reducible when restricted to $G_K$. This extension can only ramify at those primes which are ramification primes 
for the residual representation $\overline{\rho}_{13}(f_{i+4})$, that is to say, $13, 43$ and $53$. Since the image of the ramification at $43$ has 
odd order, $K/\mathbb{Q}$ is unramified at $43$. Applying Lemma 3.1, we conclude that it cannot ramify at $p=13$ (the residual Serre weight is 2). 
Thus, the only possibility is that $K=\mathbb{Q}(\sqrt{53})$. But the prime $43$ is split in this extension, hence the restriction to $G_K$ of the 
residual representation $\overline{\rho}_{13}(f_{i+4})$ (which is reducible by definition of $K$) must contain the image of all the decomposition 
group at $43$, which is absolutely irreducible.
This contradiction proves that the residual image is large.

Next we reduce again modulo $53$. Here, the possible values for the Serre weight of $\overline{\rho}_{53}(f_{i+5})$ (up to twist) are (according to 
\cite{Savitt_Corrigenda})  $26+2=28$ or $53+3-28=28$. By taking a modular minimal lift, we obtain a congruence, up to twist,  
with a modular form $f_{i+6}$ of weight $28$ and level $43^2$. This is a type A congruence. To show that the residual image 
of $\overline{\rho}_{53}(f_{i+5})$ is large, we need to exclude the reducible, exceptional and dihedral cases. The first two cases are excluded because of the presence of the MGD prime, as seen above. To exclude the dihedral case, first note that  ramification at $43$ has odd order, thus the only dihedral case that we have to consider is the so-called bad dihedral case, corresponding to the quadratic extension $K/\mathbb{Q}$ ramifying only at the residual characteristic $p=53$. It will not suffice to use Lemma 3.1 of \cite{Di-Symm5}, since in this case we have the numerical equation $53=2\cdot 28-3$. Luckily, Lemma 3.2 of loc.~cit. applies, ruling out the particular case $p=2k-3$. Hence the residual image must be large.

\small

  \begin{equation*}
  \hskip-2cm\xymatrix@C-1pc{f_{i+3}\in S_{48}(43^2) \ar@/_2pc/ @{-}[r]_{\text{mod } 53\atop \raisebox{-3ex}{\text{\small TYPE A}}} & 
  f_{i+4}\in S_2(\Gamma_1(43^2\cdot 53))
  \ar@/_2pc/ @{-}[r]_{\text{mod } 13\atop \raisebox{-3ex}{\text{\small TYPE C}}} &
  f_{i+5}\in S_{2}(\Gamma_1(43^2\cdot 53))\ar@/_2pc/ @{-}[r]_{\text{mod } 53\atop \raisebox{-3ex}{\text{\small TYPE A}}} & f_{i+6}\in S_{28}( 43^2)\\}
 \end{equation*}

\normalsize

To reduce the weight further, we proceed as in Step 3NEW and repeat the above steps, using different auxiliary primes. First, we reduce modulo $p=29$, and consider a modular weight two lift of level $29\cdot 43^2$, and nebentypus $\omega^{26}$, which is a character of conductor $29$ and order $14$. This is a type A congruence, and the residual image is large because the reducible and exceptional cases are excluded by the presence of the MGD prime $43$, and the dihedral case is ruled out by Lemma 3.1 of loc.~cit., as we have seen in some of the other steps above.

Next, we reduce modulo $7$ and take a minimal modular lift corresponding to a modular form $f_{i+8}$ of weight $2$ and level $29 \cdot 43^2$, 
with nebentypus $\omega^{14}$, a quadratic character. This is  a type C congruence, since the local type at $29$ was modified. To show that the 
 image of $\overline{\rho}_{7}(f_{i+7})$ is large, we discard as usual the reducible and exceptional cases via the presence of the MGD prime $43$. 
 For the case of dihedral image, assume that $K/\mathbb{Q}$ is the quadratic extension such that the residual representation 
 $\overline{\rho}_{7}(f_{i+7})$ 
 becomes reducible when restricted to $G_K$. Like in other cases above, the ramification of $K/\mathbb{Q}$ is limited by the ramification of the 
 residual representation $\overline{\rho}_{7}(f_{i+7})$, namely to the primes $7, 29$ and $43$. The extension cannot be ramified at $43$ because the image of inertia at $43$ has odd order (namely $11$). Lemma 3.1, together with the knowledge of the Serre weights of the residual representation, rules out the possibility that $K$ ramifies at $7$. Hence the only possibility that remain is $K=\mathbb{Q}(\sqrt{29})$. In \cite[Section 3.7]{Di-Symm5}, this possibility is excluded by an argument using class field theory. The interested reader can look up the full discussion in pag. 640  in \cite[Section 3.7]{Di-Symm5}.

Finally, we reduce modulo $29$, and applying \cite{Savitt_Corrigenda} we see that the residual Serre weight, up to twist, is $16$. Thus we get a 
congruence with a modular form $f_{i+9}$ of weight $16$ and level $43^2$. This is a type A congruence.  Regarding the residual image of
$\overline{\rho}_{29}(f_{i+8})$, the MGD prime $43$ ensures that it is not reducible or exceptional, and since the order of ramification at $43$ is odd, the only dihedral case that has to be considered is the bad dihedral case. By Lemma 3.1 of \cite{Di-Symm5}, we have that either $p=2k-1$ or $p=2k-3$, where $k=16$ is the residual Serre weight. Since in this case $29= 2 \cdot 16 -3$, an application of Lemma 3.2 of \cite{Di-Symm5} provides a contradiction. Thus the residual image is large.

 \small

  \begin{equation*}
  \hskip-2cm\xymatrix@C-1pc{f_{i+6}\in S_{28}(43^2) \ar@/_2pc/ @{-}[r]_{\text{mod } 29\atop \raisebox{-3ex}{\text{\small TYPE A}}} & f_{i+7}\in S_2(\Gamma_1(29\cdot 43^2))
  \ar@/_2pc/ @{-}[r]_{\text{mod } 7\atop \raisebox{-3ex}{\text{\small TYPE C}}} &
  f_{i+8}\in S_{2}(\Gamma_1(29\cdot 43^2))\ar@/_2pc/ @{-}[r]_{\text{mod } 29\atop \raisebox{-3ex}{\text{\small TYPE A}}} & f_{i+9}\in S_{16}( 43^2)\\}
 \end{equation*}

\normalsize

\subsection{Step 8 (cf.~\cite[Section 3.8]{Di-Symm5})}

The purpose of this step is to introduce a nebentypus of order $8$ at $17$, which is necessary for the ad hoc reasonings that enable
the congruence occurring in Step 9.
Thus, we consider a congruence modulo $17$ with a modular form of weight $2$ and level $17\cdot 43^2$, with nebentypus of conductor $17$ 
given by the character $\omega^{14}$ of order $8$ (cf.~\cite[Section 3.8]{Di-Symm5} for the existence of this modular form). 
This is a type A congruence (between a Fontaine-Laffaille and a potentially Barsotti-Tate representation), and the residual 
image is large (cf. Section 3.8 of loc.~cit.).

\subsection{Step 9 (cf.~\cite[Section 3.9]{Di-Symm5})}

The goal of this step is to remove the MGD prime $43$ from the level of the modular form produced in Step 8. We start with a modular form $f_i$ 
of weight $2$ 
and level $17\cdot 43^2$ and nebentypus of order $8$, ramified at $17$. The residual representation modulo $11$ attached to $f_i$ has large 
image, as proved in page 642, \cite[Section 3.9]{Di-Symm5}. Moreover, it is known that the residual representation $\overline{\rho}_{11}(f_i)$ 
will have either trivial 
or unipotent ramification at $43$ (cf. Section 3.9 of loc.~cit.). In both cases (by level raising in the first case), there exists a modular lift 
corresponding to a weight $2$ modular form $f_{i+1}$ of level $17\cdot 43$, Steinberg locally at $43$, and has a nebentypus of order $8$ at 17 
(cf.~Section 3.9 of loc.~cit.). We have changed the local type at $43$, and since we are doing a mod 11 congruence and 11 is not in the level, 
both representations are Barsotti-Tate. Hence it is a type C congruence. After Step 10, we will draw a common diagram for the congruences in Step 9 
and 10.

\subsection{Step 10 (cf.~\cite[Section 3.10]{Di-Symm5})} 

In this step we consider the mod $43$ representation attached to the modular form $f_{i+1}$ provided by the previous step. Using direct computations, 
the residual representation $\overline{\rho}_{43}(f_{i+1})$ has been checked to be irreducible (cf.~\cite[Section 3.10]{Di-Symm5}), and using 
information on the ramification it 
can also be proved that the residual image is large (cf. Section 3.10 of loc.~cit.). We consider the Hida family containing this residual 
representation (the modular form $f_{i+1}$ is known to be ordinary at $43$ because it has weight $2$ and is Steinberg locally at $43$). We 
specialise at weight $44$, obtaining a congruence with a modular newform $f_{i+2}$ of weight $44$, and level $17$, with nebentypus of order $8$. 
This congruence is of type B. 

We draw a diagram for Steps 9 and 10 together.

 \small

   \begin{equation*}
  \xymatrix@C-1pc{f_{i}\in S_2(\Gamma_1(17\cdot 43^2)) \ar@/_2pc/ @{-}[r]_{\text{mod } 11\atop \raisebox{-3ex}{\text{\small TYPE C}}} & f_{i+1}\in S_2(\Gamma_1(17\cdot 43))
  \ar@/_2pc/ @{-}[r]_{\text{mod } 43\atop \raisebox{-3ex}{\text{\small TYPE B}}} & f_{i+2}\in S_{44}(\Gamma_1(17))}
 \end{equation*}

\normalsize

\subsection{Step 11 (cf.~\cite[Section 3.11]{Di-Symm5})} 

This step consists in checking by direct computation that the space of newforms of weight $44$, level $17$ with any nebentypus of order $8$ at $17$ consists of a single orbit, under Galois conjugation on coefficients, of Hecke eigenforms. Thus, since Galois conjugation is a valid move, the chain constructed so far links the inital modular form with any of the newforms in this space.

This computation was originally carried out in the computer system MAGMA (see \cite[Section 3.11]{Di-Symm5} for details on the computation). We have also checked this computation in the open-source computer systems SageMath and Pari/GP at the request of the referee; we provide the computer code and the output in \cite{Dieulefait-Website}.\footnote{Incidentally, the space of newforms that we are interested in is also computed in the \emph{L-functions and modular forms database (LMFDB)}. In particular, these computations show that there is a single Galois orbit; see \texttt{http://www.lmfdb.org/ModularForm/GL2/Q/holomorphic/17/44/d/}}

The goal of the last four steps is to connect a CM modular form of weight $2$ and level $27$ with some modular form in the space of newforms of 
weight $44$, level $17$ with any nebentypus of order $8$ at $17$. Each step consists of a single congruence. We will illustrate the four congruences together 
in a diagram below.
We describe these four steps in reverse order, starting from the CM form, since we do not know a priori which of the newforms in 
the space of newforms of weight $44$, level $17$ with any nebentypus of order $8$ at $17$ we will hit (and it is not relevant, since any pair of modular 
forms in this space is linked by a Galois conjugation).

\subsection{Step 15 (cf.~\cite[Section 3]{Di-mfold})}\label{Step:15}

We start from the modular form $g_1:=g\in S_2(27)$ with CM by the quadratic imaginary field $K=\mathbb{Q}(\sqrt{-3})$ and rational coefficients. We check that there is a congruence modulo $13$ with a newform $g_2$ of level $27\cdot 43$, weight $2$, Steinberg at $43$ (the condition for level raising is satisfied, cf.~\cite[Section 3]{Di-mfold}). By direct computation\footnote{This computation has been made in the open-source software SageMath and runs on a personal laptop. The code is available in \cite{Dieulefait-Website}.} in the space of newforms of weight $2$ and level $27\cdot 43$, we identify the form $g_2$ that satisfies this congruence. It has coefficient field $\mathbb{Q}_{g_2}$ of degree $10$. 
The congruence between $g_1$ and $g_2$ is a type C congruence, since the local type at $43$ has changed, and both forms are Barsotti-Tate locally at $13$. 

Since one of the forms in the congruence has CM by the quadratic imaginary field $K=\mathbb{Q}(\sqrt{-3})$, 
the projectivisation of the
residual image of $\overline{\rho}_{13}(g_1)\equiv
\overline{\rho}_{13}(g_2)\pmod{13}$ is a dihedral group, induced from a character of $G_K$, of order prime to $13$. Therefore $\overline{\rho}_{13}(g_1)(G_{ \mathbb{Q}(\zeta_{13}) })$ is absolutely irreducible (since $\mathbb{Q}(\zeta_{13})$ and $K$ are linearly disjoint). Lemma \ref{lem:ForDihedral} implies that the image $\overline{\rho}_{13}(f)(G_{ \mathbb{Q}(\zeta_{13}) })$ is adequate$+$.

\subsection{Step 14 (cf.~\cite[Section 3]{Di-mfold})}

The $3$-adic Galois representations attached to $g_2$, corresponding to the primes of the coefficient field $\mathbb{Q}_{g_2}$ of
$g_2$ above $3$, are potentially Barsotti-Tate. Moreover, there is a prime $\mathfrak{p}$ of $\mathbb{Q}_{g_2}$ dividing $3$, 
with inertial degree $3$. 
Direct computations, together with Dickson's Classification Theorem, easily give that the image of $\overline{\rho}_{\mathfrak{p}}(g_2)$ 
contains $\mathrm{SL}_2(\mathbb{F}_{27})$ (cf.~\cite[Section 3]{Di-mfold}). By elementary group theory arguments, the same holds for the 
restricton of $\overline{\rho}_{\mathfrak{p}}(g_2)$ to $\mathbb{Q}(\zeta_{3})$, thus by Lemma \ref{rem:image-chain} this restriction has adequate$+$ image.

By direct computations\footnote{This computation has also been made in SageMath using a personal laptop. The code is also available in \cite{Dieulefait-Website}.}, we check that this residual representation has a modular lift corresponding to a newform of level $43$, 
Steinberg at $43$ and weight $4$, which we call $g_3$. This newform has field of coefficients $\mathbb{Q}_{g_3}$ of degree $6$. 
Moreover, it can be shown (cf.~\cite[Section 3]{Di-mfold}) that locally at $3$, the $3$-adic Galois representation attached to 
$g_3$ is ordinary (this is seen by the explicit computation of $g_3$, since $3$ does not divide $a_3(g_3)$). Thus it is ordinary and crystalline (3 is not in the level), hence potentially diagonalizable. Thus we are in a type A congruence.

\subsection{Step 13 (cf.~\cite[Section 3]{Di-mfold})}

In this step we perform a congruence mod $17$ between the form $g_3$ and a form $g_4$ of weight $2$, level $17\cdot 43$ and nebentypus $\omega^2$, a character of conductor $17$ and order $8$. This is a type A congruence. In \cite[Section 3]{Di-mfold}, it is shown that the residual image is large.

\subsection{Step 12 (cf.~\cite[Section 3]{Di-mfold})}

Recall that the newform $g_4$ is a weight $2$ form, of level  $17\cdot 43$, Steinberg locally at $43$, and a nebentypus at $17$ of order $8$. These are exactly the same conditions that appear at the beginning of Step 10. Thus, we can reproduce the same arguments than in Section 10, to conclude that we have a congruence with a form $g_5$ of weight $44$, level $17$ and nebentypus of conductor $17$ and order $8$, that the residual image in the congruence is large and that the congruence is of type B. Note that we have landed in the single-orbit space of newforms described in Step 11, so the form $g_5$ is conjugated to the form obtained in Step 10. 

We draw a diagram illustrating the last four steps, starting from Step 12 and ending in Step 15. Note that this is the point where there is a type B congruence (in Step 12) that is not preceeded by a type C congruence (see Steps 10 and 11). 
    
\small

   \begin{equation*}
  \hskip-2cm \xymatrix@C-1pc{g_5 \in S_{44}(\Gamma_1(17)) \ar@/_2pc/ @{-}[r]_{\text{mod } 43\atop \raisebox{-3ex}{\text{\small TYPE B}}} & g_4 \in S_{2}(\Gamma_1(17\cdot 43)) 
  \ar@/_2pc/ @{-}[r]_{\text{mod } 17\atop \raisebox{-3ex}{\text{\small TYPE A}}} & g_3\in S_{4}(43)
   \ar@/_2pc/ @{-}[r]_{\text{mod } 3\atop \raisebox{-3ex}{\text{\small TYPE A}}} & g_2\in S_2(3^343)
   \ar@/_2pc/ @{-}[r]_{\text{mod } 13\atop \raisebox{-3ex}{\text{\small TYPE C}}} & g_1\in S_2(3^3) }
 \end{equation*}

\normalsize

\begin{rem}\label{rem:UnPasoAntes} Note that at any step $f_i\equiv f_{i+1} \pmod{p_i}$ in the skeleton chain, if $p_i$ does not belong to the
finite list $\mathcal{L}=\{7, 11, 13, 23, 43, 47, q_0, t_0\}$, the congruence is of type A. \end{rem}

\begin{rem}\label{rem:ExtendedChain} In Section \ref{sec:Diagrams}, we will need to enlarge the skeleton chain recalled above. In fact,
at the steps where a type B or type C congruence occurs in the $2$-dimensional skeleton chain, we will not be able to proceed directly in our reasonings in Section \ref{sec:Diagrams},
 and we will need to enlarge the skeleton chain by adding several new links. More precisely, at each type C congruence between $f_i$ and $f_{i+1}$ modulo $p_i$,
we will also need modular forms $f'$ and $f''$, of weight $j'$ belonging to the interval $[2n, 2n + p_i^2-1]$, 
satisfying that $f_{i}\equiv f'$ modulo $p_{i}$ with $\rho_{p_i}(f_i)$ and $\rho_{p_i}(f')$ connecting locally at all primes $p\not=p_i$ and 
$f''\equiv f_{i+1}$ modulo $p_i$ with
$\rho_{p_i}(f_{i+1})$ and $\rho_{p_i}(f'')$ connecting locally at all primes $p\not=p_i$.
We will add new links to the skeleton chain as follows:
     \begin{equation*} \hskip-2cm \xymatrix{f_{i}\ar@/_2pc/ @{-}[r]_{\text{mod } p_i} & f_{i+1} } \qquad \rightsquigarrow \qquad
\xymatrix{f_{i}\ar@/_2pc/ @{-}[r]_{\text{mod } p_i} & f'\ar@/_2pc/ @{-}[r]_{\text{mod } p_i} & f_{i+1} \ar@/_2pc/ @{-}[r]_{\text{mod } p_i}
 & f'' \ar@/_2pc/ @{-}[r]_{\text{mod } p_i}& f_{i+1}} \end{equation*}

At each type B congruence, say between $f_{i+1}$ and $f_{i+2}$ modulo $p_{i+1}$, 
except for the one in Step 12, we have noted (cf.~Remark \ref{rem:CongruenceTypes}) that it is preceeded by a Type C congruence
$f_{i}\equiv f_{i+1}$ modulo $p_{i}$. 
In addition to the modifications required because of the presence of the type C congruence, we will also need two extra modular
forms, $f'''$ and $f^{iv}$, with weight comprised in the interval $[2n, 2n + p_{i+1}^{2}-1]$, such that both $\rho_{p_{i+1}}(f''')$
and $\rho_{p_{i+1}}(f^{iv})$ are potentially diagonalisable locally at $p_{i+1}$ and such that: $f'''$ and $f_{i}$ are congruent modulo
$p_{i+1}$ with $\rho_{p_{i+1}}(f''')$ connecting locally at all $v\nmid p_{i+1}$ to $\rho_{p_{i+1}}(f_i)$, and $f^{iv}$ and $f_{i+2}$ are 
congruent modulo $p_{i+1}$ with $\rho_{p_{i+1}}(f^{iv})$ connecting locally at all $v\nmid p_{i+1}$ to $\rho_{p_{i+1}}(f_{i+2})$.
The modifications in the skeleton chain corresponding to both congruences $f_i\equiv f_{i+1}$ modulo $p_i$ and $f_{i+1}\equiv f_{i+2}$ modulo 
$p_{i+1}$ are the following:
 
\begin{equation*} \hskip-2cm \xymatrix{f_{i}\ar@/_2pc/ @{-}[r]_{\text{mod } p_i} & 
f_{i+1} \ar@/_2pc/ @{-}[r]_{\text{mod } p_{i+1}} &
f_{i+2} } \qquad \rightsquigarrow \qquad\end{equation*}
\begin{equation*}\xymatrix{f_{i}\ar@/_2pc/ @{-}[r]_{\text{mod } p_{i+1}} & 
f'''\ar@/_2pc/ @{-}[r]_{\text{mod } p_{i+1}} &
f_{i}\ar@/_2pc/ @{-}[r]_{\text{mod } p_i} &
f'\ar@/_2pc/ @{-}[r]_{\text{mod } p_i} 
& f_{i+1} \ar@/_2pc/ @{-}[r]_{\text{mod } p_i}
& f'' \ar@/_2pc/ @{-}[r]_{\text{mod } p_i}
& f_{i+1} \ar@/_2pc/ @{-}[r]_{\text{mod } p_{i+1}} 
& f_{i+2} \ar@/_2pc/ @{-}[r]_{\text{mod } p_{i+1}} 
& f^{iv} \ar@/_2pc/ @{-}[r]_{\text{mod } p_{i+1}}
& f_{i+2} } \end{equation*}

At Step 12 we have a Type B congruence between the modular forms $g_5$ and $g_4$ modulo the prime $p=43$.
For this step, we need a modular form $f'$ of weight $j'$ belonging to the interval $[2n, 2n +43^2-1]$, 
with $g_5\equiv f'\pmod{43}$, and $\rho_{43}(f')$ potentially diagonalisable locally at $43$ and connecting
locally at each $v\nmid p$ to $\rho_{43}(g_5)$. We modify the skeleton chain as follows:

     \begin{equation*} \hskip-2cm \xymatrix{g_{5}\ar@/_2pc/ @{-}[r]_{\text{mod } 43} & g_{4} } \qquad \rightsquigarrow \qquad
\xymatrix{g_{5}\ar@/_2pc/ @{-}[r]_{\text{mod } 43} & f'\ar@/_2pc/ @{-}[r]_{\text{mod 43} } & g_{4}} \end{equation*}

Using Theorem \ref{thm:BigWeight}, we can construct the required modular forms  for each of the
steps where we have a congruence of type B or C. Fix such modular forms for each of these steps. We will call $\widetilde{\mathcal{C}}$ the \emph{extended} $2$-dimensional skeleton chain,
obtained from $\mathcal{C}$  by adding  all these new links as indicated above.
\end{rem}

\section{Preliminary steps}\label{sec:PreliminarySteps}

Let us start by introducing some notation for the Galois representations attached to automorphic forms. As in \cite{BLGGT}, if $F$ is either a CM field or a totally real field, and if $\pi$ is a RACP automorphic representation of $\mathrm{GL}_n(\mathbb{A}_F)$, then we can attach to it a Galois representation $r_{\ell,\iota}(\pi)\colon G_F \rightarrow \mathrm{GL}_n(\overline{\mathbb{Q}}_\ell)$ for each rational prime $\ell$ and isomorphism $\iota\colon\overline{\mathbb{Q}}_\ell\simeq \mathbb{C}$. 

Because we have fixed our choices of isomorphisms $\iota\colon\overline{\mathbb{Q}}_\ell\simeq \mathbb{C}$ from the beginning, it makes sense to denote these representations just by $r_\ell(\pi)$. Moreover, when $n=2$, we will denote it by $\rho_\ell(\pi)$ (to be most consistent with the rest of this paper and with previous work of the authors).

With that out of the way, we start the proper content of this section. Assume we are given $f\in S_k(1)$ a cuspidal modular eigenform and an
$n$-dimensional RACP representation $\pi$ of $\mathrm{GL}_n(\mathbb{A}_{\mathbb{Q}})$, which is automorphic,  with level coprime to $3$. Let $\{\rho_p(f)\}$ and
$\{r_p(\pi)\}$ denote the corresponding compatible systems of Galois representations of $G_{\mathbb{Q}}$.
We want to prove that the compatible system $\{(\rho(f)\otimes r(\pi))_p\}$ is automorphic, by linking it
through a ``safe'' chain to the compatible system $\{(\rho(g)\otimes r(\pi^{(\textrm{end})}))_{p}\}$, where $g$ is a certain modular
form with CM (cf.~Section \ref{sec:2-dim Chain}) and $\pi^{(\textrm{end})}$ is a suitable RACP automorphic representation.

In these preliminary steps we want to
 modify the automorphic
representation $\pi$ and replace it by another automorphic representation, say $\pi'$,
defined over a suitable solvable CM extension of $\mathbb{Q}$,
which has some good properties, and a good behavior with respect to the modular forms belonging to the extended $2$-dimensional skeleton chain $\widetilde{\mathcal{C}}$
(cf.~Section \ref{sec:2-dim Chain}). For example, for any modular form $f_i\in \widetilde{\mathcal{C}}$, we will want that the compatible system corresponding to the tensor product $\{(\rho(f_i)\otimes r(\pi^{'}))_{p}\}$ is
regular. To ensure this condition, we will impose the requirement that the weights of $\pi'$ are $C_0$-very spread (cf.~Definition \ref{def:C-very-spread}).
And most importantly, we need to make sure that automorphy of the new compatible system $\{(\rho(f)\otimes r(\pi'))_{p}\}$
implies automorphy of our original system $\{(\rho(f)\otimes r(\pi))_{p}\}$.
Actually, we will proceed in three steps, passing first to a suitable $\pi'$, then to $\pi''$, and finally to $\pi'''$.

For the rest of the section, fix a modular form $f\in S_k(1)$ and an $n$-dimensional RACP representation $\pi$ of $\mathrm{GL}_n(\mathbb{A}_{\mathbb{Q}})$ with Hodge-Tate weights $\{0, w_1, \dots, w_{n-1}\}$,
satisfying the assumptions of Theorem \ref{thm:main}. Let us recall the setting.
First of all, we may assume that $f$ has no complex multiplication. Thus, by Ribet's theorem, the residual image
$\overline{\rho}_p(f)(G_{\mathbb{Q}})$ is large (i.e.~contains $\SL_2(\mathbb{F}_p)$) for $p$ sufficiently large.
Recall that we say the prime $p$ is \emph{exceptional} for $f$ if the image of $\overline{\rho}_p(f)$ is not large.
Moreover, the compatible system $\{r_{p}(\pi)\}$ is irreducible (where we use the definition in \cite[Section 5, pag. 571]{BLGGT}, i.e.~there exists a positive density set of rational primes
$\ell$ such that $r_{\ell, \iota}(\pi)$ is irreducible).
In addition, the compatible system $\{(\rho(f) \otimes r(\pi))_{p}\}$ is irreducible and regular.

As explained in Section \ref{sec:2-dim Chain}, there exists a ``safe'' chain connecting $f$ to a
modular form $g\in S_2(27)$ with CM by $\mathbb{Q}(\sqrt{-3})$, via a series of congruences (up to twist and Galois conjugation) modulo different prime numbers.
We fix such a $2$-dimensional skeleton chain $\mathcal{C}$, and further we fix an extended $2$-dimensional skeleton chain
$\widetilde{\mathcal{C}}$ (cf.~Remark \ref{rem:ExtendedChain}).

First we define several constants that depend on the skeleton chain $\mathcal{C}$, the extended skeleton chain $\widetilde{\mathcal{C}}$, the dimension $n$ and the weights $\{0, w_1, \dots, w_{n-1}\}$ of the
representation $\pi$.
We adopt the notation in Section \ref{sec:2-dim Chain}, in particular, we denote by $q_0$ the good dihedral prime of the chain $\mathcal{C}$ and $t_0$
denotes the order of the character giving the ramification at $q_0$.

\begin{itemize}

\item $\mathcal{L}=\{7, 11, 13, 23, 43, 47, q_0, t_0\}$;

\item $\mathcal{P}_1:=\{p \text{ prime such that at (at least one) step of the chain, a congruence occurs modulo $p$}\}$;

\item $\mathcal{P}_2:=\{p \text{ prime such that $p$ is exceptional for some $f\in \widetilde{\mathcal{C}}$}\}$;

\item $k_{\mathrm{max}}:=\max\{k: \text{ some modular form }f\in \mathcal{C} \text{ has weight }k\}$;

\item $w_{\mathrm{max}}=\max\{w_i: i=1, \dots, n-1\}$;

\item $N_{\mathrm{max}}:=\mathrm{lcm}\{N: N \text{ is the level of some } f\in \widetilde{\mathcal{C}}\}$;

\item $N_{\pi}$ is defined as the product of the primes dividing the level of $\pi$;

\item $C_0:=\max\{k_{\mathrm{max}}, \ell^2 + 2n: \ell\in\mathcal{L}\}$;

\end{itemize}

For future reference, we also fix natural numbers $k'$, $\widehat{k}$ satisfying:
 \begin{itemize}
\item $k'>k_{\max} + w_{\mathrm{max}}$;
  \item $\widehat{k}\geq 2C_0$
\item $\widehat{k}\equiv 1 \pmod{42}$;

 \end{itemize}

 The number $\widehat{k}$ will only be used in Section \ref{sec:Diagram3}, whereas $k'$ will already play its role in the 
 first preliminary step below.

\subsection{First Preliminary step}

 For technical reasons that will be made clear in the second preliminary step,
we need that the dimension  of our automorphic representation is even.
Therefore, the first modification we will do is to tensor our
compatible system $\{r_{p}(\pi)\}$  with a suitable compatible system of $2$-dimensional, dihedral
representations $\{D_{p}\}$, replacing $\{r_{p}(\pi)\}$ by the tensor product $\{(r(\pi)\otimes D)_{p}\}$.

However, we want to avoid dimension $8$, in order to be able to apply Proposition \ref{prop:adequate} to the image 
of the residual Galois representation attached to $\pi$. Thus, if the dimension of $\pi$ is $4$, we skip this step and
proceed directly to the second preliminary step.

Assume that the dimension of $n$ is different from $4$. First of all, we
choose a prime $\ell$ satisfying:

\begin{itemize}

\item $\ell> k_{\mathrm{max}} + w_{\mathrm{max}} + k'n$

\item $\ell > \max\{p: p\text{ divides }N_{\pi} N_{\mathrm{max}}\}$.

 \item $(\ell+1)\nmid (k' -1)n$

 \item $\ell> 2(4n+1)$

\item The image of the restriction of $\overline{\rho}_{\ell}(f)\otimes \overline{r}_{\ell}(\pi)$ to the cyclotomic extension $\mathbb{Q}(\zeta_{\ell})$
is irreducible.
\end{itemize}

Note that the irreducibility of the system $\{(\rho(f)\otimes r(\pi))_{p}\}$, together with
\cite[Proposition 5.3.2]{BLGGT}, ensures that the last condition
is satisfied at a density one set of primes $\ell$.

We construct an induced representation of $G_{\mathbb{Q}}$ as in
Step 1 in the proof of Theorem 7.5 in \cite{PRIMS}. Following their construction, we choose an imaginary quadratic
field $M$ which is linearly disjoint to the compositum of $\mathbb{Q}(\zeta_{\ell})$ with the fixed field of the residual representation $\overline{\rho}_{\ell, \iota}(f)\otimes \overline{r}_{\ell, \iota}(\pi)$ (not just the factor $\overline{r}_{\ell, \iota}(\pi)$), in which $\ell$  and $3$ split, and choose $b>k'$. We consider the character
$\theta_{\ell}:G_M\rightarrow \overline{\mathbb{Q}}_{\ell}^{\times}$, chosen like in Step 1 in the proof of Theorem 7.5.  Note that we can choose this character to be unramified at the primes above 3. This character belongs to a compatible system $\{\theta_{p}\}$, and
we let $\{D_{p}\}=\{\mathrm{Ind}^{G_{\mathbb{Q}}}_{G_M}\theta_{p}\}$.
By quadratic base change, we have that the compatible system $\{(r(\pi)\otimes D)_{p}\}$ is automorphic.

\begin{rem}\label{rem:From_Pi_to_Pi1} From now on, we will work with the tensor product $\{(r(\pi)\otimes D)_{p}\}$ instead of $\{r_{p}(\pi)\}$.  That is to say, 
we replace $\pi$ by the automorphic representation $\pi'$  corresponding to the compatible system $\{(r(\pi)\otimes D)_{p}\}$. In order to make this replacement without loss of generality, we need to check that if $\{(\rho(f)\otimes r(\pi'))_{p}\}$ is automorphic, then $\{(\rho(f)\otimes r(\pi))_{p}\}$ is automorphic. This follows by solvable base change as in step 3 of the proof of Theorem 7.5 of \cite{PRIMS}.
\end{rem}

The automorphic representation $\pi'$ satisfies the following properties:

\begin{enumerate}
 \item The dimension of the compatible system of Galois representations $\{r_{p}(\pi')\}$ is even, equal to $2n$.
 
 \item  The level of $\pi'$ is coprime to $3$.
 \item The  image of $\overline{\rho}_{\ell, \iota}(f)\otimes \overline{r}_{\ell, \iota}(\pi')\vert_{\mathbb{Q}(\zeta_{\ell})}$ is irreducible: this follows from the last condition on $\ell$ and the choice of the field $M$.
 \item $r_{\ell, \iota}(f)\otimes r_{\ell, \iota}(\pi)\otimes D_{\ell}$  is regular (because of the choice of $k'$).
 \item The prime $\ell$ satisfies that $\ell>k_{\max} + w_{\max} + k'n$, and moreover $\ell$ is greater than the primes dividing the levels of $f$ or $\pi$. Thus,
 the representation $\rho_{\ell, \iota}(f)\otimes r_{\ell, \iota}(\pi)\otimes D_{\ell}$ is
 Fontaine-Laffaille at $\ell$.
 \item Furthermore, $\ell>2\cdot(4n+1)$ and $(\ell+1)\nmid (k'-1)n$.
\end{enumerate}

\subsection{Second preliminary step}
Next, we want to
 modify the automorphic
form $\pi'$ and replace it by another automorphic form $\pi''$, such that the residual representations are generically very big (not just
irreducible) at a density one set of rational primes, that the weights are suitably spread, and that the primes of $\mathcal{L}$ are Steinberg in $\pi''$, to mention
just a few desirable properties that will be exploited in
subsequent sections. In order to modify $\pi'$, we will use existence theorems and ALT from \cite{BLGGT}. We need to replace the field $\mathbb{Q}$ by a CM field to be able to apply these theorems.
 Recall that for the first preliminary step we had fixed a large prime $\ell$.  Denote by $N_{\pi'}$ the product of the rational primes dividing the conductor of $\pi'$, and define
$$S=\{\ell\}\cup\mathcal{L}\cup\{p:p\vert N_{\pi'}\}.$$

\begin{lem}There exists a solvable CM extension $F/\mathbb{Q}$  such that:

 \begin{itemize}

\item $F$ is linearly disjoint over $\mathbb{Q}$ from the field
$$\mathbb{Q}_{\text{avoid}}:=\mathbb{Q}(\sqrt{-3}, \zeta_{13}, \zeta_{\ell}, \overline{\mathbb{Q}}^{\ker{\overline{\rho}_{\ell,\iota}(f)\otimes \overline{r}_{\ell,\iota}(\pi')}}, \zeta_p: p\in \mathcal{L}).$$

\item For all places $v$ of $F^+$ above a prime in $S$, $v$ splits completely in $F/F^+$.

\item For all places $u$ of $F$ above $\ell$, the completion $F_u$ of $F$ at $u$ contains $\zeta_{\ell}$ and $\overline{r}_{\ell, \iota}(\pi')\vert_{G_{F_u}}$ is trivial.

\item For all $p\in \mathcal{L}$, for all places $w$ of $F$ above $p$, $\overline{r}_{\ell, \iota}(\pi')\vert_{G_{F_w}}$ is trivial and the cardinality $\kappa_p$ of the residue field of $\mathcal{O}_F/w$ satisfies $\kappa_p\equiv 1 \pmod{\ell}$.

\item For all $p\in N_{\pi'}\setminus\mathcal{L}$, for all places $w$ of $F$ above $p$, $\overline{r}_{\ell, \iota}(\pi')\vert_{G_{F_w}}$ is trivial.

\item  $F$ is unramified at $3$.

 \end{itemize}

\end{lem}

\begin{proof} First we construct the totally real subfield $F^{+}$ of $F$, by means of \cite[Lemma A.2.1 in the appendix]{BLGGT}. Namely,
we require: (1) $F^{+}/\mathbb{Q}$ is linearly disjoint from $\mathbb{Q}_{\text{avoid}}/\mathbb{Q}$; (2) for all places $u$ of $F^{+}$ above $\ell$, $F^{+}_u$ contains $\zeta_{\ell}$ and $\overline{r}_{\ell, \iota}(\pi')\vert_{G_{F^{+}_u}}$ is trivial; (3) for all $p\in \mathcal{L}$, for all places $w$ of $F^{+}$ above $p$, $\overline{r}_{\ell, \iota}(\pi')\vert_{G_{F^{+}_w}}$ is trivial and the cardinality $\kappa_p$ of the residue field of $\mathcal{O}_F^{+}/w$ satisfies $\kappa_p\equiv 1 \pmod{\ell}$; (4) for all $p\in N_{\pi'}\setminus\mathcal{L}$, for all places $w$ of $F^{+}$ above $p$, $\overline{r}_{\ell, \iota}(\pi')\vert_{G_{F^{+}_w}}$ is trivial  and (5) for all places $w$ of $F^{+}$ above $3$, $F^{+}_w/\mathbb{Q}_3$ is trivial.

From this point, we can argue as in Corollary A.2.3 of loc.~cit. We pick an imaginary quadratic field $K$ such that, if we define
$F:=K\cdot F^{+}$, $F$ is linearly disjoint from $\mathbb{Q}_{\text{avoid}}$, and all primes in the set $S\cup\{3\}$ are split in $K/\mathbb{Q}$.
It is clear that such a $K$ exists, for the first condition it is enough to pick $K$ linearly disjoint from $F^{+}\cdot \mathbb{Q}_{\text{avoid}}$.
Then, the primes in $S\cup \{3\}$ are split in $F/F^{+}$, and all the required local conditions are satisfied over $F$ since we know that they hold over
$F^{+}$. \end{proof}

We consider the base change $\pi'_F$ to
the field $F$ of the automorphic representation $\pi'$. $(\pi'_F, \mu)$ is a polarised automorphic representation for a certain algebraic character $\mu$. The next step is to apply
Theorem 4.4.1 of \cite{BLGGT} to the compatible system attached to $(\pi'_F, \mu)$,
to produce an automorphic lifiting of $\overline{r}_{\ell, \iota}(\pi'_F)$ with suitable level and weight. Note that the prime $\ell$ satisfies that $\ell>2(2n+1)$ (with $2n$ equal to the dimension
of $\overline{r}_{\ell, \iota}(\pi'_F)$),
$\zeta_{\ell}\not \in F$,  all primes $\lambda\vert \ell$ of $F$ are split in
$F/F^{+}$.

Label each
possible combination of  Jordan blocks in $\mathrm{GL}_{2n}$, indexed by a  finite set $I$, and for each $i\in I$, choose a different prime $p_i>\sup\{C_0, p: p\vert N_{\pi'}\}$ which is split in $F/F^{+}$, and such that: $\overline{r}_{\ell,\iota}(\pi'_F)|_{G_{F_{v}}}$ is trivial at all places $v$ of $F$ above $p_i$, and  the cardinality $\kappa_{p_i}$ of the residue field of $\mathcal{O}_F/v$ satisfies $\kappa_{p_i}\equiv 1 \pmod{\ell}$. (The existence of such primes is guaranteed by Chebotarev Density Theorem applied over $\mathbb{Q}$).

Let $\mathcal{Z}:=\{p_i: i\in I\}$, and define
\begin{equation}\label{eq:setS}S:=\{\ell\}\cup \mathcal{Z}\cup \mathcal{L}\cup \{p: p\vert N_{\pi'}\}\end{equation}

Note that  all $\lambda\vert \ell$ are  split in $F/F^+$, and that, for all $v\vert p\in S$, $v$ is split in $F/F^{+}$. For each $p\in S$, we choose a place $\tilde{v}$ above $p$.

The residual representation
$\overline{r}=\overline{r}_{\ell, \iota}(\pi'_F)$ satisfies that $(\overline{r}, \overline{\mu})$ is an automorphic, potentially diagonalisable polarised representation
mod $\ell$, unramified outside $S$, and $\overline{r}\vert_{G_{F(\zeta_{\ell})}}$ is absolutely irreducible.

We prescribe the following local conditions:

\begin{itemize}

\item At $v\vert \ell$,  we prescribe the ordinary, crystalline lift $\rho_v=1 \oplus \varepsilon_{\ell}^{-a} \oplus \cdots \oplus \varepsilon_{\ell}^{(1-2n)a}$, where
$a$ is some chosen integer greater than $C_0$ and divisible by $2n$.

\item At $v\vert p\in \mathcal{L}$, we prescribe $\pi''_v$ to be an unramified twist of the Steinberg representation of $\mathrm{GL}_{2n}(F_v)$. Via the local
Langlands correspondence, we are prescribing a  lift $\rho_v$ which is unipotent with a unique Jordan block.

\item At $v\vert p_i$ for $i\in I$, we choose the lift $\rho_v$ as an unipotent with the Jordan block decomposition corresponding to $i\in I$ (this lift exists because we asked $\overline{r}_{\ell,\iota}(\pi'_F)|_{G_{F_{v}}}$ to be trivial).

\item At $v\in N_{\pi'}\setminus \mathcal{L}$ we choose the lift $\rho_v=\mathrm{Id}$ (recall that by definition of $F$, $\overline{r}_{\ell, \iota}(\pi')\vert_{G_{F_v}}$ is trivial).

\end{itemize}

By Theorem 4.4.1 of \cite{BLGGT}, there exists a RACP automorphic representation $(\pi'', \mu)$ of $\mathrm{GL}_{2n}(\mathbb{A}_F)$
such that the compatible system $\{r_{p}(\pi'')\}$ satisfies:

\begin{itemize}

 \item The Hodge-Tate weights of $\{r_{p}(\pi'')\}$ with respect to $\tau\colon F \hookrightarrow \overline{\mathbb{Q}}_\ell$ are
 $\{0,a,\ldots,(2n-1)a\}$. In particular, they do not depend on $\tau$ and they are $C_0$-very spread.

\item The primes in $\mathcal{L}$ are Steinberg in $\{r_{p}(\pi'')\}$.

 \item At the primes of $\mathcal{Z}$, $\{r_{p}(\pi'')\}$ has unipotent ramification with all possible Jordan block decompositions.

 \item $\{r_{p}(\pi'')\}$ is unramified outside the primes of  $\mathcal{L}\cup \mathcal{Z}$. In particular, $\{\det\circ\; \overline{r}_{p}(\pi'')\}$ is unramified outside the
residue characteristic. In fact, for any prime $t$ not dividing the level of $\pi''$, 
 the determinant of $\overline{r}_{t, \iota}(\pi'')$ equals $\overline{\varepsilon}_t^k \cdot \alpha$, where $\alpha$ is an unramified character of 
 $G_F$, $\overline{\varepsilon}_t$ denotes the mod $t$ cyclotomic character, and $k$ is an integer divisible by $2n$.

\end{itemize}

 \begin{rem}\label{rem:level_raising}
For a detailed exposition of the level raising at the primes in $\mathcal{L}$ and $\mathcal{Z}$, we refer to \cite[Prop. 5.2]{ClozelThorneI} and
\cite[Section 6]{ClozelThorneII}, where a similar situation is considered.
\end{rem}

\begin{thm} Let $\mathcal{P}$ be the set of rational primes that are totally split in $F$. Then
 the compatible system $\{r_p(\pi'')\}$ described above satisfies that
 for each $p\in \mathcal{P}$ except for a density zero set of rational primes, the image of the residual representation $\overline{r}_{p, \iota}(\pi)$ contains $\SL_{2n}(\mathbb{F}_p)$ or $SU_{2n}(\mathbb{F}_p)$.
\end{thm}

\begin{proof}
The representations $r_{p, \iota}(\pi'')$ forming the compatible system are semisimple by Theorem 2.1.1 of \cite{BLGGT}. In particular, their restrictions to decomposition groups are also semisimple. Moreover, for each prime $\ell_i\in \mathcal{L}$ and each $v\vert \ell_i$, the restriction $r_{p, \iota}(\pi'')|_{G_{F_v}}$ is also unipotent with a unique Jordan block (by our prescription of local conditions), and so, it must be irreducible. Hence, we obtain that if $p\not=\ell_i$ for some $\ell_i\in\mathcal{L}$, then $r_{p, \iota}(\pi'')$ is irreducible.

Since the set $\mathcal{L}$ has more than one element,
 this ensures (absolute) irreducibility of $r_{p, \iota}(\pi'')$ for all primes $p$. By Proposition 5.3.2 of \cite{BLGGT}, we have that, for a density one set $\mathcal{M}$ of rational primes, the corresponding residual representation $\overline{r}_{p, \iota}(\pi'')\vert_{F(\zeta_{p})}$ is irreducible. This is
one of the conditions required to apply Theorem 4.4.1 of \cite{BLGGT}. By removing a finite set of primes from $\mathcal{M}$, we can assume that for all primes in $\mathcal{M}$, the corresponding Galois representation is Fontaine-Laffaille (which ensures potential diagonalisability, which is also required in Theorem 4.4.1 of loc.~cit.), $p>2(2n+1)$ and $\zeta_{p}\not\in F$. Let $\mathcal{M}'$ be the subset of primes of $\mathcal{M}$ which are totally split in $F/\mathbb{Q}$.

We perform now a level lowering argument (cf.~Lemma 5.2 in \cite{ArDiShWi}) to show that there is at most a finite subset $P$ of $\mathcal{M}'$ such that the unipotent ramification of $r_{p, \iota}(\pi'')$ locally at
$p_i\in \mathcal{Z}$ for $p\in P$ becomes smaller in $\overline{r}_{p, \iota}(\pi'')$.
 We detail this argument. Let $p_i\in \mathcal{Z}$ correspond to the Jordan block decomposition $J_i$, i.e.~the image of the chosen lift $\rho_{p_i}$ at the prime $p_i$ in the second
 preliminary step corresponds to $J_i$. Denote by
 $\overline{J}_i$ the Jordan block decomposition in $\GL_{2n}(\overline{\mathbb{F}}_{p})$ which has the same block structure as $J_i$.

 Assume that there exist infinitely many primes $p\in \mathcal{M}'$  such that  the Jordan block decomposition of $\overline{\rho}_{p, \iota}(I_{p_i})$ is strictly contained in $\overline{J}_i$, that is, some Jordan block in the decomposition of ${\rho}_{p, \iota}(I_{p_i})$ gets reduced into several blocks. We can
 assume, without loss of generality, that these reduced decompositions coincide for infinitely  many primes $p$. For each such prime, we prescribe a minimal lift locally at $p_i$ and we may apply Theorem 4.4.1 of \cite{BLGGT}
 to produce another automorphic representation, say $\pi^{(p)}$, such that $\{r_\ell(\pi'')\}$ and $\{r_\ell(\pi^{(p)})\}$ are connected at all primes different from $p_i$,
 but with different types at primes $v\vert p_i$.
 Since the new family of automorphic representations $\{\pi^{(p)}\}_{p}$ have the same infinitesimal character at $\infty$, fixed ramification set and fixed types at ramified primes,
 the finiteness result of Harish-Chandra imply that one of them occurs infinitely often, that is, there exists a representation $\Pi$ and  an infinite set of primes $p$ such that $\Pi=\pi^{(p)}$.
 As a consequence, we obtain infinitely many congruence conditions between the compatible system $\{r_p(\pi'')\}$ and $\{r_p(\Pi)\}$, implying that $r_{p, \iota}(\pi'')$ and $r_{p, \iota}(\Pi)$ are
 isomorphic. But their types at $p_i$ are different by construction, which is a contradiction.

 Thus, we have that, except for finitely many primes $p\in \mathcal{M}'$, the image of $\overline{r}_{p, \iota}(\pi'')$ contains unipotent elements of each possible Jordan form. Moreover, taking into account that $p>2(2n+1)$, and noting that the irreducibility of the image ensures that it does not contain a non-trivial normal subgroup consisting of unipotent elements, we conclude that Theorem \ref{thm:Guralnick} applies, yielding the desired result.

\end{proof}

\begin{rem}\label{rem:From_Pi1_to_Pi2} Note that by ALT-PD, the automorphy of $\{(\rho(f)\otimes r(\pi''))_{p}\}$ implies the automorphy of $\{(\rho(f)\otimes r(\pi'))_{p}\}$, thus we may replace $\pi'$ by $\pi''$ for our purposes. Let us explain that all conditions required for the application of this theorem are satisfied in this congruence modulo $\ell$. 
First of all, the prime $\ell$ satisfies $\ell > 2 (4n +1)$. Concerning the residual image, over $\mathbb{Q}(\zeta_\ell)$ it is absolutely 
irreducible  as noted in the list of properties of $\pi'$ (item 2), and then using the first listed property of the field $F$ we deduce that 
after restricting to the absolute Galois group of $F(\zeta_\ell)$ the residual representation is still absolutely irreducible.  
Both tensor products are regular: for the one containing $\pi'$ this was already noted in the list of properties of $\pi'$ (item 3), and for the other one it follows from the choice of Hodge-Tate weights in the local parameters at primes dividing $\ell$ in the construction of $\pi''$. Finally, both tensor products are pot. diag. at primes above $\ell$: again, for the one containing $\pi'$ this was noted (we noted that it was Fontaine-Laffaille over $\mathbb{Q}$, and  pot. diag. is preserved by base change) and for the other one it follows from the local parameters at primes dividing $\ell$ chosen when constructing $\pi''$ (together with the fact that the $2$-dimensional component is Fontaine-Laffaille over $\mathbb{Q}$).
\end{rem}

\subsection{Preliminary step 3}\label{subsection:PreliminaryStep3}

In this step we want to perform a level raising in order to  add a good monomial prime $q$ to the compatible system. That is to say, we
will introduce a specific type of ramification at a well-chosen prime $q$, that will
ensure that, at all the steps in the chain, we will have residually large image in the $2n$-dimensional component. The concept of good monomial prime, which was
first used in \cite{KLS08}, has its roots in the  good dihedral prime appearing in the proof of Serre's Conjecture. 

We need to introduce some terminology. Let $t, q$ be two different odd primes, satisfying that the order of $q$ mod $t$ is $2n$. If $\mathbb{Q}_{q^{2n}}$
denotes the maximal unramified extension of $\mathbb{Q}_q$ of degree $2n$, and $\ell$ denotes an odd prime different from $t$ and $q$,
we may consider a character $\chi_q: {\mathbb{Q}}_{q^{2n}}^{\times}\simeq \mu_{q^{2n}-1}\times U_1\times q^{\mathbb{Z}}\rightarrow \overline{\mathbb{Q}}_{\ell}^ {\times}$
such that the following properties hold:
\begin{itemize}
 \item the order of $\chi_q$ is $2t$;
 \item The restriction of $\chi_q$ to $\mu_{q^{2n}-1}\times U_1$ has order $t$;
 \item $\chi_q(q)=-1$.
\end{itemize}

Via local class field theory, we may identify $\chi_q$ as a character of the absolute Galois group $G_{\mathbb{Q}_{q^{2n}}}$.  Then the representation
$\rho_q:=\Ind^{\mathbb{Q}_q}_{\mathbb{Q}_{q^{2n}}}(\chi_q):G_{\mathbb{Q}_q}\rightarrow \mathrm{GL}_{2n}(\overline{\mathbb{Q}}_{\ell})$ is irreducible, and, for every $\ell \neq q, t$ its mod $\ell$
reduction $\overline{\rho}_q:G_{\mathbb{Q}_q}\rightarrow \mathrm{GL}_{2n}(\overline{\mathbb{F}}_{\ell})$ has image equal to a $(2n, t)$-group (cf.~Section \ref{sec:subgroups}, or 
\cite[Section 2]{KLS08} for the definition
of $(n, p)$-group).

This local representation can be exploited to ensure large image in a compatible system of Galois representations. Namely, let $(\rho_{\lambda})_{\lambda}$ be
a compatible system of $\lambda$-adic Galois representations (as in \cite[Section 5]{BLGGT}), where $\lambda$ runs through the finite places of a number field $M$,
and assume that, for all $\lambda \nmid q$, the restriction of $\rho_{\lambda}$ to a decomposition group at $q$ equals $\rho_q$.
For our applications, we want to ensure that, for all primes $\ell$ up to a certain bound, the image of $\overline{\rho}_{\lambda}$ is \emph{large} (i.e., satisfies the conclusion of
Theorem \ref{thm:KLS}).

\begin{prop}\label{prop:GoodMonomialPrime} Let $n, W\geq 2$ be natural numbers,
and let $d(2n)$, $p(2n)$ be the constants from Theorem \ref{thm:KLS}. Let $F$ be a number field, and $K/F$ be the compositum of all field extensions of $F$ of
degree at most $d(2n)+1$ which are unramified outside the primes above all $p\leq W$ and $\infty$.
Let $(t, q)$ be a couple of prime numbers satisfying that $q>W$ splits completely in $K/\mathbb{Q}$, $t\equiv 1 \pmod{2n}$, $t>\max\{d(2n)+1, p(2n), W\}$ and the order of $q$ mod $t$ is exactly $2n$.

 Let $(\rho_{\lambda})_{\lambda}$ be a compatible
 system of Galois representations $\rho_{\lambda}:G_{F}\rightarrow \GL_{2n}(\overline{\mathbb{Q}}_{\lambda})$, indexed by the set of primes $\lambda$ of a certain
 number field $M$, with ramification set $S$.
 Assume that all primes $v$ in $S$, say $v\vert p$, satisfy $p\leq W$ except when $p=q$. Furthermore, assume that for all places $\lambda$ of $M$ with $\lambda \nmid q$,
 \begin{equation*}
  \mathrm{Res}^{G_{F}}_{G_{F_{\mathfrak{q}_i}}}(\rho_{\lambda})=\mathrm{Ind}^{G_{\mathbb{Q}_q}}_{G_{\mathbb{Q}_{q^{2n}}}}(\chi_q)\otimes \alpha,
 \end{equation*}
 where:
 \begin{itemize}
  \item  $\mathfrak{q}_i$ is any prime above $q$, and we identify $F_{\mathfrak{q}_i}$ with $\mathbb{Q}_q$ via a natural isomorphism;
 \item $\chi_q:G_{\mathbb{Q}_{q^{2n}}}\rightarrow \overline{\mathbb{Z}}^{\times}$ satisfies that, via the embedding
$\overline{\mathbb{Z}}^\times\rightarrow \overline{\mathbb{Q}}^{\times}_{\ell}$ given by $\lambda$, is a character as described above;
 \item $\alpha:G_{\mathbb{Q}_q}\rightarrow \overline{M}^{\times}_{\lambda}$ is an unramified character.
\end{itemize}

Then, for all $\ell\leq W$ prime and $\lambda\vert \ell$,
 the residual image  $\overline{\rho}_{\lambda}(G_{F})\subset \GL_{2n}(\overline{\mathbb{F}}_{\ell})$ satisfies the conclusion of Theorem \ref{thm:KLS}.
\end{prop}

\begin{proof} Let us fix a place $\lambda\vert \ell$ with $\ell\leq W$, and let $L$ be the field cut out by $\overline{\rho}_{\lambda}$, so that $\Gamma:=\mathrm{Im}\overline{\rho}_{\lambda}\simeq \Gal(L/F)$.
Since $\ell\leq W <q$, we have that $\mathrm{Res}^{G_{F}}_{G_{F_{\mathfrak{q}_i}}}(\rho_{\lambda})=\mathrm{Ind}^{G_{\mathbb{Q}_q}}_{G_{\mathbb{Q}_{q^{2n}}}}(\chi_q)\otimes \alpha$, and in
 particular $\overline{\rho}_{\lambda}(G_{F_{\mathfrak{q}_i}})\subset \Gamma$ is a $(2n, t)$-group.
We will show that $\overline{\rho}_{\lambda}(G_{F_{\mathfrak{q}_i}})\subseteq \Gamma^{d(2n)+1}$ (where $\Gamma^{d(2n)+1}$ is the intersection of all normal groups of $\Gamma$ of index at most $d(2n)+1$),
and the conclusion will follow from Theorem \ref{thm:KLS}, taking $d=d(2n)+1> d(2n)$, $p=t>p(2n)$, and $\ell$, which is different from $p$  because $\ell\leq W<t=p$.

Let us fix a normal subgroup $H$ of $\Gamma$ of index at most $d(2n)+1$; it suffices to show that $\overline{\rho}_{\lambda}\vert_{G_{F_{\mathfrak{q}_i}}}\subseteq H$.
Let $L_H/F$ be the extension afforded by $H$; we have that $L_H/F$ is a finite Galois extension of degree at most $d(2n)+1$. Since $\overline{\rho}_{\lambda}$
ramifies only at the primes of $S$, we have that $L_H/F$ is unramified outside $S\cup\{\infty\}$. Moreover, since $\overline{\rho}_\lambda\vert_{I_{F_{\mathfrak{q}_i}}}$ has order $t>d(2n)+1$, we
can conclude that $L_H/F$ is unramified at $\mathfrak{q}_i$ for all $\mathfrak{q}_i\vert q$. Therefore, the extension $L_H/F$ is a subextension of $K/F$. Since $q$ is completely split
in $K/\mathbb{Q}$, it is also completely split in $L_H/F$. Thus the decomposition group of $L/L_H$ at each prime above $q$ can be identified with the decomposition group of $L/F$ at
each prime above $q$. Thus $\overline{\rho}_{\lambda}(G_{F_{\mathfrak{q}_i}})\hookrightarrow H$, as we wanted to prove.

\end{proof}

\begin{rem}
For each $n$, $W$ as above, the existence of the couple $(q, t)$ is guaranteed by \cite[Lemma 3.4]{KLS08}
\end{rem}

Given the automorphic representation $\pi''$ produced in the second preliminary step, we want to apply a level raising result to introduce a
good monomial prime $q$ that will ensure that the residual image is large at (almost) all characteristics that we need to work with. In order
to do this, we need to take  into account some quantities related to the representation $\pi''$, and some further quantities related to the $2$-dimensional skeleton chain.
Recall that in Remark \ref{rem:ExtendedChain} we had defined the extended $2$-dimensional skeleton chain $\widetilde{\mathcal{C}}$. Moreover, before the first
preliminary step in Section \ref{sec:PreliminarySteps}, we also
fixed a constant $\widehat{k}$. We further define

\begin{itemize}

\item $C_1:=\sup(\{k_{\mathrm{max}}, w_{\mathrm{max}}, 2(4n+1), (2n-1)\widehat{k}, 2n C_0\}\cup\mathcal{P}_1\cup\mathcal{P}_2\cup \{n + \ell^2: \ell\in \mathcal{L}\}\cup \mathcal{Z}) + 1$.

\end{itemize}

Let $B>C_1$ be such that there exist at least 15 auxiliary primes $t_i$ between $B$ and $2B$ which are totally split in $F/\mathbb{Q}$.
For future reference, denote this set of auxiliary primes as $\mathcal{T}_{\mathrm{aux}}$.
Let us also define $\mathrm{Comp}$ as the compositum, inside a fixed algebraic closure of $F$, of all extensions of $F$ which are of degree at most the maximum
between $d(2n)+1$ and $p(2n)$ and ramified only above the primes $p\leq 2B$ and $\infty$.

We need to specify a couple of primes $(t, q)$ in order to add the good monomial prime to $\pi''$ with respect to the bound $W=2B$. More precisely,
we want to obtain another automorphic representation $\pi'''$ such that, for $\ell\nmid tq$, locally at a prime $\mathfrak{q}$ above $q$, the representation $r_{\ell, \iota}(\pi''')$ contains a $(2n, t)$-group, while
at all primes $v$ of $F$ which are not above
 $q$, $r_{\ell, \iota}(\pi'')\vert_{G_{F_v}}$ and $r_{\ell,\iota}(\pi''')\vert_{G_{F_v}}$ are connected. This level raising will be achieved through an application
of Theorem 4.4.1 of \cite{BLGGT}, via a congruence modulo a prime $\mathfrak{t}$ above $t$. In order to apply this theorem,
we need that $t$ satisfies the following hypotheses:

\begin{enumerate}
 \item $t>2(4n+1)$ (in fact, $ t > 2 (2n +1)$) would be enough, but we impose a stronger condition that will be required later);
 \item $\zeta_t\not\in F$;
 \item $t$ is completely split in $F/\mathbb{Q}$;
 \item $r_{t, \iota}(\pi)$ is potentially diagonalisable (this will be implied by inequality 6 below);
 \item $\overline{r}_{t, \iota}(\pi)\vert_{G_{F(\zeta_t)}}$ is irreducible, and, moreover, it satisfies the stronger condition in the conclusion of Theorem 5.4.
\end{enumerate}

Besides, we will require several additional conditions on the size of $t$, that will be used later in our constructions:

\begin{enumerate}
\setcounter{enumi}{5}
\item $t> 2C_1$.

\item $t>\max\{d(2n)+1, p(2n)\}$ from Theorem \ref{thm:KLS}.

\item $t$ big enough so that if, for any $p<2B$, for $G\in \{\SL_{2n}(\mathbb{F}_{p^r}), \mathrm{SU}_{2n}(\mathbb{F}_{p^r}),$ $\Sp_{2n}(\mathbb{F}_{p^r}), \Omega^{\pm}_{2n}(\mathbb{F}_{p^r})\}$
satisfying that $G$ contains an element of order $t$, then the exponent $r$ is big enough so that $G$ is adequate$+$. (cf.~Theorem \ref{thm:adequate} and Proposition \ref{prop:adequate}).

\item The cardinality of $\mathrm{PSL}_{2n}(\mathbb{F}_t)$ is greater than the cardinality of the Galois group of $\mathrm{Comp}$ over $F$.
\end{enumerate}

Furthermore, we will ask a good behaviour of the prime $t$ with respect to $F$, $\mathrm{Comp}$, and $\pi''$, namely

\begin{enumerate}
\setcounter{enumi}{9}
 \item $t$ is completely split in $H/\mathbb{Q}$, where $H$ is the Hilbert class field of $F$;
 \item $t$ is completely split in the coefficient field of $\pi''$ over $\mathbb{Q}$.
 \item $t$ is unramified in $\mathrm{Comp}/F$.
\end{enumerate}

Finally, we need the condition
\begin{enumerate}
\setcounter{enumi}{12}
 \item $t\equiv 1 \pmod{2n}$ in order to apply Proposition \ref{prop:GoodMonomialPrime}.
\end{enumerate}
The existence of a prime $t$ satisfying all these conditions follows from the following lemma, which we label
for future reference.

\begin{lem}\label{lem:Existst} There exists a positive density of primes $t$ satisfying the 12 conditions listed above.

\end{lem}

\begin{proof} The proof follows from Chebotarev's Density Theorem, together with Lemma 5.2.

\end{proof}

\begin{rem}
Condition 11 in Lemma \ref{lem:Existst} implies that the image of $\overline{r}(\pi''')_{t, \iota}$ will be contained in $\GL_{2n}(\mathbb{F}_t)$. This property will be
 exploited in Lemma \ref{lem:ExistenceOfq} below.
\end{rem}

Next we want to specify the prime $q$, depending of a given prime $t$ as in the conclusion of Lemma \ref{lem:Existst}. Let us collect all conditions we want to impose to it (cf.~Proposition \ref{prop:GoodMonomialPrime}).

\begin{itemize}
 \item $q> 2B$;
 \item $q$ completely split in $\mathrm{Comp}/\mathbb{Q}$;
 \item the order of $q$ mod $t$ is exactly $2n$;
\end{itemize}

Note that this last condition can be reformulated as follows: $\mathrm{Frob}_q\in \Gal(\mathbb{F}_q(\zeta_t)/\mathbb{F}_q)$ has order exactly $2n$.
In fact, since $q$ will be completely split in $\mathrm{Comp}/\mathbb{Q}$, it will also be completely split in $F/\mathbb{Q}$, thus the residue field
of $F$ modulo any prime $\mathfrak{q}$ above $q$ will be $\mathbb{F}_q$; denoting by $\kappa(F_{\mathfrak{q}})$ the residue field of $F_{\mathfrak{q}}$,
we get that the last condition is equivalent to

\begin{itemize}
 \item For any prime $\mathfrak{q}$ above $q$, we have that $\mathrm{Frob}_{\mathfrak{q}}\in \Gal(\kappa(F_{\mathfrak{q}})(\zeta_t)/\kappa(F_{\mathfrak{q}}))$ has
 order exactly $2n$.
\end{itemize}

Since we want that the representation $r(\pi''')_{t, \iota}$ connects locally at $q$ with a lift providing a $(t, 2n)$-group,
we need a compatibility relation between the reduction $\overline{r}(\pi'')_{t, \iota}$ and the reduction of such a local parameter at $q$ modulo $t$. We specify this
relation. Let $M$ be the subfield of $\overline{F}$ cut out by the residual representation $\overline{r}(\pi'')_{t, \iota}$, and $L$ the subfield of $M$ cut out by the determinant
of $\overline{r}(\pi'')_{t, \iota}$. The condition we need is that
$q$ is completely split in $L/F$ and the order of the Frobenius element
at $q$ in $\Gal(M/L)$ equals $2n$.

\begin{lem}\label{lem:ExistenceOfq} Let $t$ be a prime number $t \equiv 1 \Mod {2n}$. Let $F$ be a number field such that $t$ is totally split in $F/\mathbb{Q}$ (so in particular $\zeta_t\not\in F$). Let $\overline{\rho}:G_F\rightarrow \GL_{2n}(\mathbb{F}_{t})$ be a continuous representation
such that the determinant $\det\overline{\rho}=\alpha\cdot \overline{\varepsilon}_t^k$, where $k$ is a multiple of $2n$, $\alpha:G_F\rightarrow \mathbb{F}_t^{\times}$ is an
unramified character and $\overline{\varepsilon}_t$ denotes the mod $t$ cyclotomic character.
Assume further that the image of $\overline{\rho}(G_F)$ contains $\SL_{2n}(\mathbb{F}_t)$ and that $t$ is totally split in the Hilbert class field $H$ of $F$.
Then for  any  Galois extension $F_{\text{avoid}}/F$ of order less than the cardinality of $\PSL_{2n}(\mathbb{F}_t)$, which is unramified at $t$,
there exists a positive density of rational primes $q$, completely split in $F/\mathbb{Q}$, satisfying the following: There is a prime $\mathfrak{q}\vert q$ of $F$ such that:
 \begin{enumerate}
 \item $\mathfrak{q}$ is completely split in $F_{\text{avoid}}/F$;
\item $\mathrm{Frob}_\mathfrak{q}\in Gal(F(\zeta_t)/F)$ has order exactly $2n$;
\item $\mathrm{Frob}_\mathfrak{q}\in Gal(M_{\overline{\rho}}/F)$ has order exactly $4n$, where $M_{\overline{\rho}}$ is the subfield
of $\overline{F}$ cut out by the representation $\overline{\rho}$.

 \end{enumerate}

\end{lem}

\begin{proof}
The strategy of the proof will be to prescribe the image of $\mathrm{Frob}_\mathfrak {q}\in G_F$ in several,
linearly disjoint finite extensions of $F$, by means of Chebotarev's Density Theorem.
Recall that the condition that $\mathfrak{q}$ is completely split in a finite Galois extension $F'/F$
amounts to asking that the conjugacy class of $\mathrm{Frob}_\mathfrak{q}$ in $\Gal(F'/F)$ is the class of identity.

First of all, note that we may assume, without loss of generality, that $\det\overline{\rho}=\overline{\varepsilon}_t^k$. Indeed, if this is not
the case, we replace $F$ by $H$ and $F_{\text{avoid}}$ by the compositum $F_{\text{avoid}}H$. On the one hand, let us check that all conditions of the Lemma
are still satisfied. Since $t$  is totally split in $H/\mathbb{Q}$, $\zeta_t\not\in H$. Moreover, $\overline{\rho}(G_{H})$ still contains the
quasi-simple group $\SL_{2n}(\mathbb{F}_{t})$. The extension $F_{\text{avoid}}H/H$ is unramified at $t$ and the order of $F_{\text{avoid}}H/H$
is less than the cardinality of $\PSL_{2n}(\mathbb{F}_t)$. On the other hand, assume we produce a positive density of rational primes $q$, completely split in $H/\mathbb{Q}$, such that
for some prime $\mathfrak{q}\vert q$, then
$\mathfrak{q}$ is completely split in $F_{\text{avoid}}H/H$ and
 the two conditions on $\mathrm{Frob}_\mathfrak{q}$ hold.
Then denoting by $\mathfrak{q}$ also the prime of $F$ below $\mathfrak{q}$, we have that $\mathfrak{q}$ is completely split in  $F_{\text{avoid}}/F$ (because it is completely split in $H/\mathbb{Q}$
and $F_{\text{avoid}}H/H$. The fact that $H/F$ is unramified at $t$ implies
that $F(\zeta_t)/F$ and $H/F$ are linearly disjoint, so condition 2 holds. For the third condition, note that $q$ is totally split in $H/F$,
hence the decomposition groups at $\mathfrak{q}$ of $M_{\overline{\rho}}H/H$ and $M_{\overline{\rho}}/F$
are isomorphic.

Since $t\equiv 1 \pmod{2n}$, there exists some $\sigma\in \Gal(F(\zeta_t)/F)$
of order exactly $2n$. Indeed, this follows because $t$ is also totally split in $F$, so $F$ is linearly disjoint from $\mathbb{Q}(\zeta_t)$ and $\Gal(F(\zeta_t)/F)\simeq\Gal(\mathbb{Q}(\zeta_t)/\mathbb{Q})$. We prescribe that $\mathrm{Frob}_{\mathfrak{q}}$ lies in the class of $\sigma$ (\textbf{Condition 1}).

Next, let $L\subset M_{\overline{\rho}}$ be the subfield cut out by $\det\overline{\rho}=\overline{\varepsilon}_t^k$. Then $L/F$ is a Galois
extension, and we prescribe that $\mathrm{Frob}_\mathfrak{q}=1$, that is to say, $\mathfrak{q}$ completely split in $L/F$ (\textbf{Condition 2}).
Note that this condition is compatible with the previous one. Indeed, observe that $L\subset F(\zeta_t)$ by definition, so we need to see that if we prescribe that $\mathrm{Frob}_{\mathfrak{q}}$ has order
exactly $2n$ in the extension $F(\zeta_t)/F$, then  $\mathrm{Frob}_{\mathfrak{q}}$  has order $1$ in $L/F$.
Indeed, assume that the order of $\mathrm{Frob}_{\mathfrak{q}}$ has order
exactly $2n$ in the extension $F(\zeta_t)/F$. On the one hand, $q^{2n}\equiv 1\pmod t$ since $q$ is completely split in $F$.
On the other hand, $2n\vert k$ implies that $\mathrm{Frob}_{\mathfrak{q}}\in \ker\varepsilon_t^k$. Thus $\mathrm{Frob}_{\mathfrak{q}}$
acts trivially on $L$, as desired.

Thus, conditions 1 and 2 allow us to prescribe a behaviour of $\mathrm{Frob}_\mathfrak{q}$ in the extension $L(\zeta_t)/L$ and $L/F$.
Next we want to prescribe the class of $\mathrm{Frob}_\mathfrak{q}$ in the extension $M_{\overline{\rho}}/L$. Let us check that this extension is
linearly disjoint from $L(\zeta_t)/L$. Indeed, if there were some non-trivial intersection, say $F'\neq L$, then $\Gal(M_{\overline{\rho}}/F')$ would be a proper normal subgroup of
$\Gal(M_{\overline{\rho}}/L)\simeq \SL_{2n}(\mathbb{F}_t)$. But the only proper normal subgroup is $\{\pm 1\}$, which would yield $\Gal(F'/L)\simeq \PSL_{2n}(\mathbb{F}_{t})$,
a contradiction since $L(\zeta_t)/L$ is cyclic. This reasoning shows that we may prescribe any conjugacy class for $\mathrm{Frob}_\mathfrak{q}$ in
$\Gal(M_{\overline{\rho}}/L)$, and it will be compatible with conditions 1 and 2.

Let us choose the element $\sigma$ corresponding to the permutation matrix  \begin{equation*}\begin{pmatrix}0  & 0 & \cdots & 0 & -1\\
 1 & 0 & \cdots & 0& 0\\
 \vdots & \vdots & \ddots & 0& 0\\
 0 & 0 &  \cdots & 1 & 0
  \end{pmatrix}\end{equation*}

This element has order exactly $4n$. We prescribe that
 the class of $\mathrm{Frob}_\mathfrak{q}$ in $\Gal(M_{\overline{\rho}}/L)$ belongs to the class of this element (\textbf{Condition 3}).

Conditions 1, 2 and 3 ensure that $\mathrm{Frob}_{\mathfrak{q}}$ has the desired behaviour in the extension $\Gal(M_{\overline{\rho}}(\zeta_t)/F)$.
It remains to see that this behaviour is compatible with the condition of $\mathfrak{q}$ being completely split in the extension $F_{\text{avoid}}/F$.

If $M_{\overline{\rho}}(\zeta_t)/F$ is linearly disjoint with $F_{\text{avoid}}/F$, we are done. Otherwise, assume that $F_1=M_{\overline{\rho}}(\zeta_t)\cap F_{\text{avoid}}$ is different
from $F$. Recall that $\Gal(M_{\overline{\rho}}(\zeta_t)/F)$ is a subgroup of
$\Gal(M_{\overline{\rho}}/F)\times \Gal(F(\zeta_t)/F)\subset \GL_{2n}(\mathbb{F}_{t})\times \Gal(F(\zeta_t)/F)$. 

Let us consider the quotient map  $\Gal(M_{\overline{\rho}}(\zeta_t)/F)\rightarrow \Gal(F_1/F)$. We have that $\Gal(F_1/F)$ either contains $\SL_{2n}(\mathbb{F}_t)$, $\PSL_{2n}(\mathbb{F}_{t})$, 
or is a direct product of cyclic groups $\mathbb{Z}/r\mathbb{Z}\times \mathbb{Z}/s\mathbb{Z}$, where $s\mid t-1$.
The first two cases cannot happen because $\sharp\Gal(F_1/F)\leq\sharp\Gal(F_{\text{avoid}}/F)< \sharp (\PSL_{2n}(\mathbb{F}_t))$, where the last
inequality holds by hypothesis.

Hence $\Gal(F_1/F)\simeq \mathbb{Z}/r\mathbb{Z}\times \mathbb{Z}/s\mathbb{Z}$, which implies the inclusion $F_1\subset L\cdot F(\zeta_t)$. Since $L\subset F(\zeta_t)$, we deduce that $F_1\subset F(\zeta_t)$. The following diagram summarizes the field inclusions.

\begin{equation*}
\begin{tikzcd}
F_{\text{avoid}} \arrow[rd, no head] \arrow[rdd, no head, bend right] &                                            & F(\zeta_t) \\
                                                               & F_1 \arrow[ru, no head] \arrow[d, no head] &            \\
                                                               & F \arrow[ruu, no head, bend right]         &            \\
                                                               & \mathbb{Q} \arrow[u, no head]              &           
\end{tikzcd}
\end{equation*}
By hypothesis $t$ is unramified in $F_{\text{avoid}}/F$, hence it is unramified in $F_1/F$. On the other hand, $t$ is totally split in $F/\Q$, hence totally ramified in $F(\zeta_t)/F$. 
Consequently, $t$ should also be totally ramified in $F_1/F$, and we reach a contradiction. This proves that that $M_{\overline{\rho}}(\zeta_t)/F$ is linearly disjoint with $F_{\text{avoid}}/F$, concluding the proof.

\end{proof}

We consider a prime $\ell=t$ as in Lemma \ref{lem:Existst}, and for such a $t$ the residue representation
$\overline{r}_{t, \iota}(\pi''):G_F\rightarrow \GL_{2n}(\mathbb{F}_t)$. Note that we can apply Lemma \ref{lem:ExistenceOfq}
to this data, since $\pi''$ was constructed so that the determinant of $\overline{r}_{t, \iota}(\pi'')$ equals $\overline{\varepsilon}_t^k\cdot \alpha$, 
for some unramified character $\alpha$ and a certain exponent $k$ divisible by $2n$ (see the list of conditions just before Remark \ref{rem:level_raising}).

Thus, we can choose a prime $q$ as in Lemma \ref{lem:ExistenceOfq} for the extension $F_{\text{avoid}}=\mathrm{Comp}$. That we meet all the required hypothesis to apply Lemma \ref{lem:ExistenceOfq} follows from the list of conditions Lemma \ref{lem:Existst} makes $t$ satisfy. For example, condition 9 guarantees that $\sharp\Gal(\mathrm{Comp}/F)<\sharp\PSL_{2n}(\mathds{F}_t)$. Finally, we are now ready to apply Theorem 4.4.1 of \cite{BLGGT} and produce the desired automorphic representation $\pi'''$.
The set $S$ consists of the union of $\mathcal{L}$, $\mathcal{Z}$ and all primes of $F$ above $t$ and $q$.  Recall that, by Condition 5, the image of $\overline{r}_{t, \iota}(\pi'')$,
even when restricted to $G_{F(\zeta_t)}$, is irreducible.
For each prime $v$ above $q$, we fix a lifting
$\rho_v:G_{F_v}\rightarrow \GL_{2n}(\mathcal{O}_{\overline{\mathbb{Q}}_t})$ given by the induction from $G_{\mathbb{Q}_{q^{2n}}}$ to $G_{\mathbb{Q}_q}$
of a character $\chi_q$ of order $2t$ as in the statement of Proposition \ref{prop:GoodMonomialPrime}. For every prime $v\vert t$ and $v\in \mathcal{L}\cup \mathcal{Z}$ we propose
the same local parameter $\rho_v$  than the one obtained from $\pi''$.
Note that $t$ is unramified in $F$, and by Condition 6 (i.e. $t>2C_1$),
$t$ is greater than the weights of $\pi''$, so we are in a Fontaine-Laffaille situation, and in particular $\rho_v$ is potentially diagonalisable for
all $v\vert t$. Applying Theorem 4.4.1 of loc. cit. we obtain another automorphic representation $\pi'''$, satisfying the following list
of conditions:

\begin{enumerate}
\item The Hodge-Tate weights of $\pi'''$ agree with those of $\pi''$. That is, given $\tau\colon F \hookrightarrow \overline{\mathbb{Q}}_t$, the weights of $\pi'''$ with respect to $\tau$ are
 $\{0,a,\ldots,(2n-1)a\}$. In particular, they do not depend on $\tau$ and they are $C_0$-very spread. 
 
\item $\pi'''$ is unramified outside of $\{q\}\cup \mathcal{L}\cup\mathcal{Z}$;
\item The primes of $\mathcal{L}$ are Steinberg in $\{r_{p}(\pi''')\}$;
 \item $\{r_{p}(\pi''')\}$ has a good monomial prime at $q$ (i.e., the local parameter at $\mathfrak{q}$ is as described at the beginning of
 this subsection);
 \item For all primes $p<2B$, an application of Theorem \ref{thm:KLS} and \ref{prop:GoodMonomialPrime} yields that the image of $\overline{r}_{p, \iota}(\pi''')$ contains $\SL_{2n}(\mathbb{F}_{p^r})$,
 $\mathrm{SU}_{2n}(\mathbb{F}_{p^r})$, $\mathrm{Sp}_{2n}(\mathbb{F}_{p^r})$ or $\Omega^{\pm}_{2n}(\mathbb{F}_{p^r})$ and is contained in the normaliser of this group, where the exponent $r$  is big 
 enough so that the image is adequate$+$ by Proposition \ref{prop:adequate} (which follows from the size of $t$, see condition (8) in the choice of $t$).
\end{enumerate}

\begin{rem}\label{rem:From_Pi2_To_Pi3} Note that by ALT-PD applied in residual characteristic $t$, the automorphy of $\{\rho_{p}(f)\otimes r_{p}(\pi''')\}$ implies the automorphy of $\{\rho_{p}(f)\otimes r_{p}(\pi'')\}$, thus we may replace $\pi''$ by $\pi'''$ for our purposes. The conditions required to apply this theorem are easy to check. The weights on the two sides of this congruence are the same, and regularity of the side containing $\pi''$ was already noted in the previous subsection. As observed in the previous paragraph, the $t$-adic Galois representations corresponding to $\pi''$ and $\pi'''$ are pot.~diag., and since the same holds for the $2$-dimensional component (because of Condition 6, $t > 2 C_1$, thus it lies in the Fontaine-Laffaille range) the two tensor products representations are pot.~diag. By Condition 1 on $t$ we know that $t > 2 (4n+1)$. Finally, for the residual image, one component has residual image that lies between $\SL_{2n}(\mathbb{F}_t)$ and $\GL_{2n}(\mathbb{F}_t)$ with $n > 1$, and the other component has an image that contains $\SL_{2}(\mathbb{F}_t)$ because of Condition 6 on $t$. Thus the residual image of the tensor product is clearly sufficiently large. We conclude that ALT-PD can be applied.
\end{rem}

\begin{rem}\label{rem:GoodMonomialPrime}
 We constructed the representation $\pi'''$ in such a way that, for all primes $p\leq 2B$, we can ensure that the image of $\overline{r}_{p, \iota}(\pi''')$
 is large. Moreover, any other representation $\pi^{iv}$ which has the same type as $\pi'''$ at $q$, and such that the primes in the ramification set
 of $\pi^{iv}$ are smaller than $2B$, will also have large image residual image at the primes $p\leq 2B$. At some places in Section \ref{sec:Diagrams}, we
 will be forced to modify the representation $\pi'''$ and replace it by another one; as long as we can guarantee that the ramification set of the new representation only contains primes smaller than $2B$,
 we will be able to make use of Proposition \ref{prop:GoodMonomialPrime} to ensure that the image is still large.

\end{rem}

\section{The  $2\times 2n$-dimensional chain}\label{sec:Diagrams}

Recall that in Section \ref{sec:2-dim Chain}, we considered the chain of modular forms
\begin{equation}\label{2ndimchain}
 \mathcal{C}=\{f_1, \dots, f_h\},
\end{equation}
which connects the given modular form $f=f_1$ with a modular form $f_h$ with CM by $\mathbb{Q}(\sqrt{-3})$. 

Recall that we are also given an automorphic representation $\pi$ of  $\GL_{2n}(\mathbb{A}_{\mathbb{Q}})$. Our aim is to prove that the
compatible system of Galois representations $\{(\rho(f)\otimes r(\pi))_{p}\}=\{(\rho_{p}(f)\otimes r_{p}(\pi))\}$ is automorphic. \\
The key idea is to consider a chain of representations
\begin{equation}\label{eq:BigChain}\{\rho(f_1)\otimes r(\pi), \dots, \rho(f_i)\otimes r(\pi^{(j)}),    \dots, \rho(f_h)\otimes r(\pi^{(\textrm{end})}) \},\end{equation}
 such that we can propagate automorphy from the right to the left.

By Remarks \ref{rem:From_Pi_to_Pi1}, \ref{rem:From_Pi1_to_Pi2} and \ref{rem:From_Pi2_To_Pi3}, it follows that if $\{(\rho(f)\otimes r(\pi'''))_{p}\}$ is automorphic, then the compatible system $\{(\rho(f)\otimes r(\pi))_{p}\}$ is also automorphic. Therefore we can assume, without loss of generality, that the leftmost compatible system is $\{(\rho(f_1)\otimes r(\pi'''))_{p}\}$. In particular, it will have dimension $2\times 2n$. Observe that at this point (see previous section), objects are only defined over a suitable solvable extension $F$ of $\mathbb{Q}$.

\begin{rem}
 Recall that, at each step of the chain $\mathcal{C}$, it holds that if $f_i$ and $f_{i+1}$ are two consecutive terms, then there exists a prime $\mathfrak{p}_i \vert p_i$ such that $f_i$ is congruent modulo $\mathfrak{p}_i$ to $f_{i+1}$, a twist $\chi\otimes f_{i+1}$ of $f_{i+1}$ by a finite order character or a Galois conjugate $f^{\sigma}_{i+1}$ of $f_{i+1}$ (cf.~Section \ref{sec:2-dim Chain}). In the corresponding step in the $2\times 2n$-dimensional chain, we want to propagate automorphy from $\{(\rho(f_{i+1})\otimes r(\pi^{(i+1)}))_{p}\}$  to $\{(\rho(f_{i})\otimes r(\pi^{(i)}))_{p}\}$, where $\pi^{(i)}$ and $\pi^{(i+1)}$ are two automorphic representations that are congruent modulo a prime $\mathfrak{P}_i$ in the Galois closure of the compositum of the fields of definition of the compatible systems above $\mathfrak{p}_i$. In order to propagate automorphy, we will apply one of the four ALT presented in Section \ref{sec:ALT}.

If $f_i$ is congruent to $f_{i+1}$, we will have a congruence between the representations  $\rho_{\mathfrak{P}_i}(f_{i})\otimes r_{\mathfrak{P}_i}(\pi^{(i)})$ and $\rho_{\mathfrak{P}_i}(f_{i+1})\otimes r_{\mathfrak{P}_i}(\pi^{(i+1)})$; we need to check that the conditions of one of the ALT hold in this situation. These conditions involve regularity and irreducibility of the systems, local conditions at certain primes of $F$ and conditions on the shape of the residual image at $\mathfrak{P}_i$. In this section we will check carefully that they are satisfied at all steps.

If $f_i$ is congruent to a Galois conjugate of $f_{i+1}$, we have to take extra care. Whenever we reach a step $i$ in the chain where $f_i$ is congruent to the Galois conjugate of $f_{i+1}$ by $\sigma$, we need to modify the $2\times 2n$-dimensional chain as follows: For all $j>i$, we replace the $2n$-dimensional component $\pi^{(j)}$ by its Galois conjugate $\pi^{(j), \sigma^{-1}}$. In this way, at the $i$-th step we will need to propagate automorphy from $\rho_{\mathfrak{P}_i}(f^{i+1})\otimes r_{\mathfrak{P}_i}(\pi^{(i+1), \sigma^{-1}})$ to $\rho_{\mathfrak{P}_i}(f_{i})\otimes r_{\mathfrak{P}_i}(\pi^{(i)})$. Now we will assume that the representation $\rho_{\mathfrak{P}_i}(f^{i+1})\otimes r_{\mathfrak{P}_i}(\pi^{(i+1), \sigma^{-1}})$ is automorphic, and since automorphy is preserved by Galois conjugation, we will also have that $\rho_{\mathfrak{P}_i}(f^{\sigma}_{i+1})\otimes r_{\mathfrak{P}_i}(\pi^{(i+1)})$ is automorphic. We can use the congruence modulo $\mathfrak{P}_i$ to conclude that $\rho_{\mathfrak{P}_i}(f_{i})\otimes r_{\mathfrak{P}_i}(\pi^{(i)})$ is automorphic, provided that we can apply an ALT.

We need to check that the Galois conjugacy that we applied to all the steps from the $i+1$-th to the last is harmless; more precisely, we need to ensure that for all $j>i$, all the necessary conditions to apply an ALT will hold at the step $j$, provided they are satisfied by the original representation $\rho_{\mathfrak{P}_j}(f_{j})\otimes r_{\mathfrak{P}_j}(\pi^{(j)})$. Indeed, the level and weights of $\pi^{(j), \sigma^{-1}}$ and $\pi^{(j)}$ coincide, and the local types at primes of $F$ are preserved. Concerning the residual image, we already observed in Remark \ref{rem:GoodMonomialPrime} that any automorphic representation with the same type as $\pi'''$ at $q$, such that the primes in its ramification set are smaller than $2B$, will have large residual image at all primes $p\leq 2B$. 

Finally, since twisting a representation by a finite order character preserves both automorphy and all the conditions to apply the ALT we presented in Section \ref{sec:ALT}, we will be able to propagate automorphy from  $\{(\rho(f_{i+1})\otimes r(\pi^{(i+1)}))_{p}\}$  to $\{(\rho(f_{i})\otimes r(\pi^{(i)}))_{p}\}$ in the case that $f_i$ is congruent to a twist of $f_{i+1}$.

 Thus, without loss of generality we may assume that, at all steps in the chain, $f_i$ and $f_{i+1}$ are congruent modulo a prime above $p_i$, and in the remainder of the section we will not indicate either the Galois conjugations nor the twists in the notation for the elements of the $2\times 2n$-dimensional chain. This convention is also followed in previous papers of one of the authors (cf.~\cite[Section 8]{Di-mfold}).

\end{rem}
 
The general strategy is to move the $2$-dimensional component in the safe chain that ends at the CM form $f_h$ leaving the other component unchanged, but this is only possible at steps where the congruence in the $2$-dimensional component is of Type A (cf.~Section \ref{sec:2-dim Chain}). The problem with the other cases (congruence of Type B or C) is that the available automorphy lifting theorems (e.g.~Theorems \ref{thm:mixto} and \ref{thm:pd}) do not apply to propagate automorphy of the tensor product (from the right to the left). The idea to cope with this problem is to produce a
refinement of the chain \ref{2ndimchain} and then to introduce several auxiliary steps at each of the conflicting spots, steps consisting in congruences where either the $2$-dimensional, the $2n$-dimensional, or both components will change,  in such a way that some ALT can be
applied in each of these intermediate steps. The congruences in these additional steps
will occur modulo primes belonging to the set $\mathcal{P}_1$ of primes $p$ such that some of the congruences in the skeleton chain occurs modulo $p$, or belonging to the set 
$\mathcal{T}_{\text{aux}}$ of auxiliary primes fixed in the third preliminary step (subsection \ref{subsection:PreliminaryStep3}).This will be explained in full detail in the diagrams appearing in this section.\\

Of course, in order to conclude automorphy of these tensor products, this should be known at the rightmost one, i.e., we
need to know that  $\{(\rho(f_h) \otimes r(\pi^{(\textrm{end})}))_{p}\}$ is automorphic. Since automorphy is preserved by solvable base change
(cf. \cite[Lemma 2.2.2]{BLGGT}),
 we may look at the restriction of the system to $G_{F(\sqrt{-3})}$. Since $\rho(f_h)$ is the induction of a Hecke character from $\mathbb{Q}(\sqrt{-3})$,
  automorphy follows from the automorphy of $\{r_{p}(\pi^{(\textrm{end})})\}$ and \cite[Lemmas 2.2.1 and 2.2 4]{BLGGT}. \\

Before starting the analysis of each congruence, let us make some general remarks on the residual image
of the representations that appear in the chain. For the three Preliminary Steps, it was already checked in the previous section that a suitable ALT can be applied, in particular, residual images were taken care of. So we can start with the automorphic form $\pi'''$ in the $2n$-dimensional component. 

\begin{rem}\label{rem:AdequateImage} Let $p$ be a prime belonging to the set $\mathcal{P}_1$
of primes such that some of the congruences in the skeleton chain occurs modulo $p$, or to the set 
$\mathcal{T}_{\text{aux}}$ of auxiliary primes fixed in the third preliminary step (subsection \ref{subsection:PreliminaryStep3}).

Concerning the residual image of the $2$-dimensional component,
we already discussed in Section \ref{sec:2-dim Chain} that, if $p\in\mathcal{P}_1$, 
$\overline{\rho}_{p, \iota}(f)\vert_{G_{\mathbb{Q}(\zeta_p)}}$ is adequate$+$ (cf.~Remark \ref{rem:image-chain}).
We need to examine the restriction of $\overline{\rho}_{p, \iota}(f)$ to $G_{F(\zeta_p)}$. Recall that, except at Step 15, we have that
$\overline{\rho}_{p, \iota}(f)(G_{\mathbb{Q}(\zeta_p)})$ contains $\SL_2(\mathbb{F}_{p^r})$, where the exponent $r$ is big enough to ensure that the group
is adequate$+$. Since $F^+/\mathbb{Q}$ is totally real, we have
that $\overline{\rho}_{p, \iota}(f)(G_{F^+})$ also contains $\SL_2(\mathbb{F}_{p^r})$, with the same $r$ as before (cf.~\cite{Di-BaseChange}). Finally, since
the extension $F(\zeta_p)/F^+$ is abelian, we obtain that the image $\overline{\rho}_{p, \iota}(f)(G_{F(\zeta_p)})$ also
contains $\SL_2(\mathbb{F}_{p^r})$ for the same exponent $r$. As a conclusion, the residual image of the restriction
of $\rho_{p, \iota}(f)$ to the absolute Galois group of $F(\zeta_p)$ is also adequate$+$. Note that, in particular, this applies at the
Step 14, where the congruence takes place modulo $3$ and the image of the residual representation contains
the subgroup $G=\SL_2(\mathbb{F}_{3^3})$, which is adequate $+$ by Lemma \ref{lem:SL2}.
In Step 15, we have a mod $13$ congruence with a modular form $g$ with CM by $\mathbb{Q}(\sqrt{-3})$.
Since the field $F$ was chosen in such a way  that $F/\mathbb{Q}$ is linearly disjoint from $\mathbb{Q}(\sqrt{-3}, \zeta_{13})$, we
can ensure that the image of the restriction of $\overline{\rho}_{13, \iota}(f)$ to $G_{F(\zeta_{13})}$ is still a dihedral group, hence
adequacy$+$ is preserved under this operation. 

If $p\in \mathcal{T}_{\text{aux}}$, in particular it holds that $p>B$, in particular $p$ is greater than all the exceptional primes of
all modular forms that appear in the extended $2$-dimensional chain $\widetilde{C}$ (such primes where collected in the set $\mathcal{P}_2$
in Section \ref{sec:PreliminarySteps}). Thus the image of $\overline{\rho}_{p, \iota}(f)$ contains $\SL_2(\mathbb{F}_{p})$ and $p>3$, and from this again we can deduce that the image 
$\overline{\rho}_{p, \iota}(f)(G_{F(\zeta_p)})$     also contains  $\SL_2(\mathbb{F}_{p})$, and this guarantees (since $p>5$)
that this image is adequate$+$ by Lemma \ref{lem:SL2}.

Concerning the $2n$-dimensional component, as we already observed, the Good
Monomial prime occurring in $\{\overline{r}_{p}(\pi''')\}$ ensures that, for all primes $p\leq 2B$  (in particular, for our prime $p$), the
residual image of $r_{p, \iota}(\pi''')$ satisfies the conditions of Proposition \ref{prop:adequate} and is thus adequate$+$. The restriction to $G_{F(\zeta_p)}$ remains adequate$+$, since this restriction clearly still satisfies the conditions of Proposition \ref{prop:adequate}. Finally, for all but the last step, knowing that the image in the $2$-dimensional component  contains $\SL_2(\mathbb{F}_{p^r})$ and is contained in $\GL_2(\mathbb{F}_{p^w})$, and the image in the $2n$-dimensional component ($2n>2$) contains a classical group $\mathrm{SL}_{2n}(\mathbb{F}_{p^k})$, $\mathrm{SU}_{2n}(\mathbb{F}_{p^k})$,  $\mathrm{Sp}_{2n}(\mathbb{F}_{p^k})$
or $\Omega^{\pm}_{2n}(\mathbb{F}_{p^k})$ and is contained in its normaliser in $\GL_{2n}(\mathbb{F}_{p^w})$, they can only intersect in the quotient corresponding to the determinant.
Remark \ref{rem:TensorProduct} ensures that the image of the tensor product $(\overline{\rho}_{p, \iota}(f) \otimes \overline{r}_{p, \iota}(\pi'''))\vert_{G_{F(\zeta_p)}}$
is adequate$+$. \\

In the last step, at the $2n$-dimensional component we will have an automorphic representation $\pi^{(\mathrm{end})}$, which may not coincide with $\pi'''$, but which will share all the relevant features with it (see Remark \ref{rem:piend}). In particular, its level will be coprime to $3$. Moreover, the good monomial prime for $\pi'''$ is still a good monomial prime for $\pi^{(\mathrm{end})}$. Therefore, if we denote by $\overline{r}_{13, \iota}:G_F\rightarrow \GL_2(\overline{\mathbb{F}}_{13})$  the Galois representation attached to $\pi^{(\mathrm{end})}$ at the prime $13$, the image of $\overline{r}_{13, \iota}\vert_{G_{F(\zeta_{13})}}$ contains a classical group of the form $\SL_{2n}(\mathbb{F}_{13^s})$, $\SU_{2n}(\mathbb{F}_{13^s})$, $\Sp_{2n}(\mathbb{F}_{13^s})$ or $\Omega^{\pm}_{2n}(\mathbb{F}_{13^s})$, 
and is contained in its normaliser in $\GL_{2n}(\mathbb{F}_{13^w})$ for some (big enough) exponents $r, s, w$.

At the $2$-dimensional component we have the modular form $g\in S_2(27)$ obtained at step 15 of the chain (cf.~\ref{Step:15}). Let $\overline{\rho}_{13, \iota}:G_{\mathbb{Q}}\rightarrow \GL_2(\overline{\mathbb{F}}_{13})$ be the Galois representation attached to $g$. As we explained above, the image of the restriction $\overline{\rho}_{13, \iota}\vert_{G_{F(\zeta_{13})}}$ is adequate$+$.

Consider the extension $M$ of $F(\zeta_{13})$ such that the image of $\overline{r}_{13, \iota}$ is a classical group of the form $\SL_{2n}(\mathbb{F}_{13^s})$, $\SU_{2n}(\mathbb{F}_{13^s})$, $\Sp_{2n}(\mathbb{F}_{13^s})$ or $\Omega^{\pm}_{2n}(\mathbb{F}_{13^s})$.  Since the level of $\pi^{(\mathrm{end})}$ is coprime to $3$ and $F$ is unramified at $3$, $M$ is linearly disjoint from $F(\zeta_{13}, \sqrt{-3})$ over $F(\zeta_{13})$.  Therefore, the image of the restriction of $\overline{\rho}_{13, \iota}$ to $G_M$ is still a dihedral group; by Lemma \ref{lem:ForDihedral} it is adequate+. We can apply Proposition \ref{prop:DirectProduct} to conclude that the image of the tensor product of $(\overline{\rho}_{13, \iota}\vert_{G_M})\otimes (\overline{r}_{13, \iota}\vert_{G_M})$ is adequate+. Note that  Proposition \ref{prop:adequate} ensures that $[M: F(\zeta_{13})]$ is coprime to $13$; thus we can apply Lemma \ref{lem:PrimeToP} to conclude that the image of $(\overline{\rho}_{13, \iota}\otimes \overline{r}_{13, \iota})\vert_{G_{F(\zeta_{13})}}$ is adequate+.

\end{rem}

This remark is extremely important because it ensures that in what follows in each congruence the residual image of the tensor product is large enough for the required ALT to be applied (from right to left).

\subsection{Standard step: type A congruence}\label{sub:type_A}

To simplify notation, we call $\pi=\pi'''$  the automorphic form obtained after the three preliminary steps (cf.~Section \ref{sec:PreliminarySteps}). \\
To fix ideas, we will first describe how can we propagate automorphy from right to left when in the
two-dimensional skeleton chain we have a type A congruence. The situation is the following: we have
a congruence between two modular forms, say $f$ and $f'$ modulo some $p\in \mathcal{P}_1$, and we know that $\{(\rho(f')\otimes r(\pi))_{p}\}$
is automorphic. Thus we have that the residual representation $\overline{\rho}_{p,\iota}(f')\otimes \overline{r}_{p,\iota}(\pi)\equiv \overline{\rho}_{p, \iota}(f)\otimes \overline{r}_{p, \iota}(\pi)$
is automorphic. We want to apply Theorem (ALT-HARRIS) with $R=\overline{\rho}_{p, \iota}(f)$, $R'=\overline{\rho}_{p, \iota}(f')$, 
$S=S'=\overline{r}_{p, \iota}(\pi)$ to prove that the compatible system $\{(\rho(f)\otimes r(\pi))_{p}\}$ is automorphic.
Let us check that all hypotheses in this theorem are satisfied.

First of all, note that all primes in $\mathcal{P}_1$ (primes where the congruences of the chain take place); in particular, our prime $p$, are odd.
Furthermore, by the choice of $F$, we know that $\zeta_{p}\not\in F$. We have $n_0=2$, $n_1=2n$, therefore the condition $p\nmid n_0$ is satisfied.
Moreover:

\begin{enumerate}

\item The compatible system $\{(\rho(f')\otimes r(\pi))_{p}\}$ is regular. Indeed, recall that the Hodge-Tate weights of $\pi$ with respect to any embedding $\tau\colon F \hookrightarrow \overline{\Q}_p$ are
\begin{equation*}
 \mathrm{HT}_\tau(\pi)=\{0,a,\ldots,(2n-1)a\}
 \end{equation*}
 where $a>C_0$ (see the list before Remark \ref{rem:From_Pi2_To_Pi3}). In particular, they are independent of the choice of $\tau$, and $C_0$-very spread, so we don't lose regularity after tensoring with $\{\rho_{p}(f')\}$, whose Hodge-Tate weights are smaller than or equal to $C_0$ (see Remark \ref{rem:Regularity_tensor}). Moreover, as explained in Remark \ref{rem:Constant2}, the extra regularity condition required in this ALT is satisfied. The same applies to $\{(\rho(f)\otimes r(\pi))_{p}\}$.

\item The assumption of a type A congruence  implies that $\rho_{p, \iota}(f)$ and $\rho_{p, \iota}(f')$ are potentially diagonalisable at all
places $v\vert p$.

\item Furthermore, for each finite places $v\nmid p$ of $F$, the representations $\rho_{p, \iota}(f)\vert_{G_v}$ and $\rho_{p, \iota}(f')\vert_{G_v}$ are connected.

\item  Concerning the residual image,
we already discussed in Remark \ref{rem:AdequateImage} above that the  image of  $(\overline{\rho}_{p, \iota}(f)\otimes  \overline{r}_{p, \iota}(\pi))\vert_{G_{F(\zeta_p)}}$
is adequate$+$.

\end{enumerate}

Thus all conditions in Theorem (ALT-HARRIS) hold, and we can conclude that the compatible system $\{(\rho(f)\otimes r(\pi))_{p}\}$ is automorphic.

In the moves where the 2-dimensional skeleton chain presents a congruence of types B or C between $f$ and $f'$, the previous reasoning does not work. We will
introduce some new steps $$\{(\rho(f)\otimes r(\pi))_{p}\}\sim \{(\rho(f^{(1)})\otimes r(\pi^{(1)}))_{p}\}\sim \cdots \sim \{(\rho(f^{(r)})\otimes r(\pi^{(r)}))_{p}\}\sim \{(\rho(f')\otimes r(\pi^{(r)}))_{p}\}$$ 
 by performing level and weight changes in the $2n$-dimensional component, the $2$-dimensional component, or both. All these changes (except at the point where the $2$-dimensional representation corresponds to Step 12, Section 4.15, which will be explained in detail in Section 6.4) will be reversed, meaning that
after these exceptional moves, we return to the compatible system $\{r(\pi)_{p}\}$ in the $2n$-dimensional component (that is to say, $\pi^{(r)}=\pi$).

To perform the level and weight changes in the $2n$-dimensional component, two essential tools will be used. In the ordinary case, we will make use of Hida families to change the weight (cf.~\cite{Ger10}).
Moreover, in the potentially diagonalisable situation, we will make use of Theorem \ref{thm:Weight0}.

\subsection{Diagram 1-Congruence of type C}

\noindent \textbf{Set-up:}
We consider a move in the $2$-dimensional skeleton chain involving a type C congruence modulo a prime $p$ between two modular forms $f_1, f_2$.
We want to construct a $2\times 2n$ ``safe chain'' connecting $\{(\rho(f_1)\otimes r(\pi))_{p}\}$ with $\{(\rho(f_2)\otimes r(\pi))_{p}\}$.
Observe that we cannot apply (ALT-HARRIS), since the two compatible systems are not connected at all places (above a prime different from $p$) because of the change of level in the $2$-dimensional
component.

We will distinguish two situations: if the next link in the $2$-dimensional skeleton chain is a type B congruence, we skip this step and proceed
directly to Diagram 2 (cf.~Subsection \ref{sub:type_B}). Otherwise, we proceed as described below.

Our aim in this situation is to use Hida families to connect $\{(\rho(f_1)\otimes r(\pi))_{p}\}$ with $\{(\rho(f_2)\otimes r(\pi))_{p}\}$. In order
to do this, we first need to connect $r(\pi)$ with an ordinary automorphic representation. There is a way to do this, consisting of the following steps: first, we apply Theorem \ref{thm:Weight0} to
change $\pi$ so that its weights with respect to any embedding $\tau$ are $\{0, 1, \dots, 2n-1\}$; then, we consider it modulo a Steinberg prime. By \cite{Ger10},
we will obtain an ordinary $2n$-dimensional representation.

However, if we perform the movements described above in the $2n$-dimensional component $\{r_{p}(\pi)\}$
while leaving the $2$-dimensional component $\{\rho_{p}(f_1)\}$
unchanged, we may lose regularity (recall that regularity was guaranteed because the weights of $\{r_{p}(\pi)\}$ were $C_0$-very spread, see Remark \ref{rem:Regularity_tensor}). The challenge is to
perform safe moves in \emph{both} components at the same time.

\begin{rem}\label{rem:OrdinarityOfOtimes} Note that, in general, the tensor product of two ordinary automorphic representations is not necessarily ordinary. However, in our situation we can ensure that we obtain an ordinary automorphic representation, because we arranged the weights so that they are sufficiently spread (cf. to Lemma \ref{lemma:regularity}). Indeed, we will have that either the Hodge-Tate weights
of the $2$-dimensional components are separated by a greater distance than the difference of any two of the Hodge-Tate weights of the $2n$-dimensional component (which will have Hodge type 0), or else the weights of the $2n$-dimensional components are $C_0$-very spread. In particular, the difference
between any two weights is greater than the weight of the modular form occurring in the $2$-dimensional component (a similar argument is applied inside the proof of Theorem 2.1 in \cite{Di-mfold}).
\end{rem}

By construction, the Hodge-Tate weights of $\rho_{p, \iota}(f_1)$ and $r_{p, \iota}(\pi)$ with respect to any choice of embedding $\tau$ are, respectively, $\{0, j_1\}$ and $\{0, a, \dots, (2n-1)a\}$ (where $j_1+1$ is the weight of the modular form $f_1$, and $a>C_0$). We will introduce three additional links between $\{(\rho(f_1)\otimes r(\pi))_{p}\}$ and $\{(\rho(f_2)\otimes r(\pi))_{p}\}$.
We illustrate the new chain via the following diagram. 

At the left column we write the weights of the modular forms involved, and at the right column
we write the weights of the $2n$-dimensional automorphic representations. Because these Hodge-Tate weights are independent of the choice of embedding $\tau$, we omit it from the notation. At the left hand side of each arrow, we indicate the prime modulo which the congruence occurs, and at the right hand side we specify which ALT will enable us to propagate modularity from bottom to top.

$$\xymatrix{[\{0, j_1\}, \{0, a, \dots, (2n-1)a\}] \ar[d]_{\pmod{p}}^{(ALT-HARRIS)}\\
            [\{0, j'_1\}, \{0, a, \dots, (2n-1)a\}] \ar[d]_{\mod{t_{\mathrm{aux}}}}^{(ALT-PD)}\\
            [\{0, j'_1\}, \{0, 1, 2, \dots, 2n-1\}] \ar[d]_{\pmod{p}}^{(ALT-MIXED)} \\
            [\{0, j_1\},  \{0, m_1, \dots, m_{2n-1}\}] \ar[d]_{\pmod{p}}^{(ALT-MIXED)} \\
            [\{0, j'_1\}, \{0, 1, 2, \dots, 2n-1\}]\ar[d]_{\pmod{t_{\mathrm{aux}}}}^{(ALT-PD)}\\
            [\{0, j'_1\}, \{0, a, \dots, (2n-1)a\}] \ar[d]_{\pmod{p}}^{(ALT-HARRIS)}  \\
            [\{0, j_1\}, \{0, a, \dots, (2n-1)a\}] }$$

 We explain in detail each of the steps in the diagram:

 \begin{enumerate}
 \item In the first congruence, we want to raise the weight $j_1$ to $j'_1\in[2n, 2n + p^2-1]$. If our initial weight $j_1$ already belongs
  to this interval, we skip this congruence. Otherwise, Remark \ref{rem:k_distinto_de_p} enables us to
  apply Theorem \ref{thm:BigWeight} to the modular
  form $f_1$ to produce another modular form $f'$ of weight $j'_1$, such that $\rho_{p, \iota}(f')$ is potentially diagonalisable and connects to $\rho_{p, \iota}(f_1)$ at
  all places $v\nmid p$.
  The new link of the chain will be $\{(\rho(f_1)\otimes r(\pi))_{p}\}\sim\{ (\rho(f')\otimes r(\pi))_{p}\}$. Note that this link is ``safe'', since
  from the automorphy of the right-hand side we can conclude automorphy of the left-hand-side by using (ALT-HARRIS). Indeed, since $C_0$ was constructed in
  such a way that $2n + p^2\leq C_0$ for all $p\in \mathcal{L}$, the tensor product $\{(\rho(f')\otimes r(\pi))_{p}\}$ is regular (see Remark \ref{rem:Regularity_tensor}). The rest of the
  hypothesis required to apply (ALT-HARRIS) hold exactly as in the standard step (cf.~Subsection \ref{sub:type_A}).

  \item In the second congruence we perform the weight change in the $2n$-dimensional component to achieve consecutive weights. First of all, we choose an
  auxiliary prime $t_{\mathrm{aux}}\in \mathcal{T}_{\mathrm{aux}}$. We will make a new link in the chain via a congruence modulo this auxiliary prime. Note
  that, by definition of the constant $B$, we have that $t_{\mathrm{aux}}<2B$, thus the residual image of $r_{t_{\mathrm{aux}}, \iota}(\pi)$ is adequate$+$.
  Next, we use Theorem \ref{thm:Weight0} to produce a Hodge type $0$ automorphic lift of  $\overline{r}_{t_{\mathrm{aux}}, \iota}(\pi)$, say
  $r_{t_{\mathrm{aux}}, \iota}(\pi^{(1)})$. Note that, by our choice
  of the set $\mathcal{T}_{\mathrm{aux}}$, the prime $t_{\mathrm{aux}}$ is bigger than $2n$, unramified in $F$ and all primes of $F^+$ above it are totally split in $F/F^{+}$.

  The new link of the chain will be $\{(\rho(f')\otimes r(\pi))_{p}\}\sim \{(\rho(f')\otimes r(\pi^{(1)}))_{p}\}$. Note that (ALT-PD) ensures that this link
  is ``safe''. Indeed, the regularity of the right-hand-side holds because $j'_1>2n + p^2-1$ and we can apply Lemma \ref{lemma:regularity}. Besides, the $2$-dimensional component is Fontaine-Laffaille at $t_{\mathrm{aux}}$,
  and the $2n$-dimensional component is potentially diagonalisable.

  \item In the third congruence, we perform the level raising in the $2$-dimensional component,
  moving to $f_2$. Since the weight of $f_2$ is $j_1$, we have to modify the weights of the $2n$-dimensional component to preserve regularity.
  Recall that our prime $p$ belonged to $\mathcal{L}$. Therefore it is a Steinberg prime in the level of $\pi$. This property is preserved in $\pi^{(1)}$, since in Theorem \ref{thm:Weight0} we did not change the type at any other prime different from $t_{\mathrm{aux}}$. Since $\{r_{p}(\pi^{(1)})\}$ has Hodge type $0$, we have that $r_{p, \iota}(\pi^{(1)})$ is ordinary.
  We can thus move it in the corresponding Hida family and obtain a lift  $r_{p, \iota}(\pi^{(2)})$, of $\overline{r}_{p, \iota}(\pi^{(1)})$,
  with weights $\{0, m_1, \dots, m_{2n-1}\}$ with respect to any embedding $\tau$. The only condition we need to ask to these weights is that they are sufficiently
  spread so that the tensor product $\{(\rho(f_2)\otimes r(\pi''))_{p}\}$ is regular.

  We will then have a congruence modulo $p$ between $\rho_{p, \iota}(f')\otimes r_{p, \iota}(\pi^{(1)})$ and $\rho_{p, \iota}(f')\otimes r_{p, \iota}(\pi^{(2)})$. In fact, one can
  check that the residual representations $\overline{\rho}_{p, \iota}(f_1)$, $\overline{\rho}_{p, \iota}(f')$ and $\overline{\rho}_{p, \iota}(f_2)$ are all equivalent (the first two by definition of $f'$, 
  and the first and third because they are linked by a mod $p$ congruence in the $2$-dimensional skeleton chain). Thus we may freely replace $\rho_{p, \iota}(f')$ by $\rho_{p, \iota}(f_2)$, and we obtain a congruence mod $p$ between $\rho_{p, \iota}(f')\otimes r_{p, \iota}(\pi^{(1)})$ and
  $\rho_{p, \iota}(f_2)\otimes r_{p, \iota}(\pi^{(2)})$. Note that this establishes a ``safe'' link: Theorem (ALT-MIXED) applies in this situation, ensuring that
  if the right-hand-side is automorphic, so is the left-hand-side too.

  Note that, at this point, we managed to replace $\{\rho_{p}(f_1)\}$ by $\{\rho_{p}(f_2)\}$. However,
  the $2n$-dimensional component is not the original one. The aim in the next three congruences is to undo the changes that occurred in this component.

 \item In the fourth congruence, we undo the second move in the $2n$-dimensional component. Like in the third congruence, in order
 to replace $\pi^{(2)}$ by $\pi^{(1)}$ and preserve regularity, we need to change the weight in the first component, choosing (for instance)
 the same $j'_1\in [2n, 2n + p^2-1]$. To produce a modular form $f''$ which has this weight and such that $\rho_{p, \iota}(f'')$ and $\rho_{p, \iota}(f_2)$ are
 connected locally at all places $v\nmid p$, we apply again Theorem \ref{thm:BigWeight}.
Thus the next link is a congruence modulo $p$ with
$\{(\rho(f'')\otimes r(\pi^{(1)}))_{p}\}$.
 Just like in the third congruence, Theorem (ALT-MIXED) ensures that this new link is ``safe''.

 \item In the fifth congruence, we undo the first move in the $2n$-dimensional component, linking our compatible system with
 $\{(\rho(f'')\otimes r(\pi))_{p}\}$ modulo $t_{\mathrm{aux}}$. Like in the second step of this diagram, the link is ``safe'' because of Theorem (ALT-PD).

 \item Finally, we link our compatible system to $\{(\rho(f_2)\otimes r(\pi))_{p}\}$ working modulo $p$. Theorem (ALT-HARRIS) ensures that this last congruence is ``safe''.

 \end{enumerate}

\subsection{Diagram 2: Congruence of type B}\label{sub:type_B}

 We turn now to the case when the move in the $2$-dimensional skeleton chain involves a congruence of type B modulo a prime $p$ between two
 modular forms $f, f'$. As in the previous case, we want to construct a $2\times 2n$ ``safe chain'' connecting
 $\{(\rho(f)\otimes r(\pi))_{p}\}$ with $\{(\rho(f')\otimes r(\pi))_{p}\}$.
Note that in this situation we cannot apply (ALT-HARRIS) because  the $2$-dimensional component of the left-hand side is not potentially
diagonalisable. Moreover, we cannot apply (ALT-MIXED) because in this situation the components which are ordinary are the $2$-dimensional,
and not the higher-dimensional ones.

As noted in Remark \ref{rem:CongruenceTypes}, whenever we have a type B congruence between $f$ and $f'$ (modulo $p_{i+1}$), the preceding congruence is of type C,
 and $f$ is weight $2$ and Steinberg at $p_{i+1}$ (with one exception that we will consider in the next subsection). The diagram of congruences to connect  $\{(\rho(f)\otimes r(\pi))_{p}\}$ with $\{(\rho(f')\otimes r(\pi))_{p}\}$
will rest on the previous congruence. Thus, if we meet a congruence of type C followed by a congruence of type B, instead of executing
Diagram 1, we perform the following diagram, (which includes an iteration of Diagram 1 inside). The goal of these new links is to modify
the $2n$-dimensional component to make it ordinary at $p_{i+1}$, just before applying Diagram 1. Thus once we finish performing Diagram 1, we
will be in a position where both components are ordinary, and we will be able to move using a Hida family. We explain in detail each of the steps.

\begin{enumerate}
 \item In the first step we have a modular form $f_i$ of weight $j_1$, potentially Barsotti-Tate at 
 $p_{i+1}$ (recall that this is before the previous congruence in the chain, a type C congruence modulo $p_i$, is executed). 
 What we do is a congruence modulo $p_{i+1}$. If $j_1\in[2n, 2n + p_{i+1}^2-1]$, then we can skip
 this step. Otherwise, by Theorem \ref{thm:BigWeight}, there exists
 another modular form $f'$, of weight $j'_1$ in this interval, satisfying that $\rho_{p_{i+1}, \iota}(f')$ 
 is potentially diagonalisable and connects to $\rho_{p_{i+1}, \iota}(f_i)$ at
  all primes $v\nmid p_{i+1}$. Thus the first link is a congruence modulo $p_{i+1}$, 
  between $\{(\rho(f_i)\otimes r(\pi))_{p}\}$ and $\{(\rho(f')\otimes r(\pi))_{p}\}$.
  As in the first step of Diagram 1, this new link is ``safe'', that is, automorphy propagates from right to left.

 \item We choose a prime $t_{\mathrm{aux}}\in \mathcal{T}_{\mathrm{aux}}$ which has not been used before, and using Theorem \ref{thm:Weight0}, we fix an automorphic
 representation $\pi^{(1)}$ of Hodge type 0, potentially diagonalisable at $t_{\mathrm{aux}}$, which is an automorphic lift of $\overline{r}_{t_{\mathrm{aux}}, \iota}(\pi)$.
  So at this step we have a congruence modulo $t_{\mathrm{aux}}$ between $\{(\rho(f')\otimes r(\pi))_{p}\}$ and $\{(\rho(f')\otimes r(\pi^{(1)}))_{p}\}$. Theorem (ALT-PD) guarantees that
 this new link is ``safe''. This step ensures us that the second component will be ordinary at $p_{i+1}$.

  \item In this step we want to replace $f'$ by our original modular form $f_i$, while obtaining an ordinary automorphic representation in the second
 component. We can do this by looking modulo $p_{i+1}$, since $\rho_{p_{i+1}, \iota}(f_i)$ is congruent to $\rho_{p_{i+1}, \iota}(f')$ modulo $p_{i+1}$, and $r_{p_{i+1}, \iota}(\pi')$ is ordinary.
 However, we do not want to lose regularity, hence we move in the corresponding Hida family and choose a member $r_{p_{i+1}, \iota}(\pi^{(2)})$ with weights $\{0, \tilde{m}_1, \dots, \tilde{m}_{2n-1}\}$ with respect to any embedding $\tau$.
 in such a way that they are $C_0$-very spread. Thus, we consider the link $\{(\rho(f')\otimes r(\pi^{(1)}))_{p}\}\sim \{(\rho(f_i)\otimes r(\pi^{(2)}))_{p}\}$.
 The ``safety'' of this link is guaranteed exactly as in Step 3 of Diagram 1.

$$\xymatrix{[\{0, j_1\},\{0, a, \dots, (2n-1)a\}] \ar[d]_{\pmod{p_{i+1}}}^{(ALT-HARRIS)}\\
            [\{0, j'_1\}, \{0, a, \dots, (2n-1)a\}] \ar[d]_{\mod{t_{\mathrm{aux}}}}^{(ALT-PD)}\\
            [\{0, j'_1\}, \{0, 1, 2, \dots, 2n-1\}] \ar[d]_{\pmod{p_{i+1}}}^{(ALT-MIXED)} \\
            [\{0, j_1\},  \{0, \tilde{m}_1, \dots, \tilde{m}_{2n-1}\}] \ar@{.>}[dd]!U|*+<15pt>[o][F-]{\text{DIAGRAM 1}} \\ 
            \ \\
            [\{0, 1\},  \{0, \tilde{m}_1, \dots, \tilde{m}_{2n-1}\}] \ar[d]_{\pmod{p_{i+1}}}^{(ALT-ORD)} \\
            [\{0, j_2\}, \{0, \tilde{m}_1, \dots, \tilde{m}_{2n-1}\}]\ar[d]_{\mod{p_{i+1}}}^{(ALT-MIXED)}\\
            [\{0, j'_2\}, \{0, 1, \dots, 2n-1\}] \ar[d]_{\mod{t_{\mathrm{aux}}}}^{(ALT-PD)}  \\
            [\{0, j'_2\},\{0, a, \dots, (2n-1)a\}] \ar[d]_{\mod{p_{i+1}}}^{(ALT-HARRIS)}  \\
            [\{0, j_2\}, \{0, a, \dots, (2n-1)a\}]  }$$

 \item The next step is an application of Diagram 1 to the modular form $f_i$ and the automorphic representation $\pi^{(2)}$. This enables us to
 perform the type C congruence modulo $p_i$ in the $2$-dimensional skeleton chain, thus linking $\{(\rho(f_i)\otimes r(\pi^{(2)}))_{p}\}$ with $\{(\rho(f_{i+1})\otimes r(\pi^{(2)}))_{p}\}$.

\item Now we will devote our attention to the next step in the $2$-dimensional skeleton chain, namely the type B congruence. Thanks to the preceeding steps,
 the situation is the following: on the $2$-dimensional component we have a weight-$2$ modular form $f_{i+1}$, which is both potentially diagonalisable at $p_i$ and Steinberg at $p_{i+1}$. On the
 $2n$-dimensional component we have an automorphic representation $\pi^{(2)}$, such that  $r_{p_{i+1}, \iota}(\pi^{(2)})$ is ordinary.
 In this step we link $\{(\rho(f_{i+1})\otimes r(\pi^{(2)}))_{p}\}$ to $\{(\rho(f_{i+2})\otimes r(\pi^{(2)}))_{p}\}$ via  a congruence modulo $p_{i+1}$. The ``safety'' of this
 link is guaranteed by (ALT-ORD). Note that, by Remark \ref{rem:OrdinarityOfOtimes}, both sides of the congruence are ordinary at $p_{i+1}$.

 \item The purpose of the remaining steps is to undo the changes we did in the $2n$-component, in order to obtain back the compatible system $\{r_{p}(\pi)\}$. Thus,
 we need to invert the moves we made in the first three steps of this diagram, without reverting the moves in the $2$-dimensional component. In order to do this, we
 will need to pass to Hodge type $0$ in the $2n$-dimensional component (using Theorem \ref{thm:Weight0}). Therefore, we need to modify the weights of $f_{i+2}$ to preserve regularity. Thus, applying Theorem \ref{thm:BigWeight}, we  produce
 a modular form $f'_{i+2}$, such that $\rho_{p_{i+1}, \iota}(f'_{i+2})$ is congruent to $\rho_{p_{i+1}, \iota}(f_{i+2})$, and satisfying that its weight belongs to $[2n, 2n+p_{i+1}^2-1]$.
 (if the weight of $f_{i+2}$ belongs already to this interval, we skip this step). 
Moreover, note that by definition of $\pi^{(2)}$, it holds that $r_{p_{i+1}, \iota}(\pi^{(2)})$ is congruent to $r_{p_{i+1}, \iota}(\pi^{(1)})$ modulo $p_{i+1}$. 
Thus we may link $\{(\rho(f_{i+2})\otimes r(\pi^{(2)}))_{p}\}$
 to $\{(\rho(f'_{i+2})\otimes r(\pi^{(1)}))_{p}\}$ via a congruence modulo $p_{i+1}$.
  The ``safety'' of this new link is ensured by (ALT-MIXED). 

 \item Next, we pick the same prime $t_{\mathrm{aux}}$ that we used in step 2 of this diagram. Note that $r_{t_{\mathrm{aux}}, \iota}(\pi^{(1)})$ is congruent to
 $r_{t_{\mathrm{aux}}, \iota}(\pi)$. Thus we may link $\{(\rho({f'}_{i+2})\otimes r(\pi^{(1)}))_{p}\}$ with $\{(\rho(f'_{i+2})\otimes r(\pi))_{p}\}$ via
 a congruence modulo $t_{\mathrm{aux}}$. This link is ``safe'' because Theorem (ALT-PD) applies, exactly as in step 2 of Diagram 1.

 \item Finally, we link $\{(\rho(f'_{i+2})\otimes r(\pi))_{p}\}$ to $\{(\rho(f_{i+2})\otimes r(\pi))_{p}\}$ through a congruence modulo $p_{i+1}$. As in the first step of Diagram 1,
 one can check that this link is ``safe'', thus completing the diagram.

\end{enumerate}

\subsection{Diagram 3}\label{sec:Diagram3}

In this subsection we address the last type-B congruence of the $2$-dimensional skeleton chain, the one appearing in Step 12 (see Section 4.15). The essential peculiarity of this congruence
is that, unlike all other congruences of type B, it is not preceded by a type C congruence. Therefore, we cannot apply Diagram 2 directly
to it, and we are forced to devise a special diagram for this one step.

Recall (cf.~Section \ref{sec:2-dim Chain}) that, as input for this link, we have a modular form $f=g_5\in S_{44}(\Gamma_1(17))$, and through a congruence
modulo $p=43$, it is linked to a modular form $g_4\in S_2(\Gamma_1(17\cdot 43))$. The modular form $f$ is crystalline and ordinary at $p=43$, whereas $g_4$
is ordinary (Steinberg) at $43$. 
We want to pass to a Hodge type 0 lift of $\overline{r}_{43, \iota}(\pi)$ corresponding to an ordinary automorphic representation $\pi$ which is
ordinary at $43$. As before, we have to be careful with the $2$-dimensional component in order not to lose regularity. 
We perform the usual tricks.

\begin{enumerate}

\item In the first step, we modify the weight of $f$ so that it is
contained in the interval $[2n, 2n + 43^2-1]$. If $n\leq 22$, we may skip this step.
In order to perform this weight change, we first apply Theorem \ref{thm:BigWeight}
to produce a modular form $f'$ such that $\rho_{43, \iota}(f)$ and $\rho_{43, \iota}(f')$
are congruent modulo $43$. We may thus link $\{(\rho(f)\otimes r(\pi))_{p}\}$ with
$\{(\rho(f')\otimes r(\pi))_{p}\}$. As in the first step of Diagram 1, (ALT-HARRIS) ensures
that this link is ``safe''.

\item Next we choose a prime $t_{\mathrm{aux}}\in \mathcal{T}_{\mathrm{aux}}$ which has not yet been used. 
We may then link $\{(\rho(f')\otimes r(\pi))_{p}\}$ with
$\{(\rho(f')\otimes r(\pi^{(1)}))_{p}\}$, where $\pi^{(1)}$ is an automorphic representation of Hodge type 0, 
whose existence is guaranteed by Theorem \ref{thm:Weight0}.
As in the second step of Diagram 2, this link is ``safe'' because (ALT-PD) applies.

\item The last step consists of a congruence modulo $43$, where we link
$\{(\rho(f')\otimes r(\pi^{(1)}))_{p}\}$ with $\{(\rho(g_4)\otimes r(\pi^{(2)}))_{p}\}$, where $\pi^{(2)}$ is a member 
in the Hida family at $43$ containing  $\pi^{(1)}$ with $C_0$-very spread weights $\{0, \widehat{k}, 2\widehat{k}, \ldots, (2n-1)\widehat{k}\}$ with respect to any embedding $\tau$ (for 
the number $\widehat{k}$ fixed at the beginning of Section \ref{sec:PreliminarySteps}).  In this case,
Theorem (ALT-ORD) ensures that this link is ``safe''.

\end{enumerate}

 $$\xymatrix{[\{0, 44\}, \{0, a, \dots, (2n-1)a\}] \ar[d]_{\pmod{43}}^{(ALT-HARRIS)}\\
            [\{0, j'_1\}, \{0, a, \dots, (2n-1)a\} \ar[d]_{\mod{t_{\mathrm{aux}}}}^{(ALT-PD)}\\
            [\{0, j'_1\}, \{0, 1, 2, \dots, 2n-1\}] \ar[d]_{\pmod{43}}^{(ALT-ORD)} \\
            [\{0, 1\},  \{0, \widehat{k}, 2\widehat{k}, \dots, (2n-1)\widehat{k}\}]  \\}$$

 \begin{rem}\label{rem:piend}
 Note that, at this step, we cannot undo the changes we did in the $2n$-dimensional component without
   reversing the changes in the $2$-dimensional component. Hence, after performing this step, we will work
   with a different automorphic representation $\pi^{(2)}$ instead of $\pi$. However, this will not affect
   our reasoning; in particular, automorphy will still be propagated from right to left. Let us expand on this point.
   \begin{enumerate}
   \item The weights of $\pi^{(2)}$ are $\{0, \widehat{k}, \dots,  (2n-1)\widehat{k}\}$. Since $\widehat{k}>2C_0$, they are $C_0$-very spread.

   \item The level of $\pi^{(2)}$ differs from the level of $\pi$ at most at the primes $43$ and $t_{\mathrm{aux}}\in \mathcal{T}_{\mathrm{aux}}$.
   Since the inequality $t_{\mathrm{aux}}<2B$ is satisfied,  the good monomial prime for $\pi$ is also a good monomial prime
   for $\pi^{(2)}$ (cf. Remark \ref{rem:GoodMonomialPrime}). However, the prime $43$ is no longer a prime of Steinberg type for $\pi^{(2)}$. Fortunately, this is of
   no consequence for us, since in the three remaining congruences of the skeleton chain, we do not use the prime $43$ at all (more precisely,
   we perform congruences modulo $17, 3$ and finally mod $13$).

   \end{enumerate}

 \end{rem}

 \subsection{Conclusion} After performing the three Preliminary Steps described in Section 5, and after following the algorithm described in this section which allows us to move along the safe chain of 2-dimensional representations in a way that ensures that, for the tensor product, a suitable ALT can be applied from right to left, we conclude that the automorphy of the tensor product containing the CM form $f_h$ (which follows from solvable base change, as explained at the beginning of this section) implies automorphy of the original tensor product. A priori, this is only deduced over a suitable solvable number field $F$, the CM field  where the last automorphic form appearing in the chain is defined; but a final application of solvable base change allows us to conclude automorphy over $\mathbb{Q}$. This proves Theorem \ref{thm:main}.

\bibliography{Bibliog}
\bibliographystyle{alpha}

\end{document}